\documentclass[a4paper]{article}
\usepackage[english]{babel}
\usepackage[utf8]{inputenc}
\usepackage{amsmath,amsthm}
\usepackage{amsfonts,amssymb}

\usepackage{mathtools}

\usepackage{tikz,tikz-cd} 
 \usetikzlibrary{decorations.pathmorphing} 

\newcommand{\mytop}[2]{\genfrac{}{}{0pt}{}{#1}{#2}}

\usepackage{xifthen} 

\newcommand{\C}{\mathbb{C}}

\newcommand{\D}{\mathbb{D}}
\newcommand{\G}{\mathbb{G}}

\newcommand{\Pb}{\mathbb{P}}

\newcommand{\Z}{\mathbb{Z}}

\newcommand{\Ocal}{\mathcal{O}}

\newcommand{\Dcal}{\mathcal{D}}

\newcommand{\Ccal}{\mathcal{C}}

\usepackage{titlesec}
\titleformat{\subsubsection}[runin]{\normalfont}{\thesubsubsection}{0pt}{}[.]

\renewcommand{\thesubsubsection}{\arabic{section}.\arabic{subsection}.\arabic{subsubsection}}
\csname @addtoreset\endcsname{subsubsection}{section}

\newcommand{\block}[1]
{

\par \subsubsection{} #1

\bigskip}

\newcommand{\blockn}[1]{\par #1 \bigskip}

\newcommand{\Th}[1]
	{
	\bigskip	
	\textbf{Theorem : }{\itshape #1}
		
	\bigskip
	}

\newcommand{\Conj}[1]
	{
	\bigskip	
	\textbf{Conjecture : }{\itshape #1}
		
	\bigskip
	}

\newcommand{\Prop}[1]
	{

	\bigskip
	
	\textbf{Proposition : }{\itshape #1}
		
	\bigskip
	
	}

\newcommand{\Cor}[1]
	{

	\bigskip
	
	\textbf{Corollary : }{\itshape #1}	
		
	\bigskip

	}

\newcommand{\Lem}[1]
	{

	\bigskip
	
	\textbf{Lemma : }{\itshape #1}
		
	\bigskip
	
	}

\newcommand{\Def}[1]
	{
	
	\bigskip
	
	\textbf{Definition : }{\itshape #1}
	
	\bigskip
	
	}

\newcommand{\Dem}[1]{
	
	\smallskip
	
	\textbf{Proof : } \par
	 {#1} $\square$
	 
	 \bigskip
}

\makeatletter
\newcommand*{\triplerightarrow}[1]{\mathrel{
  \settowidth{\@tempdima}{$\scriptstyle#1$}
  \mathop{\vcenter{
    \offinterlineskip\ialign{\hbox to\dimexpr\@tempdima+1em{##}\cr
    \rightarrowfill\cr\noalign{\kern.5ex}
    \rightarrowfill\cr\noalign{\kern.5ex}
    \rightarrowfill\cr}}}\limits^{\!#1}}}
\makeatother

\begin{document}

\pagestyle{plain}
\title{Regular polygraphs and the Simpson conjecture}

\date{}
\author{Simon Henry}

\maketitle

\begin{abstract}
We prove Carlos Simpson's ``semi-strictification'' (or ``weak unit'') conjecture in the case of $\infty$-groupoids.

More precisely, we introduce two precise versions of the conjecture, the ``general'' and the ``regular'' conjecture,  involving two different notions of ``non-unital categories''. The ``general'' version involve $\infty$-categories where absolutely all composition operations (horizontal, vertical and whiskering) are defined and compatible, the ``regular'' version involve $\infty$-categories where all the composition operations corresponding to ``regular'' pasting diagram are defined and compatible. In both case we construct (weak) model structures on these categories such that fibrant objects have weak units and weak inverse. We prove the regular version of the conjecture using the original strategy of Kapranov and Voevodsky, together with our previous work on polygraphs. The general version cannot be proved by these methods and is still open.

In order to do this we also study some subtle property of the combinatorics of polygraphs, and we construct a new counting function for polygraphs, inspired by previous work of Makkai.

\end{abstract}

\renewcommand{\thefootnote}{\fnsymbol{footnote}} 
\footnotetext{This work was supported by the Grant agency of the Czech republic under the grant P201/12/G028.}
\renewcommand{\thefootnote}{\arabic{footnote}}



\tableofcontents

\section{Introduction}

\blockn{In \cite{kapranov1991infty}, M.Kapranov and V.Voevodsky, while working on the homotopy hypothesis, erroneously claimed that the homotopy category of spaces is equivalent to an homotopy category of strict $\infty$-categories in which every arrow is weakly invertible.

In \cite{simpson1998homotopy}, C.Simpson proved that their statement is false. But he conjectured that it was true if strict $\infty$-category were replaced by weaker structure where both units and inverses are weak, but composition, associativity and exchange law are all strict. He moreover conjectured that the original paper of Kapranov and Voevodsky probably contained a proof of this claim, and made an other conjecture where he extended this to a similar strictification theorem for strict $\infty$-categories. To our knowledge the only works making progress on this conjecture are \cite{kock2006weak} which give one precise formulation of the conjecture (which we will not follows) and \cite{joyal2006weak} which proves this form of the conjecture for $3$-groupoids with only one object.
}

\blockn{The goal of the present paper is to answer this long standing conjecture in the case of $\infty$-groupoids.
}

\blockn{This paper follows \cite{henry2017positive} and we recommend looking at the appendix of \cite{henry2017positive} as an introduction to this one.}

\blockn{As we explained in this appendix, the main problem in \cite{kapranov1991infty} is that M.Kapranov and V.Voevodsky under estimated how complicated the notion of pasting diagram for $\infty$-category can be. In order for their strategy to work they need a well behave notion of $\infty$-categorical pasting diagrams. They decided to use the notion of pasting diagram introduce by Johnson in \cite{johnson1989combinatorics} and improved by them in \cite{kapranov1991combinatorial}, but it appears that this notion is largely insufficient for their purpose. Among other things, the composite of two such pasting diagrams is not always in this class of pasting diagram, making it impossible to have a well controlled strict $\infty$-category whose cells are parametrized by such diagrams. We refer to the appendix of \cite{henry2017positive} for a more detailed analysis of the strategy of Kapranov and Voevodsky and how this makes it fail.}

\blockn{The main result of \cite{henry2017positive} can be rephrased as the fact that if we work with ``non-unital strict $\infty$-category'' then there exists a (essentially unique) notion of pasting diagrams (the ``positive polyplexes'') that have all the desired property to be able to apply the original strategy of \cite{kapranov1991infty} to prove C.Simpson conjecture in the way intended by C.Simpson.}

\blockn{There is unfortunately a new difficulty arising: In this new, larger, class of pasting diagram (our positive polyplexes) there are some objects that are considerably more complicated than those considered by Johnson, Kapranov and Voevodsky. More precisely there are positive polyplexes (i.e. pasting diagrams) whose geometric realization is not contractible. The proof of Kapranov and Voevodsky relied heavily on this contractibility of the of the pasting diagram (which they prove for Johnson's diagram) and this failure of the contractibility for the more general notion of diagram prevent their proof to work even in this ``non-unital setting''.}

\blockn{It is a good place to mention that C.Simpson conjecture is not precise on what is an $\infty$-category with weak unit. In fact, even if we want to talk about non-unital $\infty$-category it appears that the notion is not unique: depending on the type of combinatorics one use to encode the different operation (globular, cubic, simplicial, opetopic, etc.) removing the units will not have the same effect at all.

In fact even in purely globular combinatorics, one can wonder whether a whiskering operation:

\[ \begin{tikzcd}[ampersand replacement=\&]
x \arrow[bend left=50]{r}[name=U,below]{}
\arrow[bend right=50]{r}[name=D]{} \&
y \arrow[Rightarrow,to path=(U) -- (D)]{} \arrow{r} \& z
\end{tikzcd} \]

should be defined in a non-unital $\infty$-category, or should be treated as a horizontal composition of $2$-cell with an identity $2$-cell. For another example, the simplicially based notion introduced by J.Kock in \cite{kock2006weak} seems very far from a globular formulation of non-unital $\infty$-category. It seems that the only test to see if something is a reasonable notion of ``non-unital $\infty$-category'' is to check that adding units to it give the notion of strict $\infty$-category (which is well defined and independent of the combinatorics).

}

\block{Different notion of non-unital $\infty$-categories, and different ways to give them weak units, will results in several possible form of the conjecture. We will actually state two of them the ``General Simpson conjecture'' and the ``Regular Simpson conjecture''. The first one corresponds to what we think is the strongest\footnote{We mean the one that makes the most operation into strict ones. We do not claim that other ``weaker'' form of the conjecture would automatically follows from it, although it might well be the case.} possible form of the Simpson conjecture (for $\infty$-groupoid). The second one to the strongest form of the conjecture that we will be able to prove using the methods of the original paper of Kapranov and Voevodsky. The General form of the conjecture is still open.
}

\blockn{The ``General Simpson conjecture'' refer to the version of the conjecture, involving all possible pasting diagrams.

More precisely: a non-unital $\infty$-category (see \cite{henry2017positive}) is a globular sets where all the operation $ x \#_i y$ are defined as soon as $x$ and $y$ are arrows (of possibly different dimension) whose $i$-dimensional source and target matches, with these operations satisfying all the usual axioms of associativity and exchange low. This is the structure needed to define composition of all of our ``positive polyplexes''.

We will construct in subsection \ref{subsec:FolkMS} a ``weak'' model structure on this category of ``non-unital $\infty$-categories'', whose fibrant objects are the non-unital $\infty$-categories which have weak units and weak inverse. Our ``General Simpson conjecture'' is essentially the fact that this model category is Quillen equivalent to the model category of spaces (we give a more precise statement in \ref{Conj:SimpsonGeneral}).

As we said earlier, using all possible pasting diagram involve some non-contractible pasting diagram and this keeps us from applying the strategy of \cite{kapranov1991infty} (see the appendix of \cite{henry2017positive}). At the moment this is written we have no serious leads on how to prove this ``general Simpson conjecture''. Though our result of \cite{henry2017positive} are very encouraging and show that, in some sense, there is is no ``obstruction\footnote{in the same way that the strict Eckmann-Hilton argument was an obstruction for the original statement of Kapranov and Voevodsky}'' for the conjecture to be true. }

\blockn{Our solution to answer C.Simpson conjecture despite this is to restrict the type of composition operation that we will make strict in our $\infty$-categories. We are typically not going to be able to make horizontal composition and whiskering:

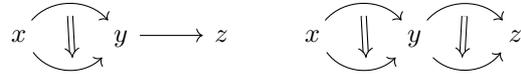
\[ \begin{tikzcd}[ampersand replacement=\&]
x \arrow[bend left=50]{r}[name=U,below]{}
\arrow[bend right=50]{r}[name=D]{} \&
y \arrow[Rightarrow,to path=(U) -- (D)]{} \arrow{r} \& z
\end{tikzcd} \qquad \begin{tikzcd}[ampersand replacement=\&]
x \arrow[bend left=50]{r}[name=U,below]{}
\arrow[bend right=50]{r}[name=D]{} \&
y \arrow[Rightarrow,to path=(U) -- (D)]{}  \arrow[bend left=50]{r}[name=V,below]{}
\arrow[bend right=50]{r}[name=W]{} \& z \arrow[Rightarrow,to path=(V) -- (W)]{}
\end{tikzcd} \]

into well defined strictly associative composition, as those appears to be sufficient to produce examples of non-contractible polyplexes. Instead, we will restrict to ``regular'' compositions operations, whose geometric realization are balls with good boundary inclusion like: 

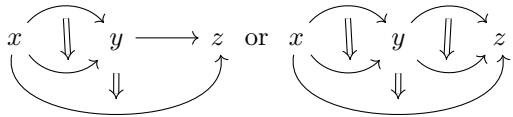
\[ \begin{tikzcd}[ampersand replacement=\&]
x \arrow[bend right=100]{rr}[name=k]{} \arrow[bend left=50]{r}[name=U,below]{}
\arrow[bend right=50]{r}[name=D]{} \&
y \arrow[Rightarrow,to path=(U) -- (D)]{} \arrow[Rightarrow,from=1-2,to=k,shorten >= 1pt,shorten <= 5pt]{}  \arrow{r} \& z 
\end{tikzcd} \text{ or }  \begin{tikzcd}[ampersand replacement=\&]
x \arrow[bend right=100]{rr}[name=k]{} \arrow[bend left=50]{r}[name=U,below]{}
\arrow[bend right=50]{r}[name=D]{} \&
y \arrow[Rightarrow,to path=(U) -- (D)]{} \arrow[Rightarrow,from=1-2,to=k,shorten >= 1pt,shorten <= 5pt]{}  \arrow[bend left=50]{r}[name=V,below]{}
\arrow[bend right=50]{r}[name=W]{} \& z \arrow[Rightarrow,to path=(V) -- (W)]{}
\end{tikzcd} \]}

\blockn{A ``Regular $\infty$-category'' will be a globular set on which every regular composition operation is defined and they are all compatible. More precisely we will define a notion of ``regular polyplexes with spherical boundary'' and a regular $\infty$-category will be a globular set where all such polyplexes can be composed. 

One can see that if one has any all regular composition operations, and all units, one can recover all composition operations. For example, the two diagrams above representing horizontal composition and whiskering can be obtain from our two examples of regular diagram by inserting units in the appropriate places.

In the exact same way as before we will construct a (weak) model category on the category of regular $\infty$-categories, whose fibrant objects are those that have weak inverse and weak units. And this time we will be able to prove, essentially following the work of Kapranov and Voevodsky, that this model category is equivalent to the model category of spaces.
}

\blockn{Finally, in his PhD thesis \cite{hadzihasanovic2017algebra} and in a very recent subsequent paper \cite{hadzihasanovic2018combinatorial}, Amar Hadzihasanovic has also independently introduced a notion of ``regular polygraphs''. His notion is introduced from a more combinatorial point of view, and his approaches has some important advantages: for example, it allows for a more explicit description of polyplexes and it makes the construction of the Gray tensor product considerably easier (though it is not proved yet that this is the same Gray tensor product as the ordinary one).

We had the opportunity to discuss about the connection between the two notions recently, and it seems not too difficult to prove that his regular polygraphs are a special case of mine. His ``globes'' are what I call regular polyplexes with spherical boundary and his ``atomic globes'' are what I call regular plexes. Unfortunately the two notions are not equivalent. This is due to a restriction Hadzihasanovic imposes on globes: He asks that non-atomic globes can always be partitioned as a composition of two sub-globes satisfying some compatibility conditions, but it appears that there are examples of regular $3$-polyplexes with spherical boundary with three $3$-dimensional cell so that is no way to form a $3$-polyplexes with spherical boundary by composing any two of the three $3$-cell. A $4$-plex with such a $3$-polyplex has its source and target would be a regular polygraph in the sense of the present paper, but not in the sense of \cite{hadzihasanovic2018combinatorial}.

But the two notions are still very similar, and it seems that the main results of the present paper can be extended without too much changes to Hadzihasanovic notion of regular polygraphs. The corresponding notion of regular $\infty$-categories would have slightly less operations: only the regular composition of two arrows would be possible (and the operation that are obtained by iterated composition of such).

}

\subsection{Sketch of the proof and plan of the paper}

\blockn{The proof relies heavily on the theory of ``plex'' and ``polyplexes'' developed in \cite{henry2017positive}, which will be absolutely essential for the present paper. The final section (dealing with all the homotopy theoretic aspect) will also relies on a companion paper, published simultaneously to this one, \cite{henry2018weakmodel}, which introduce and develop the basic theory of a so called ``weak model category'' ( a weakening of Quillen model structure which includes both left and right semi-model structure as examples). We believe that for a reader familiar with model category, a quick look at \cite{henry2018weakmodel}  ( section 2.1-2.3 and theorem 3.2.1) should be sufficient to follow these aspects. }

\blockn{As many mathematical paper, this might not be suited to be read from the first page to the last. We will explain the role of each section, and the general structure of the proof.}

\blockn{Subsection \ref{subsec:infinityCat} contains very basic preliminary on strict $\infty$-category and fix our convention. Subsection \ref{subsec:polygraph} recall the basic notion of polygraphs and the notion of plex and polyplexes from \cite{henry2017positive}. Subsection \ref{subsec:generic} contain a small extension of results in \cite{henry2017positive} and a review of the notion of generic factorization system of M.Weber (\cite{weber2004generic}) that will be useful at several point in the paper.}

\blockn{Here are the main step of the proof:

\begin{itemize}

\item The definition of the notion of ``Regular $\infty$-category'': we define first what are ``regular polygraphs'' in subsection \ref{subsec:regularPol}. And we define more generally a notion of $C$-category for $C$ a class of polygraphs satisfying certain good condition in section \ref{subsec_algebraicclass_def}. This applies both to ``positive polygraph'' to give non-unital $\infty$-category and to ``regular polygraphs'' to give regular $\infty$-category. The fact that the class of regular polygraph have the desired property basically mean that ``a regular composition of regular pasting diagram is again regular'', it sounds obvious, but it have appeared harder to proof than expected and the proof will be postponed to corollary \ref{Cor:RegularPolygraphAreAnAlgebraicClass} (discussed below).

\item Subsection \ref{subsec_Gray_tensor_product} present a well known monoidal structure (the Gray or Crans-Gray tensor product) on the category of $\infty$-category and on the category of polygraph. The end of this subsection prove that the category of regular polygraphs is stable under the Gray tensor product, which is a difficult result that involve a lot of the machinery introduce in the rest of the paper and might be skipped in a first time.

\item Section \ref{subsec_tensorstable_algebraic_class} show that if the class of polygraphs $C$ is stable under tensor product then the category of $C$-categories also admit a ``Gray tensor product''.

\item Section \ref{subsec:FolkMS} construct the ``weak model structure'' (in the sense of \cite{henry2018weakmodel}) on the categories of $C$-categories, and a closely related weak model structure on the category of $C$-polygraphs mentioned in the introduction. This is done for any good class of polygraphs (in the sense of \cite{henry2017positive}) which is stable under the Gray tensor product.

\item Section \ref{subsec:QuillenEquivCatPolyg} shows that these two model categories (of $C$-category and $C$-polygraphs) are Quillen equivalent by the natural adjunction between them.

\item Finally section \ref{subsec:QuillenequivSpace} contains the proof that this model category of regular polygraphs is equivalent to the model category of spaces, hence completing the proof of this equivalence between regular $\infty$-category and spaces. Section \ref{subsec:Conjecutres} contains some remarks about the General Simpson conjecture and the extension of the Simpson conjecture to weak $(\infty,n)$-categories. The reason why this does not work for non-unital $\infty$-category (instead of regular) is explained in the appendix of \cite{henry2017positive}.

\end{itemize}
}

\blockn{
In the sketch above, we have left aside the proof that regular polygraphs are stable under tensor product and satisfies the conditions required for the notion of regular $\infty$-category to make sense, as well as a few technical tools that will be involved in these proof as well as in subsection \ref{subsec:QuillenequivSpace}. Those have appeared to be the most difficult result of the present paper, and occupied almost half of it. In fact at the time we wrote the appendix of \cite{henry2017positive} we already had all the element mentioned above in place, but we were unable to prove these last two claims. The proof of these fact actually required some new ideas that I will present before finishing the overview of the paper:

When looking at the definition of the regular polygraphs, it sounds quite natural (for example by looking at low dimension example) that they will be stable under the Gray tensor product. But the truth is that we have very little control over the Gray tensor product of $\infty$-categories (or polygraphs). We know that the tensor product of two polygraphs $X \otimes Y$ have generating cells given by the ``$x \otimes y$'' for $x$ a generating cell of $X$ and $y$ a generating cell of $Y$, but in general it is nearly impossible to give explicit formula for the source and target of $x \otimes y$. It is hence very difficult to prove that $x \otimes y$ satisfies some regularity condition (which are condition on the shape of its iterated source and target).

The only information about the source and target of $x \otimes y$ to which we have an easy access is the number of occurrence of cells of dimension $(n-1)$ in $x \otimes y$ (where $n=dim(x)+dim(y)=dim(x\otimes y)$). This is given to us by the fact that the linearization functor attaching to each $\infty$-category its homology complex is monoidal for the Gray tensor product, and that we know how to read this number of occurrence in the homology complex. In particular, one can detect on the homology complex if these source and target are identity cells or not, and this observation allows to see that the tensor product of two positive (or non-unital) polygraphs $X$ and $Y$ is again positive.}

\blockn{The key idea in order to show that regular polygraphs are stable under tensor product (and also to show that ``regular composition of regular composition are regular''), is the introduction of a new, more subtle, ``counting functions'' which attribute to each $n$-arrow $f$ of a free $\infty$-category $X$ a formal linear combination $\delta(f)$ of the generating cells in $X$ that are involved in the construction of $f$. The $n$-dimensional generators appears in this with, as coefficient, the number of time they appears in any expression of $f$, but lower dimensional generator appears with more subtle (possible negative) coefficients. This new counting function will be sufficient to characterize the regularity of a polygraph, and it will have both good functoriality property and good compatibility with the Gray tensor product. This will allow us to both prove that stability of regular polygraph under tensor product, and this condition that ``regular composition of regular compositions are regular''.

}

\blockn{ The function $\delta$ is also somehow encoded in the homology complex, but not in a visible, easy to ready way. I first discovered this counting function as an ``artificial combinatorial trick'', but it quickly appears that it had a very natural categorical original that nicely explained all its mysterious property and its connection to the homology complex. It is defined as follows:

One knows that the category of abelian group objects in the category of $\infty$-category is equivalent to the category of chain complexes (see \cite{brown1981equivalence}). But it appears that this equivalence is non trivial, and despite these two categories being equivalent, some information seems to be lost (or hidden) when passing from group object to chain complexes.

This $\delta$ function is just the natural map from the free category $X^*$ on a polygraph $X$ to the universal abelian group object $\Z X^*$ generated by $X^*$; and $\Z X^*$ naturally identifies as a structure of $\infty$-category on the free abelian group $\Z X$ on the generating cells of $X$.*

If one see this in terms of the homology complex of $X^*$, the non-trivial equivalence between abelian group object (strict $\infty$-group) and chain complexes hide a lot of information present in $\delta$ and only let the coefficients of $n$-dimensional generators appears, which is just the ordinary counting function. But the lower dimensional generators also appears non trivially in $\delta(x)$ and they contains a lot of information.
}

\block{To finish the sketch of our proof:

\begin{itemize}

\item In subsection \ref{subsection_StrictAbGrp_vs_chaincomplexes} on review the equivalence between ``$\infty$-group'', i.e. abelian group object in $\infty$-category (strict $\infty$-groups), ''globular groups'', and chain complexes.

\item In subsection \ref{sub_sec_linearization} we construct the universal $\infty$-group attached to an $\infty$-category and define the function $\delta : X^* \rightarrow \Z X$.

\item In subsection \ref{subsec:TensorProdOfInfinityGroup} we give formula for the tensor product of globular groups/$\infty$-group coming from the equivalence with chain complexes and the usual tensor product of chain complexes.

\item Section \ref{section:Delta} is devoted to the study of the function $\delta$ in more details. Subsection \ref{Subsec_Delta_counting_general} start with its general property as a counting function. 

\item Subsection \ref{subsec_delta_counting_regular} contains the key results allowing to characterize regularity in terms of properties of the function $\delta$. One deduce in particular the condition that ``regular composition of regular composition are regular'' (Corollary \ref{Cor:RegularPolygraphAreAnAlgebraicClass}).

\item After introducing the Gray tensor product, subsection \ref{subsec_Gray_tensor_product} shows the compatibility between the function $\delta$ and the tensor product and use the result of the previous subsection to show that regular polygraphs are stable under tensor product.

\item Finally, subsection \ref{subsec:joinsAndCones} define a ``cone'' construction on polygraphs which play an important role in the homotopy theoretic property of regular polygraph (it shows for example that regular plexes form a weak test category) and which will be needed in the proof of the equivalence between regular polygraphs and spaces in subsection \ref{subsec:QuillenequivSpace}. 

\end{itemize}

}

\section{Preliminaries and definitions}

\subsection{$\infty$-categories}
\label{subsec:infinityCat}
\block{\label{Def_StreetInfinityCat}We will use Street's ``one type'' definition of $\infty$-category, the $k$-dimensional source and target of an arrow $f$ are denoted $\pi^-_k f $ and $\pi^+_k f$, compositions are in diagrammatic order and are denoted $f \#_k g$ where $k$ is the dimension of the common boundary of $f$ and $g$ along which they are composed. More precisely:

\Def{An $\infty$-category is a set $X$ together with unary operations $\pi^-_k$ and $\pi^+_k$ for $k\geqslant 0$ and partially defined binary operations $\#_n$ for $n\geqslant 0$ satisfying the following axioms:

\begin{enumerate}

\item $ x \#_n y $ is defined if and only if $\pi^+_n(x)=\pi^-_n(y)$.

\item For every $x \in X$ there exists an $n$ such that $\pi^-_n(x)=\pi^+_n(x)=x$.

\item For any $x \in X$, one has:

\[ \pi^{\epsilon}_n \pi^{\delta}_m x = \left\lbrace \begin{array}{l l} \pi^{\epsilon}_n x & \text{if } n <m \\ \pi^{\delta}_m x & \text{if } n \geqslant m\end{array} \right. \]

\item For any $x \in X$ one has $\pi^-_n(x) \#_n x = x \#_n \pi^+_n(x) = x $.

\item For all $x,y$ such that the composition below are defined,

\[ \pi^-_n(x \#_n y) = \pi^-_n(x) \]
\[ \pi^+_n(x \#_n y) = \pi^+_n(y) \]

And if $k>n$:

\[ \pi^{\epsilon}_k(x \#_n y) = \pi^{\epsilon}_k(x) \#_n \pi^{\epsilon}_k(y)\]

\item $x \#_n (y \#_n z) = (x \#_n y) \#_n z$ when either side is defined.

\item If $k <n$

\[ (x \#_n y) \#_k ( z \#_n w) = (x \#_k z) \#_n (y \#_k w) \]

when the left hand side is defined.

\end{enumerate}

A morphism of $\infty$-category is just a function commuting to all the structural functions (source, target and the partial composition).
}

Morphisms are often called functors. As usual, a $n$-arrow is an element $x\in X$ such that $\pi^{\epsilon}_n x = x$. The dimension of an arrow $x\in X$ is the smallest $n$ such that $x$ is a $n$-arrow, i.e. an arrow of dimension $n$ is a ``non-identity $n$-arrow''.
}

\block{\label{Discuss_limits}Except for the second axiom, the definition of $\infty$-category given above make them model of a partial algebraic theory. In particular, if we forget this second axioms limits of $\infty$-category exists and are computed at the level of the underlying sets, with all the operation being defined ``component wise''. One easily see that a finite limits of $\infty$-category in this sense still satisfies the second axiom as well, so finite limits of $\infty$-categories are computed at the level of the underlying set of cells (this is not true for infinite products: for those one needs to further restrict to the set of cells which satisfies this second axioms).

The category of $\infty$-categories also have all colimits and all ``free constructions'' for the same reasons (but of course they are not computed at the level of the underlying set).
}

\block{
\begin{itemize}

\item One denotes by $D_n$ the $\infty$-category freely generated by an $n$-arrow $*$ or $*_n$. Its only arrow are $*_n$ and the $\pi_k^{\epsilon} *_n$ for $k<n$ often just denoted $\pi_k^{\epsilon}$ for $k<n$. All compositions are trivial (i.e. every time $x\#_k y$ is defined axioms $4.$ applies and the result is $x$ or $y$.).

\item $D_0$ is the terminal $\infty$-category and it will sometime be denoted $*$.

\item One denotes by $\partial D_n$ the $\infty$-category freely generated by two parallel $(n-1)$-arrows $\pi^+_{n-1}$ and $\pi^-_{n-1}$. $\partial D_{n}$ is exactly the subobject of $D_n$ which contains all arrows except $*_n$.

\end{itemize}

A functor $f:D_n \rightarrow X$ is denoted by the name of the arrow $f \in X$ which is the image of $*_n$. A functor $(s,t) :\partial D_n \rightarrow X$ is denoted by the pair of parallel $n$-arrows respectively image of $\pi^-_{n-1}$ and $\pi^+_{n-1}$.
}

\subsection{Polygraphs, plexes and polyplexes}
\label{subsec:polygraph}

\block{Polygraphs are the most general notion of $n$-categorical diagram: they are the $\infty$-category freely generated by certain (explicitly given) cells.

Very roughly, an $n$-polygraph $X$ and its free $\infty$-category $X^*$ are defined inductively as follows: a $n$-polygraph $X$ is given by a $(n-1)$-polygraphs $X^{(n-1)}$ together with a set $S$ of pairs of parallel $n-1$ arrow $(s,t) : \partial D_n \rightarrow (X^{(n-1)})^*$. The free $\infty$-category $X^*$ is obtained by freely adding one arrow to $(X^{(n-1)})^*$ with prescribed source and target for each pair of selected parallel $n-1$-arrows, i.e. constructing the pushouts:

\[ (X^{(n-1)})^* \coprod_{\partial D_n} D_n \]

for each map $\partial D_n \rightarrow (X^{(n-1)})^*$ in $S$. Morphisms of polygraphs are also defined by induction in the same way, but can equivalently be defined as the functor between the free $\infty$-categories that send generating cells to generating cells.

Polygraphs were initially introduced by R.Street under the name ``computades'' in the framework of $2$-category theory (in \cite{street1976limits}) the $n$-categorical extension was alluded to in several subsequent paper of R.Street, but first appears spelled out explicitly in \cite{johnson1989combinatorics} and under the name ``polygraphs'' in \cite{burroni1993higher}.  We refer to \cite{burroni1993higher} or to any of the reference mentioned above for a more detailed an introduction to the topic. We also mention J.Penon's \cite{penon1999approche} for an equivalent, syntactic presentation of polygraphs, and M.Makkai's \cite{makkai2005word} for a lot of interesting and subtle properties of polygraphs.

We will use the following terminology:

\begin{itemize}

\item The category of polygraphs, with morphisms of polygraphs between them is denoted $\Pb$.

\item The word ``cell'' is reserved for the generators of $X$, but one also sometimes use ``generators'' or ``generating cells'' to emphasize that, those are all synonymous. A polygraph $X$ is generally thought of as the set of its generating cells, with an additional structure specifying the source and target of each cell in $X^*$. For example $x \in X$ means that $x$ is a generating cell of $X$.

\item The word ``arrow'' will denote any arrow in the free $\infty$-category $X^*$ generated by $X$, but we will (try to) avoid the use of the expression ``arrow of $X$'' and always says ``arrow of $X^*$''.

\item Similarly, we will try to forbid the expression ``cells of $X^*$'' to avoid possible confusion with arrows of $X^*$, as arrows in an $\infty$-category are sometimes called cells. We will prefer to the expression ``generating cell''.

\item A morphisms $f :X^* \rightarrow Y^*$ is said to be polygraphic if it send any generating cell to a generating cell. It is well known that polygraphic morphisms are exactly those that are induced by morphisms of polygraph, and that this is a bijection between morphisms of polygraphs and polygraphic morphisms. Also any isomorphisms between $X^*$ and $Y^*$ is polygraphic and in particular the $\infty$-category $X^*$ determine the polygraph $X$ up to unique isomorphism.

\end{itemize}}

\block{$D_n$ and $\partial D_n$ constructed above are very naturally examples of polygraphs (were every arrows is a generating cells). They are the only example were we will use the same notation of the polygraphs and the free $\infty$-category (we allows this abuse of notation as in this case they have the same underlying set). The category of polygraphs has a terminal objects $\Pb 1$, which is an extremely rich and important object characterized by the following two properties:

\begin{itemize}

\item $\Pb 1$ has only one cell $*$ in dimension $0$.

\item For each pair of parallel $n$-arrow $(s,t)$ in $(\Pb 1)^*$, $\Pb 1$ has exactly one $n+1$-cells from $s$ to $t$.

\end{itemize}

So for example, it has only one $1$-cell $f:* \rightarrow *$, its $1$-arrows are all the $n$-folds composite $f^n$ for each $n \geqslant 0$. Hence it has exactly one $2$-cell $A_{n,m}:f^n \rightarrow f^m$ for each $n,m \geqslant 0$. Cells of dimension $\geqslant 3$ very quickly become too complicated to be explicitly listed.

}

\block{ We note the following very classical lemma:
\label{Lem_2outof3_for_polygraphic}
\Lem{If one has a diagram of functors:

\[\begin{tikzcd}[ampersand replacement=\&]
X^* \arrow{dr}{h} \arrow{r}{f} \& Y^* \arrow{d}{g} \\
 \& Z^*
\end{tikzcd}\]

between free $\infty$-category on polygraphs. If $h$ and $g$ are polygraphic then so is $f$.
}

This is well known, and we give a proof only for completeness. This proof relies heavily on the results of section \ref{Subsec_Delta_counting_general} and proposition \ref{Prop_delta_count_in_top_Dim}.

\Dem{ We only need to prove that given an arrow $y \in Y^*$ such that $g(y)$ is a generating cell of $Z$, then $y$ is a generating cell of $Y$. Applying this to the cell $h(x)$ for $x$ a generating cell of $X^*$ will prove the result. So let $y$ be such an arrow in $Y^*$ on has $\delta(y) = \sum n_i \delta_{y_i}$ and $\delta(h(y))= \sum n_i \delta_{h(y_i)} = \delta_z$. If all the $n_i$ were positive one could conclude immediately that hence they are all $0$ except one, hence $\delta(y)=\delta_{y_i}$ hence $y=y_i$ is a generator by the last point of proposition \ref{Prop_delta_count_in_top_Dim}.

But thanks to proposition \ref{Prop_delta_count_in_top_Dim} again, one can use this argument to prove by downward induction on $k$ that $\delta(y)=\delta_{y_i} + u$ with $u$ of dimension at most $k$. Indeed this is clearly true for large $k$, if $\delta(y)=\delta_{y_i} + u$ for $u$ of dimension at most $k$ then proposition \ref{Prop_delta_count_in_top_Dim} show that cell of dimension $k$ or greater appears with non-negative coefficients, hence the argument above show that there is at most one non zero coefficients of dimension $k$ or greater in $\delta(y)$ which is exactly our induction hypothesis for $k-1$. Applying this to $k=0$ gives that $\delta(y)=\delta(y_i)$ and concludes the proof.

}
}

\block{We recall from \cite{henry2017positive}:

\Def{A polygraph $X$ is said to be positive if the source and target of any generating cell of dimension $n$ is a non-identity arrow of dimension $n-1$. One denotes by $\Pb^+$ the category of positive polygraphs.}

by non-identity arrows we just mean an arrow which is actually of dimension $n-1$ and not of lower dimension (note that even constructively, the dimension of an arrow in a polygraph, i.e. the smallest $n$ such that $\pi_n^+(f)=f$ is well defined).

Positive polygraphs forms a ``class of polygraphs'' in the sense of \cite{henry2017positive}. This means that there is a terminal positive polygraphs $\Pb^+ 1$, it is the subobject of $\Pb 1$ of all cells $x$, whose source and target are non-identity and themselves belongs to $\Pb^+ 1$. A polygraph is positive if and only if its unique map to $\Pb 1$ factors in $\Pb^+ 1$.
}

\block{We will now recall the main results of \cite{henry2017positive}, we only state them from the class of positive polygraph, but this kind of result holds for any ``good class of polygraph'' (see \cite{henry2017positive}):

\label{Th_def_plex_polyplex}
\Th{\begin{itemize}

\item The category of positive polygraphs is a presheaf category on a (uniquely\footnote{because the category of plexes have no non-identity idempotent} defined) category of objects called (positive) ``plex''.

\item There is a family of positive polygraph called (positive) ``polyplex'' and an isomorphism functorial in $X$ a positive polygraph:

\[ \{\text{arrows of } X^* \} \simeq \coprod_{\underline{p} : \text{polyplex}} Hom(\underline{p},X) \]

The plexes are a special case of polyplexes, they corresponds exactly to those that represents generating cells of $X$ i.e.:

\[\{\text{cells of $X$}\} \simeq \coprod_{\underline{p} : \text{Plex}} Hom(\underline{p},X)  \] 

\end{itemize}}

Note that because of the isomorphisms of the theorem, each polyplex $\underline{x}$ is canonically endowed with an arrow $x \in \underline{x}^*$, corresponding to the identity map $\underline{x} \rightarrow \underline{x}$. This arrow $x$ is called the ``universal arrow'' of the polyplex $\underline{x}$.

The isomorphisms of the theorem can be reformulated as follow: 

For any polygraph $X$, and any arrow $f \in X^*$ there exists a unique polyplex $\underline{p}$ and a unique map $\chi_f: \underline{p} \rightarrow X$ such that $f = \chi_f(p)$.

One easily see that this properties characterize polyplexes uniquely up to unique isomorphism as polygraph endowed with a chosen arrow. Indeed the property above can be rephrased as the fact the category of pairs $X \in \Pb^+$ and $f \in X^*$ is a disjoint unions of categories with initial objects, the polyplex being these local initial objects.

But it should be noted that polyplex are uniquely determined, only has objects in this category of polygraphs endowed with an arrow $d \in X^*$. In particular we have given in \cite{henry2017positive} an example of two different polyplexes which are isomorphic as polygraphs, i.e. a polygraph which is a polyplex in two different way, for two different universal arrows.

Plexes are exactly the polyplexes whose universal arrow is a generating cells. Moreover a plex has a unique generating cell of maximal dimension, so contrary to what happens with polyplexes a polygraph can be (isomorphic to) a plex in only one way.

Note that what we call here plex and polyplex were called ``positive plex'' and ``positive polyplex'' in \cite{henry2017positive}, where it was shown that there are larger class of polygraphs satisfying the above theorem and hence having also a notion of plexes and polyplexes, also in \cite{makkai2005word}, M.Makkai introduced a notion of ``computope'' that is equivalent to our notion of plex when restricted to a good class of polygraphs, but makes sense for general polygraphs. I.e. in the present paper, plex and polyplex are always implicitly assumed to be ``positive''.
}

\block{\label{discuss:InductiveConstructionOfPolyplex}
Also note that taking $X$ to be the terminal positive polygraph $\Pb^+1$ in theorem \ref{Th_def_plex_polyplex} produce bijections between the set of arrows of $(\Pb^+1)^*$ and the set of polyplexes, which restrict to a bijection between the set of plexes and the set of generating cells of $\Pb^+ 1$. We will identifies arrows of $(\Pb^+ 1)^*$ with polyplexes.

This means in particular that the set of polyplexes is equipped with the structure of an $\infty$-category. The proof of the main theorem of \cite{henry2017positive} involves and relies on a description of these operations on polyplexes in terms of the underlying polygraphs, that produces a way to describe plexes and polyplexes inductively as follows:

\begin{itemize}

\item If $\underline{p}$ is a polyplex, then $\pi_k^{\epsilon} \underline{p}$ is the polyplex representing the arrow $\pi_k^{\epsilon} p \in \underline{p}^*$. In particular one has a unique map $\pi_k^{\epsilon} \underline{p} \rightarrow \underline{p}$ sending the universal arrow of $\pi_k^{\epsilon} \underline{p}$ to $\pi_k^{\epsilon}p$.

\item If $\pi_k^+ \underline{p} = \pi_k^{-} \underline{q}$ then the underlying polygraph of $\underline{p} \#_k \underline{q}$ is:

\[ \underline{p} \coprod_{\pi_k^+ \underline{p} } \underline{q} \]

Its universal arrow is simply $p \#_k q$, which are indeed composable as arrows of this pushout.

\item Given $\underline{x}$ a generating $n+1$-cells of $\Pb^+ 1$,  $\pi_n^- \underline{x}$ and $\pi_n^+ \underline{x}$ are parallel $n$-polyplexes, i.e. they both have natural maps from $\pi_{n-1}^- \underline{x} \coprod \pi_{n-1}^+ \underline{x}$. The plex $\underline{x}$ is constructed as follows:

\[  \pi_n^- \underline{x} \coprod_{\pi_{n-1}^- \underline{x} \coprod \pi_{n-1}^+ \underline{x}}  \pi_n^+ \underline{x}  \]

together with a single additional $n+1$-cell between $\pi_{n}^- x$ and $\pi_{n}^+ x$ which have been parallel in this pushout.

\end{itemize}

}

\block{\label{discuss:polygraphic_induction}To sum up: we have an inductive formula to construct the cells and arrows of $\Pb^+ 1$, and we now how to construct the corresponding polyplexes. This allows to construct all polyplexes explicitly, but the inductive nature of this construction makes it difficult to give an explicit combinatorial description of what plexes and polyplexes are. 

In particular, most results about polyplex that we will prove will have to be proved by induction. The induction scheme is the following, if $C( \underline{p})$ is a property of a polyplex such that:

\begin{itemize}

\item $C$ holds for the unique $0$-dimensional plex.

\item If $C$ holds for all $n$-polyplex then $C$ holds for all $(n+1)$-plex. 

\item If $C$ holds for all $n$-polyplex, then, if $C(\underline{p})$ and $C(\underline{q})$ holds for two $(n+1)$-polyplexes and $k \leqslant n$ is an integer such that $ \underline{p}\#_k \underline{q}$ exists then $C(\underline{p} \#_k \underline{q})$ holds as well.
\end{itemize}

then $C$ holds for all polyplex.

}

\block{\label{Lem:syntactic_lifting}\Lem{let $f :X \rightarrow Y$ a map between two positive polygraphs. Let $v \in X^*$ an arrow and assume that $f(v)$ can be written as $u \#_k w$ in $Y^*$. Then there is a unique pair of arrows $u'$ and $w'$ in $X^*$ such that $u' \#_k w' = v$ and $f(u')=u$, $f(w')=w$. }

This type of property is closely related to the fact that the notion of polyplex exists and is well behaved. We leave as an exercise to the reader to find a counter-example to this proposition when $X$ and $Y$ are not assumed to be positive. It will of course involve the Eckmann-Hilton argument.

\Dem{Let $\underline{v} \rightarrow X$ be the polyplex representing $v \in X^*$, as $f :X \rightarrow Y$ is polygraphic, the uniqueness of the polyplex representing an arrow shows that $\underline{v} \rightarrow X \rightarrow Y$ is the polyplex representing $f(v)$ in $Y^*$.

Let $\underline{u},\underline{w} \rightrightarrows Y$ the polyplexes representing $u$ and $v$ in $Y^*$, as $u \#_k w = f(v)$ one has that $\underline{u} \#_k \underline{w} = \underline{v}$. But using this decomposition of $\underline{v}$ gives a decomposition of $v$ in $X^*$ which has the desired properties. Conversely, any other decomposition of $v$ with these properties have to give rise to the same decomposition of $\underline{v}$ and hence there is only one such pair of arrows.
}
}

\subsection{The generic factorization}
\label{subsec:generic}

\blockn{This section is mostly here to recall the notion of generic map and generic factorization of M.Weber from \cite{weber2004generic} (where it was called ``strict generic''), and to show that they applies to any good classes of polygraphs. We will only state them for the class $\Pb^+$ of positive polygraph, but all the results we will mention here still holds if one replaces $\Pb^+$ by good class of polygraph in the sense of \cite{henry2017positive}. The main results of this section is: }

\block{\label{Th_generic_factorization}\Th{The category of positive polygraphs, and functors between their corresponding free $\infty$-category has an orthogonal factorization system whose right class are the polygraphic morphisms. Element of the left class are called ``generic morphisms''.}

As explained before, polygraphic morphisms are those corresponding to morphisms of polygraphs. We also want to give some understanding of what are the generics morphisms and how this factorization is constructed and what is the relation of this with polyplexes.
}

\blockn{
It follows from our results in \cite{henry2017positive} that the free $\infty$-category monads on the category of positive polygraphs is a parametric right adjoint functor (i.e. the functor $T:\Pb^+ \rightarrow \Pb^+_{/T1}$ is a right adjoint functor). It then follows from \cite{weber2004generic} that every morphisms between free $\infty$-categories on polygraphs can be factored uniquely as a ``$T$-generic'' morphisms (in the sense of \cite{weber2004generic}) followed by a polygraphic morphisms. But it is unclear that generic morphisms are stable under composition. This would also follow formally from the results of \cite{weber2004generic} if one knew that multiplication transformation of the free $\infty$-category monad is cartesian . We don't see how to obtain this without a lot of work, but we already know that this monad has a lot of ``cartesian like'' properties, which apparently are sufficient to deduce the result, at the cost of changing the proof strategy.}

\block{
One say that a morphism $f:X^* \rightarrow Y^*$ is generic if it has the unique left lifting property against all polygraphic morphisms, i.e. $f$ is generic if for all solid square:

\[\begin{tikzcd}[ampersand replacement=\&]
X^* \arrow{d}{f} \arrow{r} \& A^* \arrow{d}{v^*} \\
Y^* \arrow[dotted]{ur} \arrow{r} \& B^* \\
\end{tikzcd}\]

with $v^*$ a polygraphic  morphisms, there is a unique dotted diagonal filling.

This is not exactly the definition that Weber is using, but it will be equivalent due to the main theorem of this subsection. We will show that any morphisms $X^* \rightarrow Y^*$ can be factored as a generic morphism followed by a polygraphic morphism. The uniqueness of such a factorization follows formally from the lifting property. and with this definition it is immediate that composition of generics is generic.

}

\block{\label{Prop_genericvspolyplex}\Prop{

\begin{itemize}

\item If $\underline{p}$ is a polyplex, then the morphism $p : D_n \rightarrow \underline{p}^*$ corresponding to its universal arrow is generic.

\item if $\underline{p}$ and $\underline{q}$ are parallel polyplex, then the map $\partial D_n \rightarrow \underline{p} \coprod_{\partial \underline{p}} \underline{q}$ corresponding to the two universal arrows is generic. 

\end{itemize}
}

 $\partial \underline{p}$ denotes the boundary of the polyplex $\underline{p}$, in the sense of:

\[ \partial \underline{p} = \pi^+_{n-1}\underline{p} \coprod_{\pi^+_{n-2} \underline{p} \coprod \pi^-_{n-1} \underline{p} } \pi^-_{n-1}\underline{p} \]

Note that the map $\partial \underline{p} \rightarrow \underline{p}$ might not be a monomorphism (for example if $\underline{p}$ is a horizontal composition).

\Dem{We start by proving that the second claim follow from the first:

Indeed, given a diagram of the form:

\[
\begin{tikzcd}[ampersand replacement=\&]
\partial D_n \arrow{d} \arrow{r} \& A^* \arrow{d}{w^*} \\
\underline{p}^* \coprod_{\partial \underline{p}} \underline{q}^* \arrow{r} \& B^*
\end{tikzcd}
\]

One can construct lifts $\underline{p}^* \rightarrow A^*$ and $\underline{q}^* \rightarrow A^*$ separately using property one, and by uniqueness of the diagonal filler these lift will agree on $\pi^+_{n-1} \underline{p}$ and $\pi^-_{n-1} \underline{q}$ hence extend into a map $\underline{p}^* \coprod_{\partial \underline{p}} \underline{q}^* \rightarrow A^*$. The uniqueness of this diagonal filler follows from the uniqueness applied to $\underline{p}$ and $\underline{q}$ separately.

We will prove the first claim by induction on $\underline{p}$.

Consider a square of the form:

\[
\begin{tikzcd}[ampersand replacement=\&]
D_n \arrow{d}{p} \arrow{r} \& A^* \arrow{d}{w^*} \\
\underline{p}^* \arrow{r} \& B^*
\end{tikzcd}
\]

with $w^*$ polygraphic. If $\underline{p}$ is of dimension $0$ then there is clearly a unique diagonal filler: the map $D_n \rightarrow A^*$ has to factor into a $0$-dimensional cell as well, and this cell give the lift. Assume the result is true for all polyplex of dimension $<n$ and let $\underline{p}$ be an $n$-plex then. Then $\partial D_n \rightarrow \partial \underline{p}$ is already known to be generic (by the induction hypothesis and the proof that the second condition follows from the first), hence one can already construct a unique lift $\partial \underline{p} \rightarrow A^*$.

$\partial \underline{p}$ is constructed from $\underline{p}$ by simply adding to it a single $n$-dimensional cell, so this cell is lifted as the image of the $n$-dimensional cell of $D_n$, and this proves both the existence and the uniqueness of the lift.

Finally, if $\underline{p} = \underline{q} \#_k \underline{r}$ with $\underline{q}$ and $\underline{r}$ satisfying the proposition. Then one factor the map $D_n \rightarrow \underline{p}$ into:

\[D_n \rightarrow D_i \coprod_{D_k} D_j \rightarrow \underline{q}^* \coprod_{\pi_k^+ \underline{p}^*} \underline{r}^* = \underline{p}^* \]

The first map is generic exactly because of lemma \ref{Lem:syntactic_lifting}, the second map is generic as a pushout of generic maps.

}
}

\block{\Cor{Any arrow $f:D_n \rightarrow X^*$ admits a generic factorization into $D_n \rightarrow \underline{p}^* \rightarrow X^*$ where $\underline{p}$ is the polyplex representing $f$. In particular $(X,f)$ is a polyplex if and only if the map $f:D_n \rightarrow X^*$ is generic.}

Indeed the previous proposition shows that $D_n \rightarrow \underline{p}$ is generic, and due to the lifting property of generic maps, generic factorization is unique up to unique isomorphisms.

}

\block{\label{Proof_Th_generic_factorization}We are now ready to prove theorem \ref{Th_generic_factorization}. Note that the proof below give an explicit description of what generic morphisms $X^* \rightarrow Y^*$ look like: $Y^*$ is constructed by iteratively gluing one polyplex along its boundary for each cell of $X^*$.

\Dem{We will prove by induction on the number of cells of a finite polygraph $X$ that any morphisms of $\infty$-category $X^* \rightarrow Y^*$ admit a generic factorization. If $X$ has zero or one cell, then this is trivial.

Assume that this is known for a given polygraph $X$ and consider $X'$ a polygraph obtained from $X$ by adding one cell $w$. Let $f:(X')^* \rightarrow Y^*$ be a morphism, one first considers the generic/polygraphic factorization of the restriction of $f$ to $X$:

\[\begin{tikzcd}[ampersand replacement=\&]
X^* \arrow[hook]{d} \arrow{r}{g_1} \& P^* \arrow{r}{u^*} \& Y^* \\
X'^* \arrow{rru}[swap]{f}
\end{tikzcd}\]

The map $D_n \overset{w}{\rightarrow} X' \rightarrow Y^* $ and its restriction to $\partial D_n$ also admit generic-polygraphic factorization due to the previous result this gives us a solid diagram:

\[\begin{tikzcd}[ampersand replacement=\&]
X^* \arrow{r}{g_1} \& P^* \arrow{dr}{u^*} \\
\partial D_n \arrow{r}{g_2} \arrow{u} \arrow{d} \& \partial \underline{p} \arrow[dotted]{d} \arrow[dotted]{u} \arrow{r}{w^*} \& Y^* \\
D_n \arrow{u}{w} \arrow{r}{g_3} \& \underline{p} \arrow{ur}{v^*} 
\end{tikzcd}\]

The dotted arrow exists because of the lifting property of generic morphisms, and they are automatically polygraphic due to lemma \ref{Lem_2outof3_for_polygraphic}. One can then form the pushout:

\[X' = X \coprod_{\partial D_n} D_n \overset{(g_1,g_2,g_3)}{\rightarrow} P^* \coprod_{\partial \underline{p}} \underline{p} \rightarrow Y^*   \]

where the first map is automatically generic and the second map polygraphic. This proves the result for all finite polygraphs, but one can easily pass to infinite one by taking directed colimits.

}
}

\block{\label{rk:Corres_generic-functorToTheTerminal}Let $f:A^* \rightarrow B^*$ be a generic morphism for $A$ and $B$ positive polygraphs, then there is a unique polygraphic morphisms from $B$ to the terminal positive polygraph $\Pb^+ 1$. This gives us a morphism $A^* \rightarrow (\Pb^+ 1)^*$, such that $A^* \rightarrow B^* \rightarrow (\Pb^+ 1)^*$ is its generic factorization. In particular $A^* \rightarrow B^*$ is uniquely determined (up to unique isomorphism) by the map $A^* \rightarrow (\Pb^+ 1)^*$, and conversely the generic factorization of any such maps gives back a $B$. This shows that there is a correspondence between generic maps $A^* \rightarrow B^*$ (up to unique isomorphisms under $A$) and morphisms $A^* \rightarrow (\Pb^+ 1 )^*$ which extend the correspondence between polyplexes and arrows of $(\Pb^+ 1)^*$.}

\subsection{Regular polygraphs, definitions and first properties}
\label{subsec:regularPol}

\block{Let $\underline{x}$ be a polyplex, and let $\pi^{\epsilon}_k \underline{x}$ be the polyplexes representing the arrows $\pi^{\epsilon}_k x$. They naturally assemble into a co-globular object:

\[
\begin{tikzcd}[ampersand replacement=\&]
 \pi^-_0 \underline{x}  \arrow{r} \arrow{ddr} \&  \pi^-_1 \underline{x}  \arrow{r} \arrow{ddr} \& \dots \arrow{r} \arrow{ddr} \&  \pi^-_{n-1} \underline{x} \arrow{dr} \& \\
             \&             \&       \&                 \& \underline{x} \\
\pi^+_0  \underline{x} \arrow{r} \arrow{uur}\& \pi^+_1 \underline{x} \arrow{r} \arrow{uur} \& \dots \arrow{r} \arrow{uur} \& \pi^+_{n-1} \underline{x} \arrow{ur} \& \\
\end{tikzcd}
\]

\Def{A polyplex $\underline{x}$ is said to have spherical boundary if the co-globular object above is ``Reedy cofibrant'', i.e. if for all $k\geqslant 0$, the maps:

\[ \pi^-_{k} \underline{x} \coprod_{\pi^-_{k-1} \underline{x} \coprod \pi^+_{k-1} \underline{x} } \pi^+_{k} \underline{x} \hookrightarrow \pi^-_{k+1} \underline{x} \]

\[ \pi^-_{k} \underline{x} \coprod_{ \pi^-_{k-1} \underline{x} \coprod \pi^+_{k-1} \underline{x} } \pi^+_{k} \underline{x} \hookrightarrow  \pi^+_{k+1} \underline{x} \]

are monomorphisms of polygraphs, with the convention that $ \pi^-_{-1} \underline{x} = \pi^+_{-1} \underline{x} = \emptyset$.
}

Note that if $v$ is a plex whose source and target are polyplexes with spherical boundary then, it also has spherical boundary. If $v$ is a polyplex with spherical boundary then its source and target are polyplexes with spherical boundary.

}

\block{\Def{A positive polygraph is said to be regular if all its plexes have spherical boundary.}

Regular polygraphs form a class of polygraphs in the sense of \cite{henry2017positive}. Indeed, the terminal regular polygraph is the sub-polygraph of the terminal polygraph whose plexes are the regular plexes (i.e. the plexes which have spherical boundary and such that all the lowed dimensional plexes appearing in them are also with spherical boundary) and a polygraph is regular if and only if its maps to the terminal polygraph factor into this subobject. 

A polyplex or a plex are said to be regular if their underlying polygraphs are regular.

Note that a plex can have a spherical boundary without being regular (if it contains plexes of lower dimension that does not have spherical boundary). And a polyplex can be regular without having itself a spherical boundary (for example any horizontal composition of two regular polyplexes).

}

\block{\label{prop_regular_mono}\Prop{\begin{itemize}

\item If $\underline{x}$ is a regular polyplex then the natural maps:
 \[ \pi^-_k \underline{x} \hookrightarrow \underline{x} \qquad \pi^+_k \underline{x} \hookrightarrow \underline{x} \]
are monomorphisms for all $k$.

\item If $\underline{x}$ and $\underline{y}$ are composable regular polyplex then the natural\footnote{They are for example the maps that comes from the explicit description of how composition of polyplex is defined, cf \ref{discuss:InductiveConstructionOfPolyplex}} maps:
\[  \underline{x} \hookrightarrow \underline{x} \#_k \underline{y} \qquad \underline{y} \hookrightarrow \underline{x} \#_k \underline{y} \]
are monomorphisms.

\end{itemize}
}

\Dem{We will prove the first property by induction on polyplexes (see \ref{discuss:polygraphic_induction}). If $\underline{x}$ is an $n$-plex then it has spherical boundary, hence the proposition holds for $k<n$ by definition of ``spherical boundary'' and for $n \leqslant k$ the maps $\pi^{\epsilon}_k \underline{x} \rightarrow \underline{x}$ are isomorphisms, in particular monomorphisms. If $\underline{x}=\underline{y} \#_n \underline{z}$ is a composite of two polyplexes $\underline{y}$ and $\underline{z}$ which satisfies the proposition (for all $k$), then:

\[ \underline{x} = \underline{y} \coprod_{\pi^+_n \underline{y}} \underline{z} \]

with $\pi^+_k \underline{y} \hookrightarrow \underline{y} $ and $\pi^-_k \underline{z} \hookrightarrow \underline{z}$ being monomorphisms.

In particular, for $k \leqslant n$, $\pi^-_k \underline{x}=\pi^-_k \underline{y} $ with:

\[ \pi^-_k \underline{y} \hookrightarrow \underline{y} \hookrightarrow \underline{x} \]

and $\pi^+_k \underline{x}=\pi^+_k \underline{z}$ with:

\[ \pi^+_k \underline{z} \hookrightarrow \underline{z} \hookrightarrow \underline{x} \]

for $k>n$ one has

\[ \pi^-_k \underline{x} = \pi^-_k \underline{y} \coprod_{\pi^+_n \underline{y}} \pi^-_k \underline{z} \]

as one has a diagram of monomorphisms:

\[
\begin{tikzcd}[ampersand replacement=\&]
\underline{y} \&  \pi^+_n \underline{y} \simeq \pi^-_n \underline{z} \arrow[hook]{r} \arrow[hook]{l} \& \underline{z} \\
\pi^-_k \underline{y} \arrow[hook]{u} \& \pi^+_n \underline{y} \simeq \pi^-_n \underline{z} \arrow[<->]{u}{=} \arrow[hook]{r} \arrow[hook]{l} \& \pi^-_k \underline{z} \arrow[hook]{u} \\
\end{tikzcd}
\]

Hence as the coproducts are taken in a presheaves category (hence objectwise in sets) on has indeed:

\[ \pi^-_k \underline{x} = \pi^-_k \underline{y} \coprod_{\pi^+_n \underline{y}} \pi^-_k \underline{z} \hookrightarrow \underline{y} \coprod_{\pi^+_n \underline{y}} \underline{z} = \underline{x}\]

the proof for $\pi^+_k \underline{x}$ is exactly the same and this concludes the proof of the first point.

The second properties follows immediately from the first because of the formula for the composition of polyplexes:

 \[ \underline{x} \#_k \underline{y} = \underline{x} \coprod_{\pi^+_k \underline{x}} \underline{y} \]
As the maps $\pi_k^+ \underline{x} \hookrightarrow \underline{x}$ and $\pi_k^+ \underline{x} = \pi_k^- \underline{y} \hookrightarrow \underline{y}$ are monomorphisms, and as pushout of monomorphisms are monomorphisms, this implies that the two maps $\underline{x} \hookrightarrow  \underline{x} \#_k \underline{y}$ and $\underline{y} \hookrightarrow  \underline{x} \#_k \underline{y}$ are indeed monomorphisms.

}

}

\block{\label{Cor:check_regular_with_inter}\Cor{For a regular $n$-polyplex $\underline{x}$, the following conditions are equivalent:

\begin{itemize}

\item \[  \pi^-_{k} \underline{x} \coprod_{ \pi^-_{k-1} \underline{x}  \coprod \pi^+_{k-1} \underline{x} }  \pi^+_{k} \underline{x}  \hookrightarrow  \pi^-_{k+1} \underline{x} \]

is a monomorphism.

\item \[  \pi^-_{k} \underline{x} \coprod_{ \pi^-_{k-1} \underline{x}  \coprod \pi^+_{k-1} \underline{x}} \pi^+_{k} \underline{x}  \hookrightarrow  \pi^+_{k+1} \underline{x} \]

is a monomorphism.

\item As sub-polygraphs of $\underline{x}$: $\pi^-_{k} \underline{x} \cap \pi^+_k (\underline{x}) = \pi^-_{k-1} \underline{x} \cup \pi^+_{k-1} \underline{x}$

\end{itemize}

In particular, it is enough to check the last condition for all $k<n$ to decide whether $\underline{x}$ has spherical boundary.

}

\Dem{This comes from the following lemma: in a presheaf category\footnote{This is actually true in any regular category.} if $A$ and $B$ are subobjects of $X$ and $C \rightarrow A \cap B$ is any morphism, then the map $A \coprod_C B \rightarrow X$ is a monomorphism if and only if the map $C \rightarrow A \cap B$ is an epimorphism.

Indeed, as in a presheaf category epimorphisms, monomorphisms, intersection, pushout and isomorphism are all constructed/characterized objectwise in the category of sets, it is enough to check this claim in the category of sets, where it is a completely trivial fact.

The category of regular polygraphs is a presheaf category hence this observation applies to $A=\pi^-_k \underline{x}$, $B=\pi^+_k \underline{x}$, $C=\pi^+_{k-1} \underline{x} \coprod \pi^-_{k-1} \underline{x}$ and $X = \pi^-_{k+1} \underline{x}$ (resp. $X = \pi^+_{k+1} \underline{x}$), which are all subobjects of $\underline{x}$ because of proposition \ref{prop_regular_mono}. This immediately gives the equivalence of the first point (resp. the second) with the third (note that $\pi^-_{k-1} \underline{x} \cup \pi^+_{k-1} \underline{x}$ is exactly the image of $\pi^-_{k-1} \underline{x} \coprod \pi^+_{k-1} \underline{x}$ in $\underline{x}$).
}
}

\block{\label{Prop_regularplex_mono}\Prop{If $\underline{p}$ is a regular polyplex and $\underline{a}$ is a plex then any polygraphic morphism $\underline{a} \rightarrow \underline{p}$ is a monomorphism.}

\Dem{We fix a $k$-plex $\underline{a}$. We will prove by induction on $\underline{p}$ that any map $\underline{a} \rightarrow \underline{p}$ is a monomorphism. As there is no such maps if the dimension of $\underline{p}$ is smaller than the dimension of $\underline{a}$, our induction start in dimension $k$.

\begin{itemize}

\item if $\underline{p}$ and $\underline{a}$ are plexes of the same dimension $k$ then any map $\underline{a} \rightarrow \underline{p}$ is an isomorphism.

\item if $\underline{p}$ is a plex of dimension $n>k$ and the hypothesis holds for any polyplex of dimension $n-1$, then as any map $\underline{a} \rightarrow \underline{p}$ factors either in the source or the target of $\underline{p}$, one has by assumption a monomorphisms $\underline{a} \hookrightarrow \pi^{\epsilon}_k \underline{p} \hookrightarrow \underline{p}$ where the second maps is a monomorphism because of proposition \ref{prop_regular_mono}.

\item If $\underline{p} = \underline{x} \#_k \underline{y}$ with $\underline{x}$ and $\underline{y}$ polyplexes satisfying the induction hypothesis, then as $\underline{a}$ is a plex, any map $\underline{a} \rightarrow \underline{p}$ factor either into $\underline{x}$ or $\underline{y}$, then one has $\underline{a} \hookrightarrow \underline{x} \hookrightarrow \underline{p}$ or $\underline{a} \hookrightarrow \underline{y} \hookrightarrow \underline{p}$ where each time the first map is a monomorphism by the induction hypothesis and the second map is a monomorphisms by proposition \ref{prop_regular_mono}.

\end{itemize}

This concludes the proof.
}
}

\block{\Def{Let $\underline{p}$ by a polyplex, a cell $x$ of $\underline{p}$ is said to be an \emph{inner} cell if it does not belong to the source or the target of $\underline{p}$.}}

\block{\label{Prop_unique_inner_cell}\Prop{Let $\underline{p}$ be a regular $n$-polyplex. Let $x$ be a $k$-cell of $\underline{p}$ for $k<n$. Then either:

\begin{itemize}

\item $x$ is in the source of $\underline{p}$.

\item or $x$ is an inner cell of the target of an $n$-plex of $\underline{p}$.
\end{itemize}

Moreover the two possibilities are exclusive and in the second case the $n$-plex is unique.

}

\Dem{We will prove the result by induction on $\underline{p}$.

 If $\underline{p}$ is a plex, then $x$ is either in the source of $\underline{p}$ or\footnote{this also holds constructively as polyplexes are finite and decidable} not in the source of $\underline{p}$. If $x$ is not in the source of $\underline{p}$, then as it is not the top cell for dimension reason, it has to be in the target of $\underline{p}$, moreover, as the source and target of the target of $\underline{p}$ are included in the source of $\underline{p}$, $x$ is an inner cell of the target of $\underline{p}$, and this two possibilities are indeed exclusive.

\bigskip

Let $\underline{p} = \underline{q} \#_k \underline{r} = \underline{q} \coprod_{\pi^+_k \underline{q}} \underline{r} $ with $\underline{q}$ and $\underline{r}$ satisfying the proposition. Note first that because of proposition \ref{prop_regular_mono} the maps $\underline{q},\underline{r} \rightrightarrows \underline{p}$ are monomorphisms, in particular, being an inner cell of the target of an $n$-plex in $\underline{q}$ or $\underline{r}$ is the same as being an inner cell of the target of the same $n$-plex in $\underline{p}$. Moreover, still using proposition \ref{prop_regular_mono} and a case disjunction, if a cell $x$ is in the source of $\underline{p}$ and belongs to $\underline{q}$ (resp. $\underline{r}$) then it is also in the source of $\underline{q}$ (resp. $\underline{r}$).

We prove that at least one of the two case is satisfied, we need distinguishes several possibilities: 

\begin{itemize}

\item $x$ is in $\underline{q}$, then either $x$ is in the source of $\underline{q}$ in which case it is in the source of $\underline{p}$ or $x$ is an inner cell of a $n$-plex of $\underline{q}$, hence of $\underline{p}$.

\item $x$ is in $\underline{r}$ and $k=n$, then either $x$ is in the source of $\underline{p}$, in which case it is also in $\underline{q}$ and this is the previous case, or it is an inner cell of the target of an $n$-plex in $\underline{r}$ and hence also in $\underline{p}$.

\item $x$ is in $\underline{r}$ and $k<n$, this is the same as above, except that now if $x$ is in the source of $\underline{r}$, then $x$ is directly in the source of $\underline{p}$.

\end{itemize}

We now check the exclusivity and uniqueness property:

If $x$ is an inner cell of two different $n$-plexes of $\underline{p}$, then these two $n$-plexes cannot be both in $\underline{q}$ or in $\underline{r}$ because of the induction hypothesis, hence one is in $\underline{q}$ and the other is in $\underline{r}$, in particular $x$ is both in $\underline{q}$ and $\underline{r}$, but cannot be in the source of $\underline{r}$, which is impossible, as $\underline{p}$ is obtained by gluing $\underline{q}$ and $\underline{r}$ together along the source of $\underline{r}$.

If $x$ is both in the source of $\underline{p}$ and an inner cell of the target of an $n$-cell $v$ of $\underline{p}$. Then $v$ is either in $\underline{q}$ or $\underline{r}$, and $x$ automatically belongs to the same component. As $x$ is in the source of $\underline{p}$ is has to be in the source of the component in which it is, but this contradicts the induction hypothesis.

}
}

\block{\label{Prop_Charsphericalboundarybytotality}\Prop{Let $\underline{p}$ be a regular $n$-polyplex. Then the following two conditions are equivalent:

\begin{itemize}

\item The source of $\underline{p}$ has spherical boundary.

\item The target of $\underline{p}$ has spherical boundary.

\end{itemize}

Moreover, assuming these conditions hold, the following conditions are also equivalent:

\begin{itemize}

\item $\underline{p}$ has spherical boundary.

\item Every inner cell of the source of $\underline{p}$ appears as an inner cell of the source of an $n$-plex in $\underline{p}$.

\end{itemize}

}

\Dem{The fact that the source of $\underline{p}$ has spherical boundary is equivalent to the fact that all the maps:

\[ \pi^-_{k} \underline{p}  \coprod_{\pi^-_{k-1} \underline{p}  \coprod \pi^+_{k-1} \underline{p}} \pi^+_k \underline{p} \hookrightarrow \pi^-_{k+1} \underline{p} \]
 for $k < n-2$ and

\[ \pi^-_{k} \underline{p} \coprod_{\pi^-_{k-1} \underline{p}  \coprod \pi^+_{k-1} \underline{p}} \pi^+_k \underline{p} \hookrightarrow \pi^+_{k+1} \underline{p} \]

for $k<n-3$ are monomorphisms. So in order to prove that the target of $\underline{p}$ also have spherical boundary one only needs to prove that the second map is also a monomorphism for $k=n-2$.

As $\underline{p}$ is regular, $\pi^-_{n-1} \underline{p} \rightarrow \underline{p}$ is a monomorphism, hence if the source of $\underline{p}$ has spherical boundary one has a monomorphism:

\[  \pi^-_{n-2} \underline{p}  \coprod_{\pi^-_{n-3} \underline{p}  \coprod \pi^+_{n-3} \underline{p}} \pi^+_{n-2} \underline{p} \hookrightarrow \pi^-_{n-1} \underline{p} \]

which fits into a commutative diagram:
\[
\begin{tikzcd}[ampersand replacement=\&]
\displaystyle \pi^-_{n-2} \underline{p}  \coprod_{\pi^-_{n-3} \underline{p}  \coprod \pi^+_{n-3} \underline{p}} \pi^+_{n-2} \underline{p} \arrow[hook]{r} \arrow{d} \& \pi^-_{n-1} \underline{p} \arrow[hook]{d} \\
\pi^+_{n-1} \underline{p} \arrow{r} \& \underline{p} \\
\end{tikzcd}
\]

Hence the left map is always a monomorphism. Which proves that the target of $\underline{p}$ also have spherical boundary. The reverse implication is obtained in the exact same way.

For the second equivalence, if $\underline{p}$ has spherical boundary, then every inner cell of the source of $\underline{p}$ is not in the target of $\underline{p}$, Hence it must be an inner cell of the source of a $n$-plex of $\underline{p}$ by the dual of proposition \ref{Prop_unique_inner_cell}.

Conversely, assume that the second condition holds. And consider $x$ and $y$ to be two cells of 

\[  \pi^-_{n-1} \underline{p}  \coprod_{\pi^-_{n-2} \underline{p}  \coprod \pi^+_{n-2} \underline{p}} \pi^+_{n-1} \underline{p} \]

which have the same image in $\underline{p}$. If $x$ and $y$ are both in $\pi^-_{n-1} \underline{p}$ or both in $\pi^+_{n-1} \underline{p}$ then $x=y$ by proposition \ref{prop_regular_mono}. So we assume that $x$ is in $\pi^-_{n-1} \underline{p}$ and $y$ is in $\pi^+_{n-1} \underline{p}$.

\begin{itemize}

\item if $x$ is in the boundary of $\pi^-_{n-1} \underline{p}$ then it is also in $\pi^+_{n-1} \underline{p}$ and we are brought back to the previous case.

\item if $x$ is an inner cell of $\pi^-_{n-1} \underline{p}$, then by assumption, $x$ is an inner cell of the source of some $n$-plex of $\underline{p}$, hence it cannot belong to the target of $f$, hence $x$ and $y$ cannot be equal in $f$ which contradicts our assumption.

\end{itemize}

In all cases, one can conclude that $x=y$ in $\displaystyle \pi^-_{n-1} \underline{p}  \coprod_{\pi^-_{n-2} \underline{p}  \coprod \pi^+_{n-2} \underline{p}} \pi^+_{n-1} \underline{p} $ which proves that the map to $\underline{p}$ is a monomorphisms and hence that $\underline{p}$ has spherical boundary.
}
}

\block{\label{Cor_nplex_cover}\Cor{In a regular $n$-polyplex $\underline{p}$ with spherical boundary, every cell $x$ of $\underline{p}$ factors into the inclusion of $y \rightarrow \underline{p}$ where $y$ is an $n$-plex.}

\Dem{We will prove it by descending induction on the dimension of $x$. If $x$ has dimension $n$ then $x$ itself is the inclusion of an $n$-plex. Assume that $x$ has dimension $k$ and that the result is known for all cells $x$ of dimension greater then $k$. Let $u$ be the smallest integer such that $x$ is in $\pi^-_u \underline{p}$ or $\pi^+_u \underline{p}$, one has that $k \leqslant u \leqslant n$, moreover $x$ is an inner cell of  $\pi^-_u \underline{p}$ or of $\pi^+_u \underline{p}$ (otherwise it would belong to $\pi^-_{u-1} \underline{p}$ or $\pi^+_{u-1} \underline{p}$). We need to treat two cases:

\begin{itemize}
\item If $u=n$, then by proposition \ref{Prop_unique_inner_cell}, as $x$ is not in the source of $\underline{p}$, $x$ is an inner cell of the target of a unique $n$-plex, in particular it factor into an $n$-plex of $\underline{p}$.

\item If $u<n$, then as $y$ has spherical boundary, $\pi^-_{u+1} \underline{p}$ has spherical boundary as well, and proposition \ref{Prop_Charsphericalboundarybytotality} implies that $x$ must be an inner cell of the source of a (unique) $u+1$-plex of $\pi^-_{u+1} \underline{p}$. By induction, as $u+1 >k$, this unique $u+1$-plex factor into an $n$-plex of $Y$ and hence this concludes the proof.

\end{itemize}

}
}

\blockn{We conclude this section with the following lemma that will only be useful during the proof of one of our key lemma (\ref{Lemma_Key_weakequivalence}) but which fits more naturally in the present section.}

\block{\label{Lemma_Stuff}\Lem{Let $\underline{p}$ be a regular $n$-polyplex which has a single $n$-cell $x$. Let $\underline{p'}$ be the polygraph obtained from $\underline{p}$ by: removing $x$, removing all the inner cells of $\pi_{n-1}^{+} x$ and $\pi_{n-1}^{-}x$ and adding a single $n-1$-cell between $\pi_{n-2}^- x$ and $\pi_{n-2}^+ x$.

Then $\underline{p'}$ is indeed a polygraph, and it is in fact a regular $n-1$-polyplex satisfying $\pi_{n-2}^{\epsilon} \underline{p'} = \pi_{n-2}^{\epsilon} \underline{p}$.

}

\Dem{We proceed by induction on the construction of the $n$-plex $\underline{p}$. If $\underline{p}$ is a plex, then it has spherical boundary, hence removing its top cell $p$ as well as all the inner cell of $\pi^{-}_{n-1} p$ and $\pi_{n-1}^{+} p$ leaves us with $\pi_{n-2}^- \underline{p}$ and $\pi_{n-2}^+ \underline{p}$ glued together along their common $(n-3)$-dimensional boundary. Adding a single $(n-1)$-cell between $\pi_{n-2}^- p$ and $\pi_{n-2}^+ p$ gives us a $(n-1)$-plex with the correct boundary, and it is regular as $\pi_{n-2}^{\epsilon} p$ are regular with spherical boundary.

If $\underline{p}$ is a composite $\underline{p} = \underline{q} \#_k \underline{r}$ then as $\underline{p}$ have to contain a single $n$-cell either $\underline{q}$ or $\underline{r}$ are of dimension $<n$. Hence one can freely assume that $k<n-1$ otherwise one already has $\underline{p}=\underline{q}$ or $\underline{p}=\underline{r}$. One assumes that, for example, the unique $n$-cell $x$ of $\underline{p}$ is in $\underline{q}$. One can then form $\underline{q'}$ inductively, and as $\pi_k^{\epsilon} \underline{q'} = \pi_k^{\epsilon} \underline{q}$ for $k<n-1$, the composite $\underline{q'} \#_k \underline{r}$ exists and is a $n-1$-polyplex. But applying the construction described in the lemma to $\underline{p}$ exactly gives the polygraph $\underline{p'}=\underline{q'} \#_k \underline{r}$.
}
}

\section{Algebraic classes of polygraphs}
\label{section_algebraic_class}

\blockn{Specifically for this subsection, it would have been more convenient to use to the ``globular'' definition of $\infty$-categories. I.e. when we talk of an arrow of an $\infty$-category it has an assigned dimension and an arrow is different from its identity arrow. One will keep the definition given in \ref{Def_StreetInfinityCat}, but for the time being whenever we say ``arrow'' of an $\infty$-category, we mean a pair $(x,n)$ where $n$ is an integer and $x$ is an $n$-arrow. The identity of $(x,n)$ is $(x,n+1)$ the $k$-source and $k$-target are $(\pi_k^{\epsilon} x ,k)$, etc...}

\blockn{The notion of class of polygraphs was introduced in \cite{henry2017positive}. Informally, a class of polygraphs is a full subcategory of the category of polygraphs whose objects are determined by a condition on the shape of their cells. More precisely, it is a full subcategory of the category of polygraphs of the form $\Pb/X$ where $X$ is a sub-polygraph of the terminal polygraph. Again informally a cell of $\Pb 1$ corresponds to a possible shape of a cell in a polygraph and $X$ select the shape of cells that are admissible for our class of polygraphs. The typical example are the class of all polygraphs, of positive polygraphs, of polygraph of dimension smaller than $n$, of opetopic or ``many-to-one'' polygraphs, and of regular polygraphs. 

If $\C$ is a class of polygraphs, its terminal object (denoted $X$ above) is denoted $\C 1$.
}

\blockn{The idea of ``algebraic class of polygraph'' is the following: one has a very close connection between the notion of ``non-unital $\infty$-category'' that we introduced in \cite[A.5]{henry2017positive} and the notion of positive polygraphs. Essentially the non-unital polygraph are exactly the ``polygraphs'' for non-unital $\infty$-categories. Algebraic classes of polygraphs are in general supposed to be the class of polygraphs that are obtained this way as ``polygraphs for a weakened notion of $\infty$-category''. The main example we are interested in are positive polygraphs and regular polygraphs, but in the case of regular polygraphs it is not completely clear what should be the corresponding notion of ``regular $\infty$-category'': they should be some kind of $\infty$-categories where only ``regular'' compositions are defined (instead of all compositions) but the notion is considerably harder to formalize than non-unital $\infty$-categories. So instead we will give a different definition of what an algebraic class should be and explain how to construct the corresponding notion of $\infty$-category. We will not prove at this point that regular polygraphs form an algebraic class of polygraphs: the goal of this section is only to clarify what exactly needs to be proved to show that the notion of regular $\infty$-category make sense. The key ideas behind the definition of $C$-categories will be explained in \ref{discuss:DefOfCcategoriesMotivation}.}

\blockn{In subsection \ref{subsec_tensorstable_algebraic_class}, we will show that if an algebraic class $\C$ satisfies an additional condition of stability under the Gray tensor product (The Gray tensor will be introduced properly only in section \ref{subsec_Gray_tensor_product}) then one can also define a Gray tensor product on the category of $\C$-categories as well. Here again our goal in this section is mostly to understand what is needed in order to have a tensor product on the category of regular $\infty$-categories, and the appropriate condition on the class of regular polygraphs will be proved later (in \ref{Prop_tensor_product_of_regular}).
}

\subsection{Definition and $C$-categories}
\label{subsec_algebraicclass_def}

\block{\Def{Given a class of polygraphs $\C$, and $\C 1 \subset \Pb 1$ its terminal object

\begin{itemize}

\item An arrow of $(\Pb 1)^*$ is said to be $\C$-admissible if it is the source or the target of a cell of $\C 1$ (in particular, it is in $\C 1$).

\item The class $\C$ is said to be pre-algebraic if every cell of $\Pb 1$ whose source and target are $\C$-admissible is in $\C$.

\item An arrow of a $C$-polygraph is said to be $\C$-admissible if its image in $\C 1$ is $\C$-admissible.

\item If $X$ and $Y$ are two $C$-polygraphs, a $\C$-morphism from $X$ to $Y$ is functor $f:X^* \rightarrow Y^*$ which send every generating cell of $X$ to a $\C$-admissible arrow in $Y$. 

\item The class $\C$ is said to be algebraic if it is pre-algebraic, that every generator of $\C1$ is $\C$-admissible, and if the image of a $\C$-admissible arrow by a $\C$-morphism is $\C$-admissible.

\end{itemize}

}

The intuitive idea is as follow: A pre-algebraic class is a class such that one has a certain class of ``nice arrows'' (the $\C$-admissible arrow) such that the cells that can be added to a $\C$-polygraph are exactly those whose both source and target are $\C$-admissible. For example, for the class of regular polygraphs, the plexes that can be added are exactly those whose source and target are regular with spherical boundary. For the class of positive polygraphs, the type of plex that one can add are those whose source and target are non-identity. Hence both the class of regular and positive polygraph are pre-algebraic. The class of opetopic polygraphs (or many-to-one polygraphs), i.e. the class of polygraph such that the target of each generator is a generator is not pre-algebraic as the condition is not symmetric on the source and the target.

An algebraic class is one where moreover $\C$-morphisms form a category. Informally it means that a ``$\C$-admissible composition of $\C$-admissible arrow is a  $\C$-admissible arrow'', where a $\C$-admissible composition mean something like\footnote{The description we are giving here is formally correct is the case of a ``good class'' in the sense of \cite{henry2017positive}, so for example if it is a sub-class of the class of positive polygraphs.} a composition indexed by a polyplex such that the corresponding arrow in $\C 1$ is $\C$-admissible. For example, for a positive polygraphs, it means that taking a composition of some non-identity arrows along a diagram of dimension $n$ with at least one $n$-arrow in the composition give a non-identity arrow of dimension $n$.

Here are our examples of such algebraic classes:

\begin{itemize}

\item The class $\Pb$ of all polygraphs: every arrow is $\Pb$-admissible and $\Pb$-morphisms are just morphisms of $\infty$-categories.

\item The class $\Pb^+$ of positive polygraphs: $\Pb^+$-admissible arrow are non-identity arrow, $\Pb^+$-morphisms corresponds to the morphisms of non-unital $\infty$-category we introduced in the appendix of \cite{henry2017positive} and it is relatively easy to check that it is indeed an algebraic class. It is also not very hard too see that it is the largest good (pre)algebraic class of polygraphs.

\item The class of regular polygraphs: One easily see that it is a pre-algebraic class with admissible arrow being the arrows whose representing polyplex have a spherical boundary. We will see in \ref{Cor:RegularPolygraphAreAnAlgebraicClass} that it is an algebraic class of polygraphs.

\item The class $\G$ of globular polygraphs (those corresponding to globular sets). The $\G$-admissible arrows are just the generators and the $\G$-morphisms are just morphisms of globular sets. It is the smallest algebraic class of polygraphs.

\end{itemize}

}

\block{\label{discuss:DefOfCcategoriesMotivation}We can now define the idea behind the definition of ``$C$-categories''. Even if we do not really know what these are supposed to be, as the corresponding notion of polygraphs is supposed to be the category of $C$-polygraph one expect $C$-categories to be monadic over $C$-polygraphs, so we only need to understand what is the corresponding monad on $C$-polygraphs. It is sufficient to understand the Kleisli category of that monad, i.e. it is sufficient to understand what are the morphisms of $C$-categories between two free $C$-categories on $C$-polygraphs. If $X$ is a $C$-polygraphs, one expect the $C$-category freely generated by $X$ to consist exactly in the $C$-admissible arrow of $X^*$, because these are those on which we are going to be able to glue new cells if one tries to reconstruct the notion of $C$-polygraph out of the notion of $C$-category. Hence it is natural to ask that the morphisms between these free $C$-categories should be exactly the $C$-morphisms, and the condition of being an algebraic class of polygraphs is exactly what we need to makes this into a category.}

\block{\label{prop:theMcMonadExist}\Prop{Let $\C$ be an algebraic class of polygraphs, then there is a unique monad $M_{\C}$ on $\C$ such that the category $\widetilde{\C}$ of $C$-polygraphs and $\C$-morphisms between them is the Kleisli category of $M_{\C}$.

The monad $M_{\C}$ commutes to directed colimits.
}

\Dem{ First observes that the category $\C$ (with polygraphic morphisms) is a locally finitely presentable category: the finite polygraphs in $\C$ are finitely presentable objects stable under finite colimits, and every object of $\C$ is canonically the directed colimits of all morphisms from finite polygraphs to it.

The inclusion functor $\C \rightarrow \widetilde{\C}$ commutes to all colimits:  If $\widetilde{\C}$ was defined using all morphisms of $\infty$-categories instead of just $\C$-morphisms, then this would follow from the fact that $X \mapsto X^*$ commutes to all colimits. Now let $X_i$ be some diagram in $\C$, the set of cells of $(\text{colim } X_i )$ is the colimit of set of cell of the $X_i$, hence a morphism of $\infty$-category $\text{colim } X_i^* \rightarrow Y$ is a $\C$-morphism if and only if its restriction at each $X_i$ is a $\C$-morphisms. Hence it is exactly the data of a compatible family of $\C$-morphisms $X_i \rightarrow Y$ which proves the claim.

The special adjoint functor theorem implies that this functor as an adjoint. As it is bijective on objects, this implies that $\widetilde{\C}$ is the Kleisli category of the monad $M_{\C}$-induced by this adjunction. The uniqueness is immediate from the Yoneda lemma.

For the commutation to direct colimit of $M_{\C}$: If $(A_i)$ is a directed diagram of polygraphs with colimits $A$ and $X$ is a finite polygraph then one easily see that the colimit of $Hom(X, M_{\C}(A_i)) = Hom_{\widetilde{\C}}(X,A_i)$ is $Hom_{\widetilde{\C}}(X,A) = Hom(X, M_{\C}(A))$: if $M_{\C}$ was defined using all morphisms of $\infty$-categories that would be just the fact that the free $\infty$-category functor is finitary, and because the maps $A_i \rightarrow A$ are polygraphic, a morphism $X^* \rightarrow A_i^*$ is a $\C$-morphisms if and only if its composite $X^* \rightarrow A^*$ is a $\C$-morphisms. Finally because the inclusion of $\C$ in the category of presheaves over $\C^f$ is fully faithful this is enough to conclude that the colimits of the $M_{\C}(A_i)$ is $M_{\C}(A)$.
}
}

\block{\Def{For $\C$ an algebraic class of polygraphs one defines the category of $\C$-categories as the category of algebras for the monad $M_{\C}$.}

Here are some examples:

\begin{itemize}

\item If $\Pb$ is the class of all polygraphs, then the $\Pb$-categories are the same as $\infty$-category. Indeed it was show by F.Metayer in \cite{metayer2016strict} that the category of $\infty$-category is monadic over the category of polygraphs and as the $\Pb$-morphisms are all morphisms of $\infty$-category the monad $M_{\Pb}$ on $\Pb$ is indeed the monad induce by the adjunction between polygraphs and $\infty$-category.

\item If $\G$ is the class of globular polygraphs, the category of $\G$-categories is the category of globular sets. Indeed, $\G$-morphisms are the same as morphisms of globular sets, so the monad $M_{\G}$ is the identity monad.

\item It is not very hard to check that $\Pb^+$-categories are the ``non-unital categories'' defined in \cite{henry2017positive}: the functor from positive polygraphs to non-unital $\infty$-categories clearly induces the monad $M_{\Pb^+}$ defined above, so the only things to be check is that the category of non-unital $\infty$-categories is indeed monadic over positive polygraphs. This is not very difficult to do, in fact, it follows from theorem $5.1$ of M.Batanin \cite{batanin1998computads}, whose proof is incorrect in general as it assume that polygraphs always form a presheaves category, but is correct in this case as positive polygraphs do form a presheaves category. 

\end{itemize}

}

\block{\label{Prop_Nervetheorem}\Prop{Let $\C$ be an algebraic class of polygraphs, and let $\C^f$ be the full sub-category of finite polygraphs. Then the monad $M_{\C}$ is a monad with arity $\C^f$ in the sense of \cite{berger2012monads}, or a ``Nervous Monad'' in the sense of \cite{bourke2018monads}.

In particular, if $\widetilde{\C^f}$ denotes the category of finite polygraphs in $\C$ with $\C$-morphisms between them, then the category of $\C$-categories identifies with the full subcategory of presheaves on $\widetilde{\C^f}$ whose restriction of $\C^f$ commutes to finite colimits.}

\Dem{We already mentioned that $\C^f$ was a dense generator of $\C$. So following definition $1.8$ of \cite{berger2012monads} one needs to check that for every object $X \in  \C$ the monad $M_{\C}$ send the diagram given by $\C^f/X$ and its canonical cocone $X$ to a colimit diagram in the category of presheaves in $\C^f$. But this is immediate: this diagram is directed because $\C^f$ has all finite colimits, hence the colimits in the category of presheaves is a colimit in $\C$, and $M_{\C}$ commutes to directed colimits. 

The second part of the claim is then Leinster's Nerve theorem (theorem $1.10$ of \cite{berger2012monads}).

}
}

\subsection{ $\otimes$-stable algebraic class of polygraphs}

\label{subsec_tensorstable_algebraic_class}

\block{\Def{An algebraic class of polygraph is said to be $\otimes$-stable, if the Gray tensor product of two $\C$-polygraphs is a $\C$-polygraph and the tensor product of two $\C$-morphisms between $\C$-polygraphs is again a $\C$-morphism.}

See subsection \ref{subsec_Gray_tensor_product} for the definition of the Gray tensor product.

}

\blockn{Let $\C$ be a $\otimes$-stable algebraic class of polygraph. We want to define the Gray tensor product of two $\C$-categories. The simplest way to do that is to apply the following theorem of B.Day:}

\block{\label{Th_DayConv}\Th{Let $\Ccal$ be a complete and co-complete category. $\Dcal$ a full subcategory of $\Ccal$, which is endowed with a monoidal structure. One assume that:

\begin{itemize}

\item $\Dcal$ contains a dense subcategory of $\Ccal$.

\item There are functors $H,H': \Dcal^{op} \times \Ccal \rightarrow \Ccal$, and isomorphisms:

\[ Hom(D,H'(D',C)) \approx Hom(D \otimes D',C) \approx Hom(D',H(D,C) )\]

functorial in $D,D' \in \Dcal$ and $C \in \Ccal$.

\end{itemize}

Then there exist a unique (up to unique isomorphism) bi-closed monoidal structure on $\Ccal$ which makes the inclusion functor $\Dcal \hookrightarrow \Ccal$ into a monoidal functor.

}

This form of the result can be found as theorem 6.3 of \cite{ara2016joint}.
}

\block{\Prop{Let $\C$ be a $\otimes$-stable algebraic class of polygraphs. Then there is a unique bi-closed monoidal structure on the category of $\C$-categories whose restriction to the category of $\C$-polygraphs (and $\C$-morphisms between them) is the Gray tensor product.}

\Dem{We apply Day's theorem (i.e. \ref{Th_DayConv}) with $\Ccal$ the category of $\C$-categories and $\Dcal$ the category of finite $\C$-polygraph and $\C$-morphisms between them.  $\Dcal$ is a dense full subcategories essentially by definition, so we just need to check that the functors $H$ and $H'$ can be defined. But given proposition \ref{Prop_Nervetheorem}, the existence of $H$ is equivalent to the fact that for any $D$ a finite $C$-polygraph and $X$ a $\C$-category the functor:

\[Hom_{C-Cat}( D \otimes \_ , X) \]

on the category of finite $\C$-polygraphs (and polygraphic maps between them) send finite colimits to finite limits.

But the functor $D \otimes \_$ from finite $\C$-polygraphs to the category of all polygraphs commutes to all finite colimits and the functor from $\C$-polygraphs to $\C$-categories commutes to all colimits as it is a left adjoint functor. The case of $H'$ is exactly the same and this concludes the proof.
}
}

\section{Abelian groups in strict $\infty$-categories}

\label{sec_Abelian_group_object}

\subsection{Strict abelian $\infty$-group and chain complexes}
\label{subsection_StrictAbGrp_vs_chaincomplexes}

\block{\Def{A Strict abelian $\infty$-group, or $\infty$-group for short, is an abelian group object in the category of $\infty$-categories.}
 Given the description of limits in the category of $\infty$-categories explained in \ref{Discuss_limits}, an $\infty$-group is precisely given by:

\begin{itemize}

\item $(X,\pi^+_k,\pi^-_k, \#_k)$ is an $\infty$-category.

\item $X$ is endowed with a structure of an abelian group.

\item $\pi_k^{+}$ and $\pi_k^{-}$ are group morphisms for all $k$.

\item $\#_k : \{ x,y\in X^2 | \pi_k^+(x) = \pi_k^-(y) \} \rightarrow X$ is a morphism of groups, i.e.

\[ (x \#_k y ) + (z \#_k w) = (x+z) \#_k (y+w) \]

whenever the left hand side is defined.

\end{itemize}

}

\blockn{It is well known that the category of $\infty$-group is equivalent to the category of chain complexes. This was proved by R.Brown and P.J. Higgins in \cite{brown2003cubical} for various kind of $\infty$-categories, which unfortunately does not include the globular definition of $\infty$-categories but only globular $\infty$-groupoids. In this paper we will anyway need a precise understanding of the relation between the two notions in order to studies combinatorial properties of these strict $\infty$-groups in terms of chain complexes, for this reason we will review and reprove this equivalence in details in this subsection.
}

\block{\Def{A reflexive globular group\footnote{The name comes from the fact that these are the same as abelian group in the category of reflexive globular sets, but we will generally just says globular group.}, or simply a globular group, is an abelian group $G$ endowed with group endomorphisms $\pi_k^+, \pi_k^-$ for all $k \geqslant 0$ satisfying the ``globular relations'':

\[\pi_i^{\delta} \pi_j^{\epsilon} = \left\lbrace \begin{array}{c c} \pi_j^{\epsilon} & \text{ if } j \leqslant i \\ \pi_i^{\delta} & \text{ if } i<j \end{array} \right. \]

And such that $\forall x \in X, \exists i, \pi_i^+ x =x $.
}
}

\block{\label{Prop_equiv_globgroup_strictinfgroup}\Prop{The functor from the category of $\infty$-groups to globular groups which forget the composition operations $\#_k$ is an isomorphism of categories. More precisely:

\begin{itemize}

\item Every globular group admit a unique family of partially defined composition operation $x \#_k y$ making it into a strict $\infty$-category. They are given by:

 \[ x \#_k y =  x +y - \pi_k^+(x) \text{ when } \pi_k^+(x)= \pi_k^-(y) \]

\item Morphisms of globular groups are the same as morphisms of strict $\infty$-groups.

\end{itemize}
}

\Dem{If $G,\pi_k^+,\pi_k^-$ is a globular group, one easily check that $x \#_k y$ as defined in the proposition makes it into a strict $\infty$-group by looking at each axioms one by one.
Conversely, if $G$ is a strict $\infty$-group, and $x,y \in G$ are such that $\pi_k^+(x)=\pi_k^{-1}(y)$ then: 

 \begin{multline*} x + y = (x \#_k \pi_k^+ x) + ( \pi_k^+ x \#_k y) = (x + \pi_k^+ x) \#_k ( \pi_k^+ x + y) \\ = (x \#_k y) + (\pi_k^+(x) \#_k \pi_k^+(x)) = x \#_k y + \pi_k^+(x) \end{multline*}

i.e.  \[ x \#_k y = x + y - \pi_k^+(x) \]

hence the uniqueness of the structure. And any morphism of globular group clearly preserve the operation $\#_k$ as defined here.
}

}

\block{The equivalence between globular groups and chain complexes, is then a form a Dold-Kan correspondence, in the sense of \cite{lack2015combinatorial}. We review it in some details: 

Let $G$ be a globular group.

\begin{itemize}

\item One defines $\pi_i^+$ and $\pi_i^-$ for $i<0$ as being $0$. Note that the globular relations keep being true for negative values as well.

\item One defines $G_n$ to be the subgroup of $n$-arrows, i.e. the range of the projection $\pi_n^+$, or equivalently the range of $\pi_n^-$. One has $G_n=0$ for $n<0$.

\item For any $n \geqslant 0$ one defines $PG_n = G_n / G_{n-1}$.

\item As $\pi_{n-1}^-$ and $\pi_{n-1}^{+}$ are projection from $G_n$ to $G_{n-1}$ their complement allows to define a section of the projection $G_n \twoheadrightarrow PG_n$, more precisely, one defines:

\[ PG_n^{\epsilon} = \{ v \in G_n | \pi_{n-1}^{\epsilon} v = 0 \} \subset G_n \]

Both $PG^-_n$ and $PG_n^{+}$ are canonically isomorphic to $PG_n$ through the restriction of the projection $G_n \twoheadrightarrow PG_n$.
Note that $PG_0^+ = PG_0^- = G_0$ so we will sometimes write that $PG_0 \subset G$ without specifying a sign, but we will never write that $PG_n \subset G$ for $n>0$ has there are two possible and equally canonical such inclusions (so we only write $PG^+_n \subset G$ and $PG_n^{-} \subset G$).

\item The natural projection $G_n \rightarrow PG_n \simeq PG_{n}^{\epsilon} \subset G_n$ is $1- \pi_{n-1}^{\epsilon}$. Moreover $\pi_n^{\delta} - \pi_{n-1}^{\epsilon}$ is a projection of $G$ onto $PG_{n}^{\epsilon}$ for any value of $\delta$. Indeed: $\pi_n^{\delta} - \pi_{n-1}^{\epsilon} = (1- \pi_{n-1}^{\epsilon}) \pi_n^{\delta}$ is the composite of $\pi_n^{\delta} G \rightarrow G_n$ and the canonical projection $G \twoheadrightarrow PG_n \simeq PG^{\epsilon}_n$.

\item One defines a differential $\partial: PG_n \rightarrow PG_{n-1}$ which takes $x \in G_n$ to $\partial(x)=\pi_{n-1}^+(x) - \pi_{n-1}^-(x)$, and indeed pass to the quotient projection as if $x \in G_{n-1}$ then $\partial(x)=x-x=0$.

Note that $\pi_{n-1}^+$ send $PG_n^{-}$  to $ PG^{-}_{n-1}$ as $\pi_{n-2}^-\pi_{n-1}^+ = \pi_{n-2}^- = \pi_{n-2}^-\pi_{n-1}^+$. And through the identification of $PG^-_n$ and $PG_n$ this is just the maps $\partial$.

Similarly, $-\pi_n^{-} : PG_n^+ \rightarrow PG^{+}_{n-1}$ is well defined and corresponds to $\partial$.

\item It follows immediately from the globular relations that $(PG_n, \partial)$ forms a chain complex.

\item At this point, a careful reader might be puzzled by the fact that the map $\partial$ is not just defined as a map $PG_n \rightarrow PG_{n-1}$ but naturally appears as a map $PG_n \rightarrow G_{n-1}$. But no information is lost when passing to the quotient because of the following observation: $PG_n^+ \cap PG_n^{-1}$ canonically identifies through the canonical projection to $\{x \in PG_n | \partial(x)=0 \}$. Hence $\partial(x)$ always is in the subset of $PG_{n-1}$ on which the two inclusions to $G_{n-1}$ coincide. In particular $\partial(x)=0$ in $PG_{n-1}$ if and only if $\pi_{n-1}^+ x = \pi_{n-1}^- x$ in $G_{n-1}$ for any representant of $x$ in $G_n$.

\item As $PG^+_n$ and $PG^-_n$ are the kernel of projections $G_n \twoheadrightarrow G_{n-1}$ one has a canonical decomposition\footnote{We wrote ``$=$'' to signify that they are decomposition of $G_n$ and $G$ as coproduct of subspaces and not abstract isomorphisms.}:

\[ G_n = PG^-_n \oplus G_{n-1} \]
\[ G_n = PG^+_n \oplus G_{n-1} \]

iterating this one has an isomorphism:

\[ G_n = PG^-_n \oplus PG^-_{n-1} \oplus \dots \oplus PG^-_0 \]
or
\[ G_n = PG^+_n \oplus PG^+_{n-1} \oplus \dots \oplus PG^+_0 \]

or actually

\[ G_n = PG^{\epsilon_n}_n \oplus PG^{\epsilon_{n-1}}_{n-1} \oplus \dots \oplus PG^{\epsilon_0}_0 \]

for any sequence of signs $\epsilon_i$.

As $G$ is the increasing unions of the $G_n$ (and the inclusion $G_{n-1} \subset G_n$ does correspond to the canonical inclusion on the decomposition above) one has:

\[ G= \bigoplus_{n \geqslant 0} PG_n^{\epsilon_i} \]

For any sequence $\epsilon_i$ of signs.

\item One can abstractly says that $G \simeq \bigoplus PG_n$ but one needs to remember that one get a different such isomorphisms for each choice of sequence of sign (and one can actually get even other decomposition in full generality).

\item The operation $\pi_i^+$ and $\pi_i^{-}$ on $G$ can be reconstructed just from the differential $\partial$ in the decomposition above. Indeed, one knows that on $PG_k^{\epsilon}$ all the $\pi_{n}^{\delta}$ for $n\geqslant k$ are the identity, $\pi_{k-1}^{\epsilon}$ and all the $\pi_n^{\delta}$ for $n<k-1$ are zero, and $\pi_{k-1}^{-\epsilon}(x)$ is given by $-\epsilon \partial(x)$ which, as $\partial^2=0$, takes values in both $PG_{k-1}^+$ and in $PG_{k-1}^-$.

For example, in the decomposition $G = \bigoplus PG_n^{-}$, the operation $\pi_i^+$ and $\pi_i^{-}$ are given as follows:

 $\pi_i^-$ is the canonical projection on $PG^-_i \oplus \dots PG^-_0$.

 $\pi_i^+$ send all the $PG^-_n$ for $n>i+1$ to zero, it sends $PG^-_{i+1}$ to $PG^-_{i}$ by $\partial$ and it is the identity on all $PG^-_n$ for $n \leqslant i$.

\item Conversely, if $K$ is any chain complex (in positive degree and with a degree decreasing differential) then $G:= \bigoplus K_n$ with $\pi_i^+$ and $\pi_i^-$ as defined above is a globular group with $PG_n^-=K_n$.

\item These two constructions induces an equivalence between the category of chain complexes and the category of globular groups.

\item A more symmetric description of $G$ in terms of its chain complex exists (We learned it from \cite{steiner2004omega}):

To a chain complex $(K_n,\partial)$ one associate the set of sequence of pairs $(k^-_n ,k_n^+) \in (K_n)^2$ which are zero after a certain rank and which satisfies the additional condition $\partial(k_n^{\epsilon}) = k^+_{n-1} - k^-_{n-1}$.

The operation $\pi_{n}^{\epsilon}$ is defined on this group of double sequences by:

 \[ \pi_n^{\epsilon} \left( \dots \,  \mytop{k_n^-}{k_n^+} \, \mytop{k_{n-1}^-}{k_{n-1}^+} \, \dots \, \mytop{k_0^-}{k_0^+}  \right) = \left( \dots \, \mytop{0}{0} \, \mytop{0}{0} \,  \mytop{k_n^{\epsilon}}{k_n^{\epsilon}} \, \mytop{k_{n-1}^-}{k_{n-1}^+} \,  \mytop{k_{n-2}^-}{k_{n-2}^+} \, \dots \, \mytop{k_0^-}{k_0^+}  \right) \]

If $G$ is a globular group, and $K_n = PG_n$ is the corresponding chain complex, then $G$ identifies with this space of double sequence by:

\[ g \mapsto \left( \dots \,  \mytop{\overline{\pi^-_n g }}{ \overline{\pi_n^+ g}} \,\mytop{\overline{\pi^-_{n-1} g }}{ \overline{\pi_{n-1}^+ g}} \, \dots \, \mytop{\overline{\pi^-_0 g }}{ \overline{\pi_0^+ g}}  \right) \]

where  $\overline{x}$ denotes the projection of $x \in G_n$ in $K_n = G_n/G_{n-1}$.

In this description of $G$, the subspace $G_n$ corresponds to the sequences that are zero for rank $>n$. the subspace $PG_n^-$ and $PG_n^+$ corresponds respectively to the elements of the form:
\[ \left( \dots \, \mytop{0}{0} \,  \mytop{k}{k} \, \mytop{0}{\partial(k)} \, \mytop{0}{0} \, \dots \, \mytop{0}{0}  \right) \]
\[ \left( \dots \, \mytop{0}{0} \,  \mytop{k}{k} \, \mytop{-\partial(k)}{0} \, \mytop{0}{0} \, \dots \, \mytop{0}{0}  \right) \]

(where $k \in K_n$ and appears in degree $n$).

In particular, the isomorphisms from $G$, given in this description, to $\bigoplus PG_n^{-}$ send a double sequence $(k_n^-,k_n^+)$ to the top sequence $(k_n^-)$.

\end{itemize}

}

\subsection{Linearization of $\infty$-categories and the $\delta$ counting function}

\label{sub_sec_linearization}

\blockn{In this subsection we will describe the ``Free abelian group'' functor $C \mapsto \Z C$ on $\infty$-category, i.e. the left adjoint to the forgetful functor from the category of strict $\infty$-groups to the category of strict $\infty$-categories.}

\block{\label{def_linearizationOfInfinityCat}\Def{Let $C$ be a strict $\infty$-category. One defines $\Z C$ as the abelian group freely generated by symbols $\delta(c)$ for $c$ any arrows of $C$ and subject to the relations:

\[ \delta(x \#_k y ) = \delta(x)+\delta(y) - \delta(\pi^+_k(x)) \]

when $x \#_k y$ exists, i.e. $\pi^+_k(x) = \pi^-_k(y)$.
}

This construction is functorial: if $f:X \rightarrow Y$ , then $f$ induces a unique function $f:\Z X \rightarrow \Z Y$ such that $f(\delta(x)) = \delta(f(x))$.

The motivation for this definition should be clear: any morphism of $\infty$-categories from $C$ to a strict $\infty$-group will factor uniquely into $\Z C$ has the operation $x \#_k y$ computed in an $\infty$-group will be given by the formula above due to proposition \ref{Prop_equiv_globgroup_strictinfgroup}.
}

\block{\label{Prop_linear_infty_group}\Prop{For $C$ a strict $\infty$-category, $\Z C$ carries a structure of globular group $\pi^{\epsilon}_n : \Z C \rightarrow \Z C$ characterized by:

\[ \pi^{\epsilon}_n(\delta(f)) = \delta(\pi^{\epsilon}_n(f)) \] 

for all $f\in C$.

\bigskip

Moreover, when $\Z C$ is endowed with the corresponding structure of $\infty$-group, $\delta : C \rightarrow \Z C$ is a morphism of strict $\infty$-categories which makes $\Z C$ the free strict $\infty$-group generated by $C$.

}

\Dem{ We first check that $\pi^{\epsilon}_n$ is well defined. One needs to check that $v: f \mapsto \delta(\pi^{\epsilon}_n(f))$ satisfies $v(x \#_k y)=v(x)+v(y)-v(\pi_k^+(x))$ so that one can apply the universal property of $\Z C$ to extend it into a linear endomorphism. One needs to deal with several cases:

\begin{itemize}

\item If $n \leqslant k$ and $\epsilon=+$.

Then $\pi^+_n(x \#_k y) = \pi_n^+(y)$ so that $v(x \#_k y) = v(y)$ but  $v(\pi_k^+(x)) = \delta(\pi_n^+\pi_k^+(x)) = \delta(\pi_n^+(x))=v(x)$ hence $v(x \#_k y ) = v(y) +v(x) - v(\pi_k^+(x))$.

\item  If $n \leqslant k$ and $\epsilon=-$

Then $\pi^-_n(x \#_k y) = \pi_n^-(x)$ so that $v(x \#_k y) = v(x)$ but  $v(\pi_k^+(x)) = v(\pi_k^-(y)) = \delta(\pi_n^-\pi_k^-(y)) = \delta(\pi_n^-(y))=v(y)$ hence $v(x \#_k y ) = v(x) +v(y) - v(\pi_k^+(x))$.

\item If $n>k$, then $\pi_n^{\epsilon}(x \#_k y) = \pi^{\epsilon}_n(x) \#_k \pi^{\epsilon}_n(y)$, hence $v(x \#_k y) = v(x)+v(y) - \delta(\pi_k^+\pi_n^- x )$; but $\pi_k^+ \pi_n^{-} x = \pi_k^+ x = \pi_n^{-} \pi_k^+ x $ hence $\delta(\pi_k^+\pi_n^- x ) = v(\pi_k^+ x)$. which proves our claim.

\end{itemize}

The fact that the $\pi_n^{\epsilon}$ defined like this on $\Z C$ form a globular structure is immediately checked on generators because the $\pi_n^{\epsilon}$ commute to $\delta$ and the globular relations are satisfied in $C$. This makes $\Z C$ into a globular group, and hence a strict $\infty$-group with the composition $x \#k y =x +y - - \pi_k^{+}(x)$. The fact that $\delta: C \rightarrow \Z C$ is a morphism of $\infty$-categories follow from the fact that $\delta$ commutes to the $\pi_n^{\epsilon}$ by the definition of $\pi_n^{\epsilon}$ on $\Z C$, and the fact that $\delta $ commutes to composition is the relation $\delta(x \#_k y) = \delta(x)+\delta(y) - \delta(\pi_k^+(x))$ which is the definition of $\Z C$.

Finally, if $f:C \rightarrow G$ is any $\infty$-functor to an $\infty$-group, then by definition its satisfies $f(x \#_k y) =f(x)+f(y)-f(\pi_k^{+} x )$ and hence it factors uniquely as a group homomorphism $\Z C \rightarrow H$ and one can check on the generators of $\Z C$ that this group homomorphism additionally commutes to the globular structure.

}
}

\block{\Prop{\begin{itemize}

\item For any $k$, let $C_k$ be the subcategory of $C$ of arrow of dimension at most $k$, and let $(\Z C)_k$ be the sub strict $\infty$-group of $\Z C$ of arrow of dimension at most $k$. Then the natural inclusion $C_k \rightarrow C$ induce an isomorphisms from $\Z (C_k)$ to $(\Z C)_k$.

\item For any $k$, $P(\Z C)_k = \Z C_k / \Z C_{k-1}$ is the abelian group generated by symbols $[f]_k$ (abbreviated to $[f]$)  for $f$ among the $k$-arrow of $C$ satisfying the relations $[ x \#_i y] = [x]+[y]$ for all $i <k$. The $[f]$ for $f$ of dimension $k-1$ are equal to zero in $P(\Z C)_k$.

\end{itemize}

}

The first point will let us remove the parenthesis and simply write $\Z C_k$.

\Dem{\begin{itemize}

\item One defines a map $P_k^{+} : \Z C \rightarrow \Z (C_k)$ such that $P_k^+ \delta(c) = \delta( \pi_k^+ c)$ , which exists because $f(c):=  \delta( \pi_k^+ c)$ satisfies the relation $f(x \#_i y) =f(x)+f(y)-f(\pi_i^+ x)$, by the exact same proof as the proof given in \ref{Prop_linear_infty_group} that $\pi_k^{\epsilon}$ is well defined. Let $i_k : \Z (C_k) \rightarrow \Z C$ be the map induce by the natural inclusion $C_k \rightarrow C$, i.e. $i_k(\delta(c)) = \delta(c)$ for all $c \in C_k$.

By evaluating these on generators, one easily see that $P_k^+ \circ i_k = Id_{C_k}$ and $i_k \circ P_k^+ = \pi_k^+$. This proves that $\Z (C_k)$ is the retract of $\Z C$ by the projection $\pi_k^+$, i.e. is the subgroup of cell of dimension at most $k$.

\item $P(\Z C)_k$ is $\Z C_k / \Z C_{k-1}$, and we define $[x]$ to be the projection of $\delta(x)$. The description of $\Z C_k$ and $\Z C_{k-1}$ as free group, immediately identifies $P(\Z C)_k$ with the abelian group generated by the $k$-arrow of $C$ with the relations $[x \#_n y ]=[x]+[y] - [\pi_n^+(x)]$ for all $n$, and $[x]=0$ when $x$ is a $(k-1)$-arrow. These relations immediately implies that $[x\#_n y] = [x]+[y]$ for all $n<k$ as $\pi_n^{+}(x)$ is a $k-1$-arrow.

Conversely, the relation $[x\#_n y] = [x]+[y]$ for all $n<k$ implies first that if $x$ is a $k-1$ arrow then $[x]=0$ as the relation $x \#_{k-1} x =x$ gives us $[x]+[x]=[x]$. Is also implies that one has $[x \#_n y ]=[x]+[y] - [\pi_n^+(x)]$ for all $n<k$ as $[\pi_n^+(x)]=0$. Finally as one only have $k$-arrows, the relation $[x \#_n y ]=[x]+[y] - [\pi_n^+(x)]$ for $n\geqslant k$ are always trivial. So we have proved that the relations $[x\#_n y] = [x]+[y]$ for all $n<k$ is equivalent to the relation defining $P(\Z C)_k$. 

\end{itemize}

}
}

\block{\label{prop:LinearizedPolygraphAreFree}\Prop{If $X$ is a polygraph, then $\Z X^*$ is freely generated as an abelian group by the $\delta(x)$ for $x$ among the generating cells of $X$.}

Because of this proposition, we will omit the ``$*$'' and simply write $\Z X$ instead of $\Z X^*$ when $X$ is a polygraph.

\Dem{ We will prove by induction on $k$ that $\Z X_k$ is free generated by the cells of $X$ of dimension less than $k$.

It is immediate in dimension $0$. We assume the result in dimension $k-1$. A morphism from $\Z X_k$ to another strict $\infty$-group $G$ is the same as a functor from $X_k^*$ to $G$. Such a functor is given by a functor from $(X_{k-1})^*$ to $G$, together with a choice of the image of each $k$-cell of $X$ in $G$ with correct source and target. As a functor from $(X_{k-1})^*$ to $G$ is the same as a morphism from $\Z X_{k-1}$ to $G$ this shows that $\Z X_k$ is freely generated as a strict $\infty$-group from $\Z X_{k-1}$ by adding $\delta(x)$ for each $k$-cell $x$ of $X$ with specified values for the $\pi_{k-1}^{\epsilon} \delta(x)$. But this makes the underlying group of $\Z X_k$ freely generated by $\Z X_{k-1}$ and the $\delta(x)$ for $x$ a $k$-cell.
}
}

\subsection{Tensor product of Strict $\infty$-group}

\label{subsec:TensorProdOfInfinityGroup}

\blockn{The category of chain complexes have a natural monoidal structure given by the tensor product $K_{\bullet} \otimes L_{\bullet} = M_{\bullet}$ where:

 \[ M_k = \bigoplus_{i=0}^k K_i \otimes L_{k-i} \]

and is endowed with the differential defined by $\partial(k \otimes l) = \partial k \otimes l +(-1)^{i} k \otimes \partial l$ where $i$ is the dimension of $k$ (i.e. $k \in K_i$).

As we have explained in subsection \ref{subsection_StrictAbGrp_vs_chaincomplexes}, the category of chain complexes is equivalent to the category of strict $\infty$-groups, hence this monoidal structure transport through this equivalence into a tensor product on the category of $\infty$-groups. Our goal in the present subsection is to give an explicit description of this tensor product of globular groups.

While the functor from strict $\infty$-group to chain complex is relatively simple, its inverse is a little more complicated and it is a little technical to use it explicitly to construct the tensor product of strict $\infty$-group. So to make the proof simpler, we will define a tensor product on strict $\infty$-group, and show that functor $G \mapsto K_{\bullet} = G_{\bullet}/G_{\bullet -1}$ is a monoidal equivalence.
}

\block{\label{Def_tensor_product_infinitygroup}\Def{Let $G$ and $H$ be two strict $\infty$-groups. We define their tensor product as the tensor product of groups $G \otimes H$ endowed with the operation:

\[ \pi_n^{\epsilon}(g \otimes h) = \sum_{i} (\pi^{\epsilon}_i-\pi_{i-1}^{\epsilon}) g \otimes \pi^{(-1)^i \epsilon}_{n-i} h  \]

}

We remind the reader of the convention that $\pi^{\epsilon}_i$ is zero when $i<0$ so that the sum is actually finite. As it will become quite systematic to have $(-1)^i$ as exponent of $\pi$ we will introduce the notation $|i|$ for $(-1)^i$ when used in such an exponent.

\Lem{ $G \otimes H$ as defined above is a strict $\infty$-group, i.e. the operation $\pi^{\epsilon}_n$ that we defined satisfies the globular relations.
}

\Dem{We start by a general computation of $\pi_{n}^{\epsilon} \pi_m^{\delta}$ and we then check the relevant cases one by one. One has:

\[ \pi_n^{\epsilon} \pi_m^{\delta}(g \otimes h) = \sum_{i} \sum_j (\pi_i^{\epsilon}\pi_j^{\delta} - \pi_i^{\epsilon}\pi_{j-1}^{\delta} - \pi_{i-1}^{\epsilon}\pi_j^{\delta} + \pi_{i-1}^{\epsilon}\pi_{j-1}^{\delta}  ) g \otimes \pi^{|i|\epsilon}_{n-i} \pi^{|j|\delta}_{m-j} h \]

Now if $i-1 \leqslant j$  then:

\[ \pi_i^{\epsilon}\pi_j^{\delta} - \pi_i^{\epsilon}\pi_{j-1}^{\delta} - \pi_{i-1}^{\epsilon}\pi_j^{\delta} + \pi_{i-1}^{\epsilon}\pi_{j-1}^{\delta} = \pi_i^{\epsilon} - \pi_i^{\epsilon}  - \pi_{i-1}^{\epsilon}  + \pi_{i-1}^{\epsilon}  = 0   \]

and if $j-1 < i$ :

\[ \pi_i^{\epsilon}\pi_j^{\delta} - \pi_i^{\epsilon}\pi_{j-1}^{\delta} - \pi_{i-1}^{\epsilon}\pi_j^{\delta} + \pi_{i-1}^{\epsilon}\pi_{j-1}^{\delta} = \pi_j^{\delta} - \pi_{j-1}^{\delta} - \pi_{j}^{\delta} + \pi_{j-1}^{\delta} = 0   \]

So the sum above is only non trivial if $i=j$ or $i=j-1$ moreover if $i=j$:

\[ \pi_i^{\epsilon}\pi_j^{\delta} - \pi_i^{\epsilon}\pi_{j-1}^{\delta} - \pi_{i-1}^{\epsilon}\pi_j^{\delta} + \pi_{i-1}^{\epsilon}\pi_{j-1}^{\delta} =\pi_j^{\delta} - \pi_{j-1}^{\epsilon} \]

and if $i=j-1$ then:

\[  \pi_i^{\epsilon}\pi_j^{\delta} - \pi_i^{\epsilon}\pi_{j-1}^{\delta} - \pi_{i-1}^{\epsilon}\pi_j^{\delta} + \pi_{i-1}^{\epsilon}\pi_{j-1}^{\delta} = \pi_{j-1}^{\epsilon} - \pi_{j-1}^{\delta} \]

Hence, in general  $\pi_n^{\epsilon} \pi_m^{\delta}(g \otimes h)$ is equal to:

\[ \left( \sum_j (\pi_j^{\delta} - \pi_{j-1}^{\epsilon}) g \otimes \pi^{|j|\epsilon}_{n-j} \pi^{|j| \delta}_{m-j} h   \right) +  \left( \sum_j (\pi_{j-1}^{\epsilon} - \pi_{j-1}^{\delta}) g \otimes \pi^{-|j|\epsilon}_{n-j+1} \pi^{|j| \delta}_{m-j} h   \right) \]

One needs to treat four cases separately:

\begin{itemize}

\item If $m  \leqslant n$ the formula above reduce to:

\[ \left( \sum_j (\pi_j^{\delta} - \pi_{j-1}^{\epsilon}) g \otimes \pi^{|j| \delta}_{m-j} h   \right) +  \left( \sum_j (\pi_{j-1}^{\epsilon} - \pi_{j-1}^{\delta}) g \otimes  \pi^{|j| \delta}_{m-j} h   \right)  \]

And one can regroup the two sums, factor the $ \otimes  \pi^{|j| \delta}_{m-j} h$, simplify the two $\pi_{j-1}^{\epsilon}$ and we obtain exactly the definition of $\pi_m^{\delta}(g \otimes h)$.

\item If $n=m-1$ and $\delta= \epsilon$ then the second terms of the sum vanish and the first term is the definition of $\pi_n^{\epsilon}$.

\item If $n=m-1$ and $\delta = - \epsilon$ then the formula reduce to:

\[ \left( \sum_j (\pi_j^{\delta} - \pi_{j-1}^{\epsilon}) g \otimes \pi^{|j|\epsilon}_{n-j} h   \right) +  \left( \sum_j (\pi_{j-1}^{\epsilon} - \pi_{j-1}^{\delta}) g \otimes  \pi^{|j| \delta}_{n+1-j} h   \right) \]

A change of variables, and replacing the appearing $-\delta$ with $\epsilon$ in the second sum gives:

\[ \left( \sum_j (\pi_j^{\delta} - \pi_{j-1}^{\epsilon}) g \otimes \pi^{|j|\epsilon}_{n-j} h   \right) +  \left( \sum_j (\pi_{j}^{\epsilon} - \pi_{j}^{\delta}) g \otimes  \pi^{|j| \epsilon}_{n-j} h   \right) \] 

regrouping and simplifying the $\pi_j^{\delta}$ gives the definition of $\pi_{n}^{\epsilon}(g \otimes h)$.

\item Finally, if $m>n+1$ then the formula reduce to:

\[ \left( \sum_j (\pi_j^{\delta} - \pi_{j-1}^{\epsilon}) g \otimes \pi^{|j|\epsilon}_{n-j} h   \right) +  \left( \sum_j (\pi_{j-1}^{\epsilon} - \pi_{j-1}^{\delta}) g \otimes \pi^{-|j|\epsilon}_{n-j+1} h   \right) \]

And a change of variables in the second sum gives:

\[ \left( \sum_j (\pi_j^{\delta} - \pi_{j-1}^{\epsilon}) g \otimes \pi^{|j|\epsilon}_{n-j} h   \right) +  \left( \sum_j (\pi_{j}^{\epsilon} - \pi_{j}^{\delta}) g \otimes \pi^{|j|\epsilon}_{n-j} h   \right) \]

regrouping and simplifying the $\pi_j^{\delta}$ gives the definition of $\pi_{n}^{\epsilon}(g \otimes h)$.

\end{itemize}

}
}

\block{\label{claim_pitens_alternate}The definition of $\pi_n^{\epsilon}$ can seems asymmetrical but one can also write it as:

\[\pi_n^{\epsilon}(g \otimes h) = \sum_i \pi_i^{\epsilon} g \otimes \left( \pi_{n-i}^{|i| \epsilon} - \pi_{n-i-1}^{-|i| \epsilon}  \right) h \]

which somehow restore a form of symmetry. Indeed:

\begin{align*} \pi_n^{\epsilon}(g \otimes h) & =  \sum_{i} (\pi^{\epsilon}_i-\pi_{i-1}^{\epsilon}) g \otimes \pi^{|i| \epsilon}_{n-i} h \\
 & =   \sum_{i} \pi^{\epsilon}_i g \otimes \pi^{|i| \epsilon}_{n-i} h -  \sum_{i} \pi_{i-1}^{\epsilon} g \otimes \pi^{|i| \epsilon}_{n-i} h \\
 & =   \sum_{i} \pi^{\epsilon}_i g \otimes \pi^{|i| \epsilon}_{n-i} h -  \sum_{i} \pi_{i}^{\epsilon} g \otimes \pi^{-|i| \epsilon}_{n-i-1} h \\
 & =   \sum_i \pi_i^{\epsilon} g \otimes \left( \pi_{n-i}^{|i| \epsilon} - \pi_{n-i-1}^{-|i| \epsilon}  \right) h
 \end{align*}

}

\block{\Prop{This tensor product defines a monoidal structure on the category of globular groups. The unit object is $\Z$ with all $\pi_i^{\epsilon}=Id$, and all the structure isomorphisms (associativity and unit law) are those of the tensor product of groups.}

\Dem{one just needs to check that the structure isomorphism of associativity and unit law are compatible to the globular structure. For the unit law:

\[ \pi_n^{\epsilon}(g \otimes  1) =  \sum_i \pi_i^{\epsilon} g \otimes \left( \pi_{n-i}^{|i| \epsilon} - \pi_{n-i-1}^{-|i| \epsilon}  \right) 1 \]

but the left side of the tensor product is zero for all value of $i$ except $i=n$ where it is $1$ (as $\pi_0=1$ but $\pi_{-1}=0$), hence one indeed obtains $\pi^{\epsilon}_n(g) \otimes 1 $. The case of $1 \otimes g$ is exactly the same.

For the associativity, the different way of computing $\pi_n^{\epsilon}(x \otimes y \otimes z)$ using the different way of associating the tensor product and the two form of $\pi_n^{\epsilon}(g \otimes h)$ can all be expended to one of the following three expressions:

\[ \sum_{u+v+w=n} \left[ \pi_u^{\epsilon} - \pi_{u-1}^{\epsilon} \right] x \otimes \left[ \pi_v^{\epsilon |u|} - \pi_{v-1}^{\epsilon |u|} \right] y \otimes \pi_w^{\epsilon|u+v|} z\]

\[ \sum_{u+v+w=n} \left[ \pi_u^{\epsilon} - \pi_{u-1}^{\epsilon} \right] x \otimes  \pi_v^{\epsilon |u|} y \otimes \left[ \pi_w^{\epsilon|u+v|} -  \pi_{w-1}^{-\epsilon|u+v|}  \right] z\]

\[ \sum_{u+v+w=n} \pi_u^{\epsilon} x \otimes \left[ \pi_v^{\epsilon |u|} - \pi_{v-1}^{-\epsilon |u|} \right] y \otimes \left[ \pi_w^{\epsilon|u+v|} -  \pi_{w-1}^{-\epsilon|u+v|}  \right] z\]

and one can pass from any of these expressions to any other using exactly the same computation as in \ref{claim_pitens_alternate}.

}
}

\block{\label{Prop_Projection_tensordecomp}\Prop{

\[ \left[ \pi_n^{\epsilon} -  \pi_{n-1}^{\delta} \right] (x \otimes y) = \sum_{i+j=n} \left[\pi_i^{\epsilon} - \pi_{i-1}^{\delta} \right] x \otimes \left[ \pi_{j}^{\epsilon|i|} - \pi_{j-1}^{\delta|i|} \right] y \]

}

\Dem{One treats the case $\delta= \epsilon$ and $\delta=-\epsilon$ separately. In the case $\delta=\epsilon$, one start by applying both the definition of $\pi_{n}^{\epsilon}(x \otimes y)$ and of $\pi^{\epsilon}_{n-1}(x \otimes y)$ and to expand them. This gives us that the left hand side is equal to:

\begin{multline*} \sum_{i+j=n} \left(  \pi_i^{\epsilon} x \otimes \pi_j^{\epsilon |i|} y - \pi_{i-1}^{\epsilon} x \otimes \pi_j^{\epsilon |i|} y \right) \\ -  \sum_{i+j=n-1} \left(  \pi_i^{\epsilon} x \otimes \pi_j^{\epsilon |i|} y - \pi_{i-1}^{\epsilon} x \otimes \pi_j^{\epsilon |i|} y \right)  \end{multline*}

One makes the change of variables $j \mapsto j-1$ in the second sum and after regrouping the sum one gets:

\[ \sum_{i+j=n} \pi_i^{\epsilon} x \otimes \pi_j^{\epsilon |i|} y - \pi_{i-1}^{\epsilon} x \otimes \pi_j^{\epsilon |i|} y  -   \pi_i^{\epsilon} x \otimes \pi_{j-1}^{\epsilon |i|} y + \pi_{i-1}^{\epsilon} x \otimes \pi_{j-1}^{\epsilon |i|} y  \]

which is exactly the right hand side in expanded form.

In the case $\delta=-\epsilon$ one instead start from the alternate definition of \ref{claim_pitens_alternate} for both $\pi_n$ and $\pi_{n-1}$. The left hand side is equal to:

\begin{multline*} \sum_{i+j=n} \left(  \pi_i^{\epsilon} x \otimes \pi_j^{\epsilon |i|} y - \pi_{i}^{\epsilon} x \otimes \pi_{j-1}^{-\epsilon |i|} y \right) \\ -  \sum_{i+j=n-1} \left(  \pi_i^{-\epsilon} x \otimes \pi_j^{-\epsilon |i|} y - \pi_{i}^{-\epsilon} x \otimes \pi_{j-1}^{\epsilon |i|} y \right)  \end{multline*}

this times one makes the change of variable $i\mapsto i-1$ in the second sum and one gets after regrouping:

\[ \sum_{i+j=n}   \pi_i^{\epsilon} x \otimes \pi_j^{\epsilon |i|} y - \pi_{i}^{\epsilon} x \otimes \pi_{j-1}^{-\epsilon |i|} y -    \pi_{i-1}^{-\epsilon} x \otimes \pi_j^{\epsilon |i|} y + \pi_{i-1}^{-\epsilon} x \otimes \pi_{j-1}^{-\epsilon |i|} y  \]

which, once we replace all the $-\epsilon$ by $\delta$, is the expanded form of the right hand side of the equality in the proposition.

}
}

\block{\label{Lem_pi_tens_degminus1}\Prop{Let $g \in G_i$ and $h \in H_j$, then for any $n \leqslant i+j$ :

\[ \pi_n^{\epsilon}(g \otimes h) = \sum_{\mytop{a+b=n}{a \leqslant i , b \leqslant j}} \pi^{\epsilon}_a(g)\otimes \pi^{\epsilon |a|}_b(h) - \sum_{\mytop{a+b=n-1}{a < i , b < j}} \pi_a(g) \otimes \pi^{-\epsilon |a|}_b(h) \]

in particular:

\begin{itemize}

\item For $n=i+j$:  $\pi_n^{\epsilon}(g \otimes h)= g \otimes h$, i.e. $g \otimes h \in (G \otimes H)_{i+j}$. In particular one also has $\pi^{\epsilon}_n(g \otimes h) = g \otimes h$. for all $n \geqslant i+j$.

\item For $n=i+j-1$ :

\[ \pi_{i+j-1}^{\epsilon}(g \otimes h) = \pi^{\epsilon}_{i-1}(g) \otimes h + g \otimes \pi_{j-1}^{\epsilon |i|}(h) -  \pi^{\epsilon}_{i-1}( g ) \otimes \pi_{j-1}^{\epsilon |i|}(h)  \]

\end{itemize}

}

More generally, if $n = i+j-k$ , the first sum as $k+1$ terms and the second sum has $k$ terms exactly.

\Dem{

\[ \pi_n^{\epsilon}(g \otimes h) = \sum_{a+b = n} (\pi_a^{\epsilon} - \pi_{a-1}^{\epsilon} ) g \otimes \pi_b^{|a| \epsilon} h  \]

For $a >i$, one has $\pi^{\epsilon}_a(g)= g = \pi^{\epsilon}_{a-1} g$ hence the left side of the tensor product is zero, so that one can restrict the sum to $a \leqslant i$:

\[ \pi_n^{\epsilon}(g \otimes h) = \sum_{\mytop{a+b = n}{a \leqslant i}} (\pi_a^{\epsilon} - \pi_{a-1}^{\epsilon} ) g \otimes \pi_b^{|a| \epsilon} h  \]
 
 For $b \geqslant j$ one has $\pi_b^{|a| \epsilon} h = h$, so if we split the sum into two one gets:

\[\pi_n^{\epsilon}(g \otimes h) =  \sum_{\mytop{a+b = n}{\mytop{a \leqslant i}{b < j }}} (\pi_a^{\epsilon} - \pi_{a-1}^{\epsilon} ) g \otimes \pi_b^{|a| \epsilon} h + \sum_{\mytop{a+b = n}{\mytop{a \leqslant i}{b \geqslant j }} }  (\pi_a^{\epsilon} - \pi_{a-1}^{\epsilon} ) g \otimes h. \]

As $n \leqslant i+j$, the condition $a \leqslant i$ in the second sum can be removed, so one has a telescoping sum whose value is $\pi_{n-j}^{\epsilon} g \otimes h$:

\[\pi_n^{\epsilon}(g \otimes h) =  \left( \sum_{\mytop{a+b = n}{\mytop{a \leqslant i}{b < j }}} \pi_a^{\epsilon} g \otimes \pi_b^{|a| \epsilon} h \right) + \pi_{n-j}^{\epsilon} g \otimes h - \left(  \sum_{\mytop{a+b = n}{\mytop{a \leqslant i}{b < j }}} \pi_{a-1}^{\epsilon} g \otimes \pi_b^{|a| \epsilon} h \right) \]

inserting the middle terms in the left sum just extend the summation to $b \leqslant j$ instead of $b <j$ (the case $b=j$ indeed corresponds to a case where $a \leqslant i$ as $n \leqslant i+j$) and the change of variable $a \mapsto a'+1$ gives exactly the second terms of the formula in the proposition so this concludes the proof of the main statement, the others follow easily.

}
}

\block{\Prop{ For $G$ and $H$ two globular groups:

\begin{itemize}

\item one has the equality:

\[ P(G \otimes H)_n^{\delta} = \bigoplus_{i+j=n} PG_i^{\delta} \otimes PH_{j}^{|i| \delta} \]

as subspaces of $G \otimes H$. 

\item If $g \in G_i$, $h \in H_j$ then the class of $g \otimes h$ in $P(G \otimes H)_{i+j}$ only depends on the class of $g$ and $h$ in $PG_i$ and $PH_j$. This induces a map $PG_i \otimes PH_j \rightarrow P(G\otimes H)_{i+j}$ which induces an isomorphism:

\[ P(G \otimes H)_n \simeq \bigoplus_{i+j=n} PG_i \otimes PH_{j} \]

Compatible to the identification above.

\item Under this identification, the differential of $P(G \otimes H)_{\bullet}$ is the usual differential of a tensor product of chain complexes.

\item This makes the equivalence between globular groups and chain complexes an equivalence of monoidal categories. 

\end{itemize}
}

\Dem{\begin{itemize}

\item This follows from proposition \ref{Prop_Projection_tensordecomp} above: $\pi_n^{\epsilon} -  \pi_{n-1}^{\delta}$ on $G \otimes H$ is a projection on $P(G \otimes H)_n^{\delta}$, and for each $i+j=n$, the terms of the sum $\left[\pi_i^{\epsilon} - \pi_{i-1}^{\delta} \right] \otimes \left[ \pi_{j}^{\epsilon|i|} - \pi_{j-1}^{\delta|i|} \right] $ is a projection of $PG_i^{\delta} \otimes PH_{j}^{|i| \delta}$. 

Now one can check that a product of the form:

\[ \left( \pi_i^{\epsilon} - \pi_{i-1}^{\delta} \right) \left( \pi_j^{\epsilon'} - \pi_{j-1}^{\delta'} \right) \]

vanish as soon as: $i \geqslant j+1$, or $i<j-1$, or $i=j-1$ and $\epsilon=\delta'$, hence the product of any two distinct terms of our sum vanish (for all value of $i$ and $j$ except $i=j-1$ and $i=j$ both side of the tensor product vanish, for $i=j-1$ one of the two side vanish, more precisely, the left side if $\epsilon=\delta$ and the right side if $\epsilon=-\delta$.
And when a projection $P$ is decomposed as a sum of projection $P_i$ whose pairwise (distinct) products are all zero, it means that the range of $P$ is the direct sum of the range of all the $P_i$, which gives the decomposition by the first point.

\item If $g'=g+a$ with $g \in G_i$ and $a \in G_{i-1}$ and $h'=h+b$ with $h \in H_j$ and $b \in H_{j-1}$ then: 

\[ g' \otimes h ' = g \otimes h +a \otimes h + g \otimes b+ a \otimes b \]

and because of the first special case of proposition \ref{Lem_pi_tens_degminus1}, all the terms except the first are in $(G \otimes H)_{i+j-1}$, hence $g' \otimes h'$ and $g \otimes h$ indeed have the same image in $P(G \otimes H)_{i+j}$.

This proves in particular that one has a well defined map:

 \[ \bigoplus_{i+j=n} PG_i \otimes PH_j \rightarrow P(G \otimes H)_{n} \]

but by definition, this map lifts into the identity map on 

\[ \bigoplus_{i+j=n} PG_i^{\delta} \otimes PH_{j}^{|i| \delta} =  P(G \otimes H)_n^{\delta}  \]

and the projection restricted to these subspaces are isomorphisms, hence the map we just constructed is also an isomorphism.

\item For $g \in PG_i$ and $h \in PH_j$, one chose the representant of $g$ and $h$ such that $g \in PG_{i}^-$ and $h \in PG_{j}^{-|i| }$

We already know that $g \otimes h  \in P(G \otimes H)_n^{-}$, so that its differential reduces to $\pi_{n-1}^+(g\otimes h)$ which we compute using the second special case of proposition \ref{Lem_pi_tens_degminus1}:

\[\pi_{n-1}^{+}(g \otimes h) = \pi^{+}_{i-1}(g) \otimes h + g \otimes \pi_{j-1}^{|i|}(h) - \pi^{+}_{i-1}(g) \otimes  \pi_{j-1}^{|i|}(h)  \]

The last terms is in $(G\otimes H)_{n-2}$ so it plays no role in the projection in $P(G\otimes H)_{n-1}$, moreover $\pi_{i-1}^+(g)$ is $\partial(g)$ in $PG_{i-1}$  and as $h \in PG_{j}^{-|i|}$, the term $\pi_{j-1}^{\delta|i|}(h)$ is $(-1)^i \partial h$ in $PH_{j}$ hence one indeed obtains the desired formula for the differential.

\item One needs to check that the isomorphism between $P(G\otimes H)_{\bullet}$  and the tensor product of chain complex $PG_{\bullet} \otimes PH_{\bullet}$ satisfies the usual compatibility axioms. This is completely immediate once the axiom and the structure map of the tensor product of chain complexes are written out.

\end{itemize}

}
}

\section{ $\delta$ as a counting function}
\label{section:Delta}

\blockn{By a ``counting function'' we mean some function which, given an arrow $f \in X^*$, define a ``number of occurrence'' of each cell $x$ in $f$. If $f$ is a $k$-arrow and $x$ is a $k$-cell then this is a well defined concept (cf \ref{Prop_delta_count_in_top_Dim} below, but this was known long before) which actually corresponds to the number of occurrence of $x$ in any expression of $f$ in terms of the generators of $X$. But when the dimension of $x$ is strictly lower than the dimension of $f$ this become a more complicated concept. In general such a counting function can be seen as a mapping from $X^*$ to the set of formal linear combination of cells of $x$.

In \cite{makkai2005word} (section 5, proposition 12), Makkai introduced such a counting function, that he calls the content of an arrow. We will denote it $c_M$, it is characterized by the following properties:

\begin{itemize}

\item If $x$ is a $k$-cell of $X$, then:

 \[ c_M(x) = e_x + c_M(\pi_{k-1}^+(x)) + c_M(\pi_{k-1}^-(x)) \]

\item  $c_M(1_v)=c_M(v)$ (Or equivalently, it is defined on arrows in the sense of the ``one sort'' definition of $\infty$-categories). 

\item $c_M(f \#_k g) = c_M(f)+c_M(g) - c_M(\pi_k^+(f))$

\end{itemize}

Those formulas allow to compute $c_M$ on any expression of $f$ in terms of the generating cells, and Makkai results is that this does not depends on the expression of $f$. The intuitive idea is the following: one wants to count all the cell that appears in $f$, including in its boundary. So $c_M( 1_v)$ should indeed by the same as $C_m(v)$, and, at least in good case the formula for the composite is natural as when one forms $c_M(f)+c_M(g)$ the cells appearing in $\pi_k^+(f) = \pi_k^-(g)$ have been counted twice. But there is a small problem with this heuristic in the definition of $c_M(x)$ for $x$ a generating cell, one indeed want to consider that $x$ appears exactly once in $x$, and all the cell of the boundary of $x$ also appeared, but the formula Makkai gave will count twice the cell that are in $\pi_{k-2}^{\epsilon} (x)$. So $c_M$ is slightly ``over counting''. But $c_M$ has a very desirable property for a counting function: one can prove (see \cite{makkai2005word} section 5 proposition 12) that $c_M(f)$ is always a linear combination with positive coefficients, and any attempt to ``fix'' this over counting, while keeping a simple inductive definition of $c_M$ seem to always destroy this nice property.
}

\blockn{Another possible attempt to define a relevant counting function is the following: if $X$ is a positive polygraph, and $v \in X^*$ is an arrow, then one can consider the classifying polyplex $\chi_v : \underline{v} \rightarrow X$ of $v$, and then simply count the size of the pre-image in $\chi_v$ of each cell of $x$, that would define:

\[ c_p(v) := \sum_{a \in \underline{v}} e_{\chi_v(a)} \]

This counting function has very good properties, and characterize polyplexes: a polygraph $X$ with an arrow $f \in X^*$ is a polyplex if and only if every cell of $X$ appear exactly once in $f$ in this sense.

In fact Makkai explicitly ask in \cite{makkai2005word} whether such a counting function exists\footnote{but this one does not totally answer Makkai's question as it is only defined on positive polygraph, or more generally in a good class of polygraph but not on all polygraphs}. But it has a big problem: there is apparently no known ``inductive'' way to compute it on an expression of $v$ in the same kind of way as Makkai functions. In particular one cannot compute it without already being able to identifies polyplexes.
}

\blockn{Another approaches would be to try to count only the cells that appears as ``inner cells'' in an arrow $f$, i.e. not in the boundary. This suggest to define instead something satisfying $c_i(f \#_k g) = c_i(f)+c_i(g)+c_i(\pi_k^+(f))$, as the cell appearing in $\pi_k^+(f)$ where not inner in either $f$ or $g$ but become inner in the composite, but definition of this kind all fail to be independent of the expression of the cell (hence they are not well defined).}

\blockn{In this section we will investigate the property of the function $\delta :X^* \rightarrow \Z X$ from subsection \ref{sub_sec_linearization}, seen as a counting function. It is closely related to Makkai's content function as it satisfies the same induction formula, except for the base case. In fact $\delta$ is the universal function which satisfies the same induction property as Makkai's content function, but with the values on generating cells left unspecified. See also \ref{prop:RelDeltaMakkai} for a more concrete relation between $\delta$ and Makkai's content function $c_M$.

The ``counting function'' $\delta$ miss one of the big property expected from a counting function: it is not positive ! For example, if $f$ and $g$ are composable $1$-cell, with $x=\pi_0^+ f = \pi_0^- g$, then $\delta(f \#_0 g) = \delta f +\delta g - \delta x$. But we will show that it has a lots of other good properties to compensate for this.

In fact, this counting function, correctly interpreted will reconcile all the different point of view mentioned above in the case of regular polygraphs:

\begin{itemize}

\item There is a nice inductive formula, in fact the same as Makkai's content function to compute it.

\item One will have a characterization of regular polyplex in terms of this counting function (see \ref{Prop_delta_sigma_Charac_polyplex}).

\item If one computes $\delta(v)$ and then change the sign of all the cells appearing in $\delta(v)$ whose dimension does not have the same parity as the dimension of $v$, then one actually obtains a counting function for inner cell of $v$, but which count ``singular\footnote{We mean by that the cells that are involved in a defect of sphericity if the polyplex representing $v$.}'' cells with some more subtle and possibly negative coefficients to obtain a good behavior under composition. See proposition \ref{Prop_pushsigma} for the precise claim that up to this change of sign $\delta$ does count cells correctly for arrows whose representing polyplex is regular with spherical boundary. 
Note that this alternating change of signs, and the fact that singular cells can be counted negatively, explain why the induction relation is of the form $\delta(f \#_k g ) =\delta(f)+\delta(g) -\delta(\pi_k^+ f)$, instead of the expected $\delta(f)+ \delta(g)+\delta(\pi_k^+ f)$ for a function that count inner cells.
\end{itemize}
}

\blockn{Subsection \ref{Subsec_Delta_counting_general} investigate some general properties of this $\delta$ function for arbitrary polygraphs. Subsection \ref{subsec_delta_counting_regular} investigate the property of $\delta$ more specific to the class of regular polygraphs, this will allow to show that regular polygraphs form an algebraic class of polygraph in the sense of subsection \ref{subsec_algebraicclass_def}.

In subsection \ref{subsec_Gray_tensor_product} we will introduce the Gray tensor product of polygraphs and show that $\delta$ behave also nicely with respect to the Gray tensor product and we will combine this with the results of subsection \ref{subsec_delta_counting_regular} to show that regular polygraph indeed form a $\otimes$-stable algebraic class of polygraph in the sense of subsection \ref{subsec_tensorstable_algebraic_class}.
}

\subsection{General properties}
\label{Subsec_Delta_counting_general}

\blockn{To emphasize the special status of the $\delta(x) \in \Z X$ for $x$ a cell of $X$ as a basis of $\Z X$ we will denote them also $\delta_x$ (but only when $x$ is a generating cell)}

\block{Let $X$ and $Y$ be two polygraphs and $f :X^* \rightarrow Y^*$ a functor. By functoriality of linearization one obtains a morphisms of strict $\infty$-group $f : \Z X \rightarrow \Z Y$, it is characterized by the property that $f(\delta_x)= \delta(f(x))$. The functoriality of this constructions already says something interesting on $\delta$ as a counting functions:

If one has $X^* \overset{f}{\rightarrow} Y^* \overset{g}{\rightarrow} Z^*$ and $x$ is a generator of $X^*$, then $f(x)$ has a certain decomposition in terms of the generators of $y$, and for each $y$ a generator of $Y$, $g(y)$ has a decomposition in terms of the generators of $Z$. The functoriality of linearization tell us that $\delta( g(f(x))$ can be computed by first computing $\delta(f(x))$ in $\Z Y$, and then sending each generator $y \in Y$ to $\delta(g(y))$. Thinking $\delta$ as a counting function it means that ``if $z$ is a composition of composed cell $y_i$, then the number of time each cell $x$ appears in $z$ is the sum over all the $y_i$ of the number of time $x$ appears in $y_i$ times the number of time $y_i$ appears in $z$'', which is also a very desirable property of a counting function.

}

\blockn{The following proposition says that $\delta$ is indeed a counting function in some range:}

\block{\label{Prop_delta_count_in_top_Dim}\Prop{ \begin{enumerate}

\item Let $x$ be an arrow of $X^*$ which admit an expression from the generators of $X$ only involving composition operations $\#_k$ for $k<n$. Then for any generator $b$ of dimension $\geqslant n$ the coefficient of $b$ in $\delta(x)$ is the number of occurrence of $b$ in any such expression of $x$. 

\item Let $x$ be an arrow of $X^*$ which admit an expression involving at most one occurrence of a single generator of dimension $>k$, then $x$ admit an expression involving at most one generator of dimension $>k$ and only composition operations $\#_l$ for $l<k$.

\item Let $x$ be an arrow of $X^*$ such that either $\delta(x)$ is in $\Z X_k$ or $\delta(x)-a \in \Z X_k$ for some generator $a$ and some integer $k$, then $x$ admit an expression only involving composition $\#_l$ for $l<k$. In particular, the coefficient of $l$-cell for $l \geqslant k$ in $\delta(x)$ corresponds to the number of occurrences in such an expression of $x$.

\item If $\delta(x) \in \Z X_k$ then $x$ is a $k$-arrow.

\item If $\delta(x)=\delta_a$ for $a$ a generator then $x=a$.

\end{enumerate}
}

\Dem{
The first point is an induction on such expressions for $x$. If $x$ is given by a generator then it is immediate. If $x = z \#_k y$ with $k<n$ then $\delta(x)= \delta(y)+\delta(z)-\delta(\pi^+_k(y))$ but $\pi^+_k(y)$ is in $X_k$ so $\delta(\pi^+_k(y))$ is in $\Z X_k$, i.e. only involves cell of dimension at most $k$, so does not contribute to the coefficient of cell of dimension $\geqslant n$. For such cells the coefficient in $\delta(x)$ is hence the sum of the number of appearances in $y$ and in $z$, i.e. the number of appearances in $x$.

For the second point. If $x$ has at most one generator of dimension $>k$, then any $\#_l$ appearing in $x$ has one side which only contains generators of dimension $\leqslant k$, hence is in $X_k$, and hence is a $k$-arrow. But if $h$ is a $k$-arrow and $l \geqslant k$ then $h \#_l x =x $ and $x \#_l h =x$ when they are defined. Hence each $\#_l$ for $l \geqslant k$ can be eliminated from the expression of $x$.

For the third point, we only need to prove that $x$ admit an expression only involving $\#_l$ for $l<k$, the rest follow from point $1.$.
We start with an arbitrary expression of $x$, it is finite so there is a $L$ such that $x$ only contains composition operations $\#_l$ for $l<L$. If $L \leqslant k$ then we are done. Otherwise one can apply point $1.$ and deduce from the assumption on $\delta(x)$ that this expression of $x$ will contain at most one occurrence of a generator of dimension $> L-1$, hence one can apply point $2.$ and deduce that $x$ has an expression involving only compositions $\#_l$ for $l<L-1$. One can repeat this argument up to the point where $L=k$.

The last two points are special case of the third.
}
}

\blockn{We will now clarify the relation between the function $\delta$ and Makkai's content function $c_M$.}

\block{\label{prop:RelDeltaMakkai}\Prop{For each generator $x \in X$ of dimension $n$, let:

 \[ m_x = \delta_x - \delta(\pi^+_{n-1} x) - \delta(\pi^-_{n-1} x ) \]

then:

\begin{itemize}

\item The $m_x$ for $x \in X$ form a basis of $\Z X$.

\item For any arrow $f \in X^*$, the decomposition of $\delta(f)$ in the basis of the $m_x$ gives the coefficients of Makkai's content function $c_M$.

\end{itemize}

}

The $m_x$ appears to be the opposite of the inverse of $\delta_x$ with respect to $\#_{n-1}$ (where $n$ is the dimension of $x$). But we have absolutely no clues on the relation between this observations the various properties of Makkai's content function.

\Dem{The $m_x$ form a basis simply because their coefficients in the basis $\delta_x$ are triangular (with unit diagonal coefficients) with respect to the filtration of $\Z X$ by the $\Z X_k$. The second claim is immediate by induction: $\delta$ and $c_M$ satisfies the same induction formula and on basis elements one has:

\[ \delta(x) = m_x + \delta(\pi^+_{n-1} x) + \delta(\pi^-_{n-1} x) \]

hence if one assume that the property is true by induction for arrow of dimension $\leqslant n-1$ (in particular for the $\pi_{n-1}^{\epsilon} x$) then it is true for generator of dimension $n$, and then by composition it is true for all arrow of dimension $n$. And $\delta_x=m_x$ for $x$ of dimension $0$ so this also works in dimension $0$.
}
}

\blockn{In \cite{makkai2005word}, Makkai proved some interesting properties of his content function (as well as the support function associated) that does not seems to follow formally from the link with the linearization functor. The following are reformulation of these properties in our language.}

\block{\Prop{ $v \in \Z X$ is said to be Makkai-positif if its expression in the basis of the $m_x$ has non-negative coefficients. One says that $v \leqslant_M w$ if $w - v$ is Makkai-positive.

\begin{itemize} 

\item For any $f \in X^*$, $\delta(x)$ is Makkai-positif.

\item For any $f \in X^*$, $\delta(\pi_k^{\epsilon} f) \leqslant_M \delta(f)$

\item For any $f,g \in X^*$ composable, $\delta(f),\delta(g) \leqslant_M \delta(f \#_k g )$.

\item For any $f \in X^*$, the set of $x \in X$ such that $m_x \leqslant_M \delta(f)$ is a sub-polygraph of $X$. It is the smallest sub-polygraph of $X$ which contains $f$.

\end{itemize}

}

Aside from this proposition, the general role and the ``raison d'être'' of this notion of positivity in the linearization of an $\infty$-category is still very mysterious to us. Moreover we will introduce in \ref{Def_alternate-positive} another notion of positivity (alternate-positivity) which will be extremely important for the theory of regular polygraphs, and which seems to have absolutely nothing to do with Makkai-positivity. It should be also noted that the second property of this proposition show that in general the $\delta(f)$ for $f \in X^*$ satisfies a stronger positivity condition than Makkai-positivity: one also have that $(1 - \pi_k^{\epsilon}) \delta(f) \geqslant_M 0$ for all $k,\epsilon$, which is not true for a general Makkai-positive elements. For example, $m_f$ for $f$ a $1$-cell satisfies $(1-\pi_0^{+}) m_f = - m_{\pi_0^- f}$.

\Dem{Within the proof one will just write $\leqslant$ for the order relation defined in proposition.

We will prove, by induction on $n$, that the first three properties holds if $f,g$ are arrows of dimensions smaller than $n$. The last property will be proved afterwards.

In dimension $0$, $\delta(f)=\delta_f=m_f$ is indeed positive and the two other properties are tautologies. Assume that all three properties holds for any arrows of dimension $<n$. We will prove the first two properties by induction on composition (for the arrow $f$), and the third will follow along the way. Note that for $k \geqslant n$, $\pi_k^{\epsilon} f =f$ so we only care about the case where $k<n$.

If $x$ is a generator, then $\delta(x)=m_x + \delta(\pi_{n-1}^+ x) + \delta(\pi_{n-1}^- x)$ is indeed positive as the terms $\delta(\pi^{\epsilon}_{n-1} x)$ are. For the same reason $\delta(x) \geqslant \delta(\pi_{n-1}^{\epsilon} x)$ and for any $k <n$: 

\[ \delta(\pi_k^{\epsilon} x ) = \delta(\pi_k^{\epsilon} \pi_{n-1}^{\epsilon} x ) \leqslant \delta(\pi_{n-1}^{\epsilon} x ) \leqslant \delta(x) \]

If $f = g \#_k h$ with $g$ and $h$ satisfying the first two properties, then:

 \[ \delta(f) = \delta(g)+\delta(h)- \delta(\pi_k^+ g) = \delta(g)+\delta(h)- \delta(\pi_k^- h)  \]

as $\delta(\pi_k^{\epsilon} h) \leqslant_M \delta(h)$ and $\delta(\pi_k^{\epsilon} g) \leqslant_M g$ this already implies that $\delta(f) \geqslant_M 0$ and $\delta(g),\delta(h) \leqslant \delta_M(f)$. We now want to show that $\delta(\pi_i^{\epsilon} f) \leqslant_M \delta(f)$, one needs to treat separately the case $i \leqslant k$ and $k <i$:

If $i \leqslant k$ and $\epsilon = +$ (the case $\epsilon=-$ being exactly similar):

\[ \delta(\pi_i^+ f) = \delta(\pi_i^+ h ) \leqslant \delta(h) \leqslant \delta(f) \]
If $i > k$ :

\begin{multline*} \delta (\pi_i^{\epsilon} f) = \delta (\pi_i^{\epsilon} g \#_k \pi_i^{\epsilon} h) = \delta(\pi_i^{\epsilon} g) + \delta(\pi_i^{\epsilon} h) - \delta(\pi_k^+ \pi_i^{\epsilon} g) \\ \leqslant \delta(g)+ \delta(h) - \delta(\pi_k^+ g) = \delta(f) \end{multline*}

This concludes the proof by induction on the first three points of the proposition.

Before moving to the fourth point, One should note two things: Saying that $m_x \leqslant \delta(f)$ is the same as saying that $m_x$ appears non-trivially in the decomposition of $f$ in the $(m_x)$ basis. And if $V \subset X$ is a sub-polygraph then the $m_v$ of $\Z V$ and the $m_v$ of $\Z X$ for $v \in V \subset X$ identifies through the inclusion of $\Z V \rightarrow \Z X$. Hence any arrow $f \in V^*$ is decomposed in the $m_v$ for $v \in V$. In particular, any polygraph $W$ such that $f \in W^*$ have to contain all the $x \in X$ such that $m_x$ appears in the decomposition of $f$ in the $m$-basis. Hence if $\{x | m_x \leqslant \delta(f) \}$ is a sub-polygraph which contains $f$, it is immediate that it is the smallest.

We now prove by induction on an arrow $f$ that $\{ x | m_x \leqslant f \}$ is indeed a sub-polygraph containing $f$. For $f$ of dimension $0$, $\delta(f)=m_f$, so only $m_f \leqslant \delta(f)$ and $\{f\}$ is indeed the smallest sub-polygraph of $X$ which contains $f$. Assume the result holds for $f$ of dimension $<n$. Let $f$ be a generator of dimension $n$, then:

 \[ \delta(f)=m_f + \delta(\pi^+_{n-1}f) + \delta( \pi^-_{n-1} f) \]

hence $m_x$ appears in $\delta(f)$ if and only if it appears in one of these three terms. So the set of $x$ such that $m_x \leqslant \delta(f)$ is formed of: a sub-polygraph containing $\pi^+_{n-1} f$ union a subpolygraph containing $\pi^-_{n-1} f$ union $\{f\}$. But the union of two sub-polygraphs is a sub-polygraph, and adding a $n$-cell $f$ to a sub-polygraph $V$ such that both $\pi_{n-1}^{\epsilon} f \in V^*$ gives a polygraph containing $f$. 

Finally if $f = g \#_k h$, then I claim that $m_x$ appears in $\delta(f)$ if and only if it appears in $\delta(g)$ or in $\delta(h)$. Indeed as $\delta(f) = \delta(g)+ \delta(h) - \delta(\pi_k^{+} g)$ a $m_x$ appearing (positively) in $\delta(f)$ have to appear in $\delta(g)$ or $\delta(h)$, and conversely as $\delta(g),\delta(h) \leqslant \delta(f)$ , if $m_x$ appears in $g$ or $h$ then it appears in $f$. Hence $\{x | m_x \leqslant \delta(f) \}$ is the union of a polygraph containing $g$ and a polygraph containing $h$, and hence it is a polygraph containing $f$.

}
}

\blockn{The last part of this proposition above shows that $\delta$ can be used to compute whether a given subset of $X$ is a subpolygraph of $X$ or not and whether a cell $f \in X^*$ belong to some polygraph or not. This characterization can be formulated more explicitly, and without mentioning the $m_x$-basis as follows:}

\block{\label{Prop_deltaSeeSubPolygraph}\Prop{
\begin{itemize}

\item Let $X$ be a polygraph and $Y \subset X$ a sub-polygraph, and let $v$ an arrow of $X^*$ then:

 \[ \delta(v) \in \Z Y \Leftrightarrow v \in Y^* \]

\item Let $P$ be a subset of cell of a polygraph $X$. Then $P$ is a sub-polygraph if and only if for all $p \in P$ of dimension $n$, $\delta(\pi^{+}_{n-1} p)$ and $\delta(\pi^-_{n-1} p)$ are linear combination of the $\delta_q$ for $q \in P$.

\end{itemize}

}
\Dem{For the first point, if $v \in Y^*$ then clearly $\delta(v) \in \Z Y$. Conversely, if $\delta(v) \in \Z Y$, then as the $m_y$ form a basis of $\Z Y$, $\delta(v)$ has a decomposition in the basis of the $m_y$ for $y \in Y$. This implies that the ``smallest polygraph containing $v$'' is included in $Y$.

The second point follows immediately from the first point: one simply proves by induction on $k$ that the set of cells of $P$ of dimension at most $k$ form a sub-polygraphs, by checking, using the first condition, that the source and target of each $k$-cell is in the polygraph of cell of dimension at most $k-1$.

}
}

\subsection{The $\delta$ counting function on regular polygraphs}
\label{subsec_delta_counting_regular}

\block{\label{Def_sigma}\Def{Let $X$ be a finite polygraph. One denotes by $\sigma_X \in \Z X$ the element:

\[ \sigma_X = \sum_{x \in X} (-1)^{\text{dim}(x)} x \]

}

In order to simplify the notation:

\begin{itemize}

\item One will sometimes omit the index $X$ in $\sigma_X$.
\item One will write $(-1)^x$ or $|x|$ instead of $(-1)^{\text{dim} x}$. We will also sometimes use some ``operations'' on those cells as exponent, i.e. we will write $(-1)^{n-x-y}$ instead of $(-1)^{n-\text{dim}(x)-\text{dim}(y)}$.

\end{itemize}

}

\blockn{The reason why $\sigma$ will be extremely interesting is because in good cases\footnote{not always unfortunately as examples \ref{CE1_PolyplexNotSatisfyingTheSigmaCond} and \ref{CE2_Regular_Generic_not_sigma} will show.}, generic morphisms will preserves $\sigma$, while on the other hand for polygraphic morphisms, preserving $\sigma$ is the same as being an isomorphisms. So $\sigma$ will help us characterize generic morphisms and generic factorization in a lot of situations in terms of the behavior of morphisms on linearization.}

\block{\Lem{Let $f :X \rightarrow Y$ be a polygraphic morphism. Then:

\[ f \sigma_X = \sigma_Y \Leftrightarrow f \text{ is an isomorphism.} \]

}
\Dem{The coefficient of $f \sigma_X$ on a $y \in Y$ is (up to the sign) the number of pre-image of $y$ by $f$. So $f \sigma_X = \sigma_Y$ exactly says that any element of $Y$ has a unique pre-image by $f$, i.e. that $f$ is a bijection on cells.}
}

\block{\label{Lem_pushsigma_1}\Lem{In a polygraph $X$, let $z=x \#_k y$ be a composed arrow, and denote $b=\pi^+_k(x)=\pi^-_k(y)$. Then for any $i$ one has:

\[\begin{array}{c c c l } \delta(\pi^+_i(z))+ \delta(\pi^-_i(z))&  =& & \left[ \delta(\pi^+_i(x))+ \delta(\pi^-_i(x)) \right] \\ & & + & \left[ \delta(\pi^+_i(y))+ \delta(\pi^-_i(y)) \right] \\ & &  - & \left[ \delta(\pi^+_i(b))+ \delta(\pi^-_i(b)) \right] \end{array} \]

}

\Dem{There are three cases to treat depending on the relative values of $k$ and $i$:

\begin{itemize}

\item If $i>k$, then $\pi^+_i(z)=\pi^+_i(x) \# \pi^+_i(y)$ hence $\delta(\pi^+_i(z)) = \delta(\pi^+_i(x))+\delta(\pi^+_i(y))-\delta(\pi^+_k(\pi^+_i(x))$ but $\pi^+_k(\pi^+_i(x))=\pi^+_k(x)=b$ and $b=\pi^+_i(b)$. The same is true for $\pi^-_i$ and hence the ``additivity property'' holds for $\pi^+_i$ and $\pi^-_i$ separately and hence as well for their sum.

\item If $i<k$ then $\pi^+_i(z)=\pi^+_i(x)=\pi^+_i(y)=\pi^+_i(b)$ hence $\delta(\pi^+_i(z))= \delta( \pi^+_i(x))+ \delta(\pi^+_i(y))- \delta(\pi^+_i(b))$ and similarly, $\pi^-_i(z)=\pi^-_i(x)=\pi^-_i(y)=\pi^-_i(b)$ hence  $\delta(\pi^-_i(z))= \delta( \pi^-_i(x))+ \delta(\pi^-_i(y))- \delta(\pi^-_i(b))$. Here again the ``additivity property'' holds for $\pi^+_i$ and $\pi^-_i$ separately and hence as well for their sum.

\item If $i=k$ then $\pi^+_i(z)=\pi^+_i(y)$, $\pi^-_i(z)=\pi^-_i(x)$ and $\pi^-_i(b)=\pi^+_i(b)=\pi^+_i(x)=\pi^-_i(y)$ hence one can simplify some terms in the right hand side of the equality in the lemma to obtain that it is equal to just $ \delta(\pi^-_i(x))+\delta(\pi^+_i(y))$ which is equal to $\delta(\pi^-_i(z))+\delta(\pi^+_i(z))$ (in this last case the additivity property is not satisfied by $\delta(\pi^+_i(z))$ and $\delta(\pi^-_i(z))$ separately).

\end{itemize}
}
}

\block{\label{Lem_pushsigma_2}\Lem{Let $X$ be a polygraph and $f$ be an $n$-arrow of $X^*$. One denotes by $f:D_n \rightarrow X^*$ the corresponding morphism, and for $m>n$ let $\widetilde{f}:D_m \rightarrow X$ be the morphism representing $f$ seen as an $m$-arrow. Then:

 \[ f\left( \sigma_{D_n} \right) = \widetilde{f} \left( \sigma_{D_m} \right) . \]
}

\Dem{It is enough to show it for $m=n+1$ and to iterate. In this case,

\[ \widetilde{f} \left( \sigma_{D_{n+1}} \right) = (-1)^{n+1} \delta(f) + (-1)^{n} (\delta(f)+\delta(f)) + \sum_{k=0}^{n-1} (-1)^k (\delta(\pi^-_k(f))+ \delta(\pi^+_k(f)) \]

but the first three terms simplifies into just $(-1)^n \delta(f)$ and hence we obtain that this is just:

\[(-1)^{n} \delta(f) + \sum_{k=0}^{n-1} (-1)^k (\delta(\pi^-_k(f))+ \delta(\pi^+_k(f)) \]

which is the definition of $f^{\delta} \sigma_{D_n}$.
}
}

\block{\label{Prop_pushsigma}\Prop{\begin{enumerate} 
\item Let $\underline{p}$ be a regular $n$-polyplex, and let $p:D_n \rightarrow \underline{p}$ be the morphism representing its universal arrow, then:

\[ p\left( \sigma_{D_n} \right) = \sigma_{\underline{p}} \]

with $\sigma$ as defined in \ref{Def_sigma}.

\item Let $\underline{p}$ be a regular $n$-polyplex with spherical boundary, and $p \in \underline{p}^*$ its universal arrow, then:

\[ \delta(p) = (-1)^n \sum_{x \text{ inner cell of } \underline{p}} (-1)^{x} \delta_x \]
\end{enumerate}
}

\Dem{We will first prove, by induction on $n$, that for a fixed $n$-polyplex $\underline{p}$ with spherical boundary, if the first point holds for $\underline{p}$ and for all $\pi^{\epsilon}_k \underline{p}$ for $k<n$ then the second point hold for $\underline{p}$. 

Note that in dimension $0$, $\delta(p)=\delta_p$ is indeed the (signed) sum of all inner cell. Assume the claim holds in dimension $<n$ and let $\underline{p}$ be a regular $n$-polyplex with spherical boundary satisfying our assumption.

A cell of $\underline{p}$ appears for a unique $k$ as an inner cell of $\pi^-_k(\underline{p})$ or $\pi^+_k(\underline{p})$ hence if one denotes by $In(\underline{p})$ the set of inner cell of a polyplex $\underline{p}$ one see that the set of all cells of $\underline{p}$ is $ In(\underline{p}) \coprod In(\pi^+_{n-1}\underline{p}) \coprod In(\pi^-_{n-1} \underline{p}) \coprod \dots \coprod In(\pi^-_0(\underline{p}))$, hence one can decompose the definition of $\sigma_{\underline{p}}$:

\[ \sigma_p = \sum_{x \in In(\underline{p}) }( (-1)^{x} x ) + \sum_{k=0}^{n-1} \left[ \sum_{x \in In(\pi^+_k \underline{p}) } ( (-1)^{x} x ) + \sum_{x \in In( \pi^-_k \underline{p}) } ( (-1)^{x} x )  \right] \]

using our induction hypothesis, and the fact that all the $\pi^{\epsilon}_k \underline{p}$ satisfies our assumptions, one can replace each $\sum_{x \in In(\pi^{\epsilon}_k \underline{p}) } ( (-1)^{x} x ) $ by $(-1)^k \delta(\pi^{\epsilon}_k p)$ and one gets that:

\[ \sigma_{\underline{p}} =  \sum_{x \in In(\underline{p}) }( (-1)^{x} x )  + \sum_{k=0}^{n-1} (-1)^k (\delta(\pi^-_k(p))+ \delta(\pi^+_k(p)) ) \]

finally using the formula of the first point for $\underline{p}$ itself one has that:

\[ \sigma_{\underline{p}} = p \sigma_{D_n} = (-1)^n \delta(p) + \sum_{k=0}^{n-1} (-1)^k (\delta(\pi^-_k(p))+ \delta(\pi^+_k(p)) ) \]

taking the difference of the two formulas concludes the proof.

We now prove the first point by induction on polyplexes. The result is trivial in dimension $0$. If the result holds for all cell of dimension $<n$ and $\underline{a}$ is a regular plex, all the $\pi^-_k \underline{a}$ and $\pi^+_k \underline{a}$ are then regular polyplexes with spherical boundary for which the $1.$ has already been proved, hence, by the first half of the proof, one can apply the formula of point $2.$ to them. Hence one gets that if $a:D_n \rightarrow \underline{a}$ is the morphism representing the universal cell of $\underline{a}$ then:

\[ a(\sigma_{D_n}) = (-1)^n\delta(a)+\sum_{k=0}^{n-1} (-1)^k \left( \delta(\pi^+_k(a)) +\delta(\pi^-_k(a)) \right) \]

But $\delta(a)=\delta_a$ and, as explained above, all $(-1)^k\delta(\pi^+_k(a))$ and $(-1)^k \delta(\pi^-_k(a))$ can be replaced by the sum over their inner cells of $(-1)^x x$. Hence using the same decomposition of the set of all cells of $\underline{a}$ into inner cells of the source and target as above, one obtains the formula for $\underline{a}$.

It remains to prove that if $x$ and $y$ are composable regular polyplexes, and if $\underline{x},\underline{y}$ and $\underline{b}=\pi^+_k \underline{x}=\pi^-_k \underline{y}$ satisfies the formula of the proposition then their composite $\underline{z} =\underline{x} \#_k \underline{y}$ also satisfies it. One denotes by $x,y,z$ and $b$ the corresponding morphism $D_n \rightarrow \underline{z}$, lemma \ref{Lem_pushsigma_2} show that as long as we only care about the value of $x( \sigma)$ (which is our case) the precise value chosen for $n$ for each of these does not matter, so we just chose the same (large enough) value of $n$ for all of them. One start by writing

\[ z( \sigma) = (-1)^n \delta(z)+\sum_{k=0}^{n-1} (-1)^k \left( \delta(\pi^+_k(z)) +\delta(\pi^-_k(z)) \right) \]

But one can replace $\delta(z)$ by $\delta(x)+\delta(y)-\delta(b)$ and each terms of the sum admit a similar decomposition using \ref{Lem_pushsigma_1} into:

\[ \left[\delta(\pi^+_k(x))+\delta(\pi^-_k(x)) \right]+\left[\delta(\pi^+_k(y))+\delta(\pi^-_k(y)) \right] - \left[\delta(\pi^+_k(x))+\delta(\pi^-_k(x)) \right]  \]

Which immediately gives is:

\[ z(\sigma)= x( \sigma) + y( \sigma ) - b( \sigma) \]

and hence by the induction hypothesis:

\[ z( \sigma ) = \sigma_{\underline{x}} + \sigma_{\underline{y}} - \sigma_{\underline{b}}\]

where $\sigma_{\underline{x}},\sigma_{\underline{y}}$ and $\sigma_{\underline{b}}$ actually denotes their image in $ \Z \underline{z}$. But as $\underline{z}$ is a regular polygraph, proposition, \ref{prop_regular_mono} implies that $\underline{z}$ is the union of $\underline{x}$ and $\underline{y}$ along their subobject $\underline{b}$, and hence $\sigma_{\underline{z}} = \sigma_{\underline{x}}+ \sigma_{\underline{y}}-\sigma_{\underline{b}}$ and this concludes the proof.

}
}

\block{\label{Prop_delta_sigma_Charac_polyplex}\Prop{Let $p$ be a finite regular polygraph and $f$ an arrow of $p^*$. We also denote by $f$ the corresponding map $D_n \rightarrow p$. The following conditions are equivalent:

\begin{itemize}

\item $p$ is (isomorphic to) a regular polyplex with universal arrow $f$.

\item $f \sigma_{D_n} = \sigma_p$.

\end{itemize}

with $\sigma$ as defined in \ref{Def_sigma}.

}

It follows from this proposition that if $f$ is an arbitrary arrow in $X^*$ for $X$ a regular polygraph, and $f:D_n \rightarrow X$ is the corresponding map, then computing $f \sigma_{D_n}$ and then alternating the sign of odd dimensional cell, gives the values of the counting function described in the introduction of this section, which count the number of time each cells of $X$ appears in the image of the polyplex representing $f$.

\Dem{One implication is given by proposition \ref{Prop_pushsigma}. Conversely, assume that $f^{\delta} \sigma_{D_n} = \sigma_p$. Let $\underline{f}$ be the polyplex representing $f$, endowed with the corresponding polygraphic map $\chi_f:\underline{f} \rightarrow p$, and let $\widetilde{f}:D_n \rightarrow \underline{f}$ the morphism representing the universal cell of $\underline{f}$, in particular one has that $\chi_f \circ \widetilde{f} = f$.

Because of proposition \ref{Prop_pushsigma} $\widetilde{f} \sigma_{D_n} =\sigma_{\underline{f}}$ hence $\chi_f \sigma_{\underline{f}} = f \sigma_{\D_n} = \sigma_p$. But $\chi_f$ being polygraphic, this means that $\chi_f$ is an isomorphisms, i.e. that $p$ was isomorphic to a polyplex.

}
}

\block{\label{Def_alternate-positive}\Def{Let $X$ be a polygraph. An element of $\Z X$ is called ``alternate-positive'' if it is a linear combination with positive coefficients of elements of the form:

\[(-1)^{x} \delta_x \]

for $x \in X$. I.e. if it is a linear combination of generators where odd dimensional generator appears with negative (or zero) coefficient and even dimensional generator appears with positive (or zero) coefficient.
}

The typical example of ``alternate-positive'' element being obviously $\sigma_X$.

}

\block{\label{Prop_Sigmacriterion_boundary}\Prop{Let $\underline{p}$ be a regular $n$-polyplex, with $p:D_n \rightarrow \underline{p}$ representing the universal arrow. The following are equivalent:

\begin{itemize}

\item $\underline{p}$ has spherical boundary.

\item The image by $ p :  \Z D_n \rightarrow \Z \underline{p} $ of all the generators of $D_n$ have disjoint support in the canonical basis of $\Z \underline{p}$.

\item $p: \Z D_n \rightarrow \Z \underline{p} $ send alternate-positive elements to alternate-positive elements.

\end{itemize}
}

\Dem{The equivalence between the second and the third point follows directly from the fact that $p \sigma_{D_n} = \sigma_{\underline{p}}$. Indeed, this means that $\sigma_{\underline{p}}$ is the sum of all the $(-1)^a p(\delta_a)$ for $a$ a generator of $D_n$. But if all the $p(\delta_a)$ are pairwise disjoint, the $(-1)^a p(\delta_a)$ also are, and hence they must already have the correct coefficient to appears in $\sigma_{\underline{p}}$. In particular, they are alternate-positive, which shows that $p$ send alternate-positive elements to alternate-positive elements. Conversely, if all the $(-1)^a p(\delta_a)$ are alternate positive, then there can be no cancellation when forming their sum, and as all the coefficients of the sum ($\sigma_{\underline{p}}$) are $1$ or $-1$ it means that all the $(-1)^a p(\delta_a)$ must be disjoint.

The fact that the first condition implies the other two follow directly from the second claim in \ref{Prop_pushsigma}.

It remains to prove that the last two claim implies the first. We will prove it by induction on dimension. A $0$-dimensional or $1$-dimensional polyplex always have a spherical boundary so there is nothing to prove in these dimension. Assume this implication holds for polyplex of dimension $<n$ and let $p:D_n \rightarrow \underline{p}$ be a $n$-polyplex satisfying the second condition. The source and target of $\underline{p}$ also satisfies this second condition, and hence the source and target of $p$ have spherical boundary. In particular the second claim in \ref{Prop_pushsigma} shows that the $\delta(\pi_{n-1}^{\epsilon} p)$ are given by the (signed) sum of their inner cell, and our assumption show that they are disjoint. In particular, $\pi^-_{n-1} \underline{p}$ and $\pi^+_{n-1} \underline{p}$ share no inner cell, so that 
\[\pi^-_{n-1} \underline{p} \wedge \pi^-_{n-1} \underline{p} = \pi_{n-2}^- \underline{p} \cup \pi_{n-2}^+ \underline{p}. \]

As $\pi^-_{n-1} \underline{p}$ and $\pi^+_{n-1} \underline{p}$ are already known to have spherical boundary this proves that $p$ has spherical boundary because of corollary \ref{Cor:check_regular_with_inter}.

}

}

\block{\Def{If $X$ is a regular polygraph, an arrow $f \in X^*$ will be called regular if the polyplex representing $f$ has spherical boundary. A morphisms $f:X^* \rightarrow Y^*$ with $X$ and $Y$ regular polygraph will be called a regular morphisms if $f$ send generating cells to regular arrows. }

it corresponds exactly to the notion of $C$-admissible arrows and $\C$-morphisms of sub-section \ref{subsec_algebraicclass_def} for $\C$ the class of regular polygraphs.

}

\block{\label{Prop_GenRegularPresSigmaAndPos}\Prop{Let $f :X^* \rightarrow Y^*$ be a generic and regular morphism between regular polygraphs then:

\[ f \sigma_X = \sigma_Y \]

and $f$ send alternate-positive elements to alternate-positive elements, and all the $f \delta_a$ have disjoint support in the $\delta$-basis.

Conversely, if $f:X^* \rightarrow Y^*$ is a regular morphism between regular polygraphs which satisfies $f \sigma_X = \sigma_Y$ then $f$ is generic.

}

\Dem{For the first claim, we will proceed by induction on the number of cells of $X^*$. If $X$ has only one cell, then a generic map $X^* \rightarrow Y^*$ is an isomorphism and the result is immediate.

We now assume that the conclusion holds for a generic and regular morphism $f:X^* \rightarrow Y^*$. And we consider an extension of this into a generic and regular morphism $f':X'^* \rightarrow Y'^*$ where $X'$ is obtained by adding a new $n$-cell $a$ to $X$, whose boundary is $(u,v):\partial D_n \rightarrow X$.

let $\chi:X'^* \rightarrow (\Pb^{reg} 1)^*$ corresponding to $f'$, i.e. the composite of $f'$ with the unique polygraphic map $Y' \rightarrow \Pb^{reg} 1$ (the restriction of $\chi$ to $X$ corresponding to the map associated to $X^* \rightarrow Y^*$). This map $\chi$ has to send cell of $X$ to a polyplex (an arrow of $(Pb^{reg}1)^*$) with spherical boundary. Following the proof of theorem \ref{Th_generic_factorization} (in \ref{Proof_Th_generic_factorization}) we know that $Y'$ is constructed as the following pushout:

\[ Y' = Y \coprod_{B} Q \]

where $B$ is obtained as the generic factorization of $(\chi(u),\chi(v)) : \partial D_n \rightarrow (\Pb^{reg} 1)^*$ and $Q$ is the generic factorization of $\chi(a):D_n \rightarrow (\Pb^{reg} 1)^*$. Note that because we are assuming that $\chi(a)$ is represented by a polygraph with spherical boundary, the morphism $B \rightarrow Q$ is a monomorphism.

The map from $X'^*$ to $Y'^*$ is given by the fact that:

 \[ X'^*= X^* \coprod_{\partial D_n} D_n \]

and the maps $f:X^* \rightarrow Y^*$, $\partial \alpha : \partial D_n \rightarrow B^*$ and $\alpha: D_n \rightarrow Q^*$. The map $f'$ send any generator of $\Z X$ to its image in $\Z Y \subset \Z Y'$, and send the generator $a\in X'$ to the image by $Q \rightarrow Y'$ of $\alpha \in \Z Q$.

But by the second point of propositions \ref{Prop_pushsigma} as $Q$ is a polyplex with spherical boundary, $\alpha \in \Z Q$ is the sum on all inner cells of $Q$ of $(-1)^x \delta_x$, but the inner cell of $Q$ are send injectively to $Y'$ and corresponds exactly the cell of $Y'$ that where not in $Y$, so this concludes the proof of the first half of the proposition..

For the ``converse'', let $f :X^* \rightarrow Y^*$ be a regular morphism between regular polygraphs which satisfies $f \sigma_X = \sigma_Y$. We consider the generic factorization $f:X^* \overset{g}{\rightarrow} Y'^* \overset{u^*}{\rightarrow} Y^*$. As polygraphic map detects and preserve regular arrows (they does not change the shape of the representing polyplexes) the regularity of $f$ implies that $g$ is also regular. In particular the first half of the proposition shows that $g \sigma_X = \sigma_{Y'}$ and hence $u \sigma_{Y'} = f \sigma_X =\sigma_Y$, but a polygraphic map preserving $\sigma$ must be an isomorphisms, which proves that $f$ is generic.
}
}

\block{\label{Cor:RegularPolygraphAreAnAlgebraicClass}\Cor{Let $f :X ^* \rightarrow Y^*$ be a regular morphism between regular polygraph. Then $f$ send regular arrow of $X^*$ to regular arrow of $Y^*$. I.e. regular polygraph are an algebraic class of polygraphs. }

\Dem{ Let $f$ be a morphism as above, $a \in X^*$ a regular arrow, and $a:D_n \rightarrow X^*$ the corresponding morphism. One considers the following diagram of generic factorization:

\[
\begin{tikzcd}[ampersand replacement=\&]
 \& X'^* \arrow{r}{g_2} \arrow{d}{u^*} \& Y'^* \arrow{d}{v^*} \\
D_n \arrow{ur}{g_1} \arrow{r} \& X^* \arrow{r} \& Y^* \\
\end{tikzcd}
\]

$X'$ is a regular polyplex with spherical boundary by assumption. The generic morphism $g_2$ send generators of $X'$ to arrow represented by polyplexes with spherical boundary. Indeed, given a generator $x' \in X'$, its image in $Y^*$ is by assumption represented by a polygraph with spherical boundary, but as $v^* : Y'^* \rightarrow Y^*$ is polygraphic, the image of $x'$ in $Y'^*$ is represented by the same polyplex.

In particular, both $g_2$ and $g_1$ are generic and regular. Hence, because of proposition \ref{Prop_GenRegularPresSigmaAndPos} it means that their action on linearization preserve $\sigma_{\bullet}$ and alternate-positivity, hence their composite also preserves $\sigma_{\bullet}$ and alternate-positivity. In particular the map $g_2 \circ g_1$ is a generic morphism $D_n \rightarrow Y'^*$ whose action on the linearization preserve alternate-positivity. Finally, propositions \ref{Prop_delta_sigma_Charac_polyplex} and \ref{Prop_Sigmacriterion_boundary} show that $D_n \rightarrow Y'^*$ makes $Y'$ into a polyplex with spherical boundary, but this is the polyplex representing the image by $f$ of the arrow $a \in X^*$. Hence $f$ send regular arrows to regular arrows.
}

}

\blockn{We now give some counter-examples:}

\block{\label{CE1_PolyplexNotSatisfyingTheSigmaCond}The following example show that the first point of proposition \ref{Prop_pushsigma} fails for positive but non-regular polyplex:

Let $X$ be the $2$-polyplex:

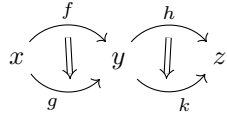
\[ \begin{tikzcd}[ampersand replacement=\&]
x \arrow[bend left=50]{r}{f}[name=U,below]{}
\arrow[bend right=50]{r}[name=D]{}[swap]{g} \&
y \arrow[Rightarrow,to path=(U) -- (D)]  \arrow[bend left=50]{r}{h}[name=V,below]{}
\arrow[bend right=50]{r}[name=W]{}[swap]{k} \& z \arrow[Rightarrow,to path=(V) -- (W)]{}
\end{tikzcd} \]

And let $Y$ be the $3$-plex corresponding to a cell for $X$ to $X$. More precesely,  $Y$ has:

\begin{itemize}

\item Three $0$-cell $x,y,z$.

\item Four $1$-cell $f,g : x \rightarrow y$ and $h,k:y \rightarrow z$.

\item Four $2$-cells $\alpha,\beta:f \Rightarrow g$ and $\gamma,\epsilon:h \Rightarrow k$.

\item One three cell $\Omega: \alpha \#_0 \gamma \rightarrow \beta \#_0 \epsilon$.

\end{itemize}

Let $\Omega:D_3 \rightarrow Y$ the morphism corresponding to the unique $3$-cell of $Y$. We claim that:

\[ \Omega(\sigma_{D_3}) = -\Omega + \alpha + \beta +\epsilon +\delta - f -g -h-k +x +z = \sigma_Y - y \]

(we have removed all the ``$\delta$'' to simplify notation) 

Hence this is an example of a (non regular) polyplex for which the map $D_n \rightarrow Y$ corresponding to the universal arrow does not satisfies $f \sigma = \sigma$ as it is the case for regular polyplex according to proposition \ref{Prop_pushsigma}.

Conversely, let $Y'$ be as $Y$ but with the $0$-cell $x$ and $y$ identified. Then the composite $f':D_3 \rightarrow Y \rightarrow Y'$ satisfies $f^{\delta} \sigma_{D_3} = \sigma_{Y'}$ but is not a polyplex (is not generic). So this gives a counter-example to the reverse implication as well.

}

\block{\label{CE2_Regular_Generic_not_sigma}Given proposition \ref{Prop_pushsigma} and \ref{Prop_GenRegularPresSigmaAndPos} it might be tempting to conjecture that for morphisms between regular polygraphs, being generic is equivalent to preserving $\sigma$. This indeed holds either when the domain is $D_n$ or when the morphism is regular, but not in general as the following examples will show:

Let $X$ be the polygraph with:

\begin{itemize}

\item One $0$-cell $x$.

\item One $1$-cell $f:x \rightarrow x$.

\item One $2$-cell $\alpha:f \Rightarrow f$.

\end{itemize}

It is a regular (even globular) polygraph. We consider the map $X \rightarrow \Pb^{reg} 1$ defined as: $x$ goes to $\bullet$, $f$ goes to $\bullet \rightarrow \bullet \rightarrow \bullet$ and $\alpha$ goes to:

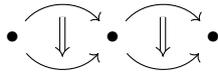
\[ \begin{tikzcd}[ampersand replacement=\&]
\bullet \arrow[bend left=50]{r}[name=U,below]{}
\arrow[bend right=50]{r}[name=D]{} \&
\bullet \arrow[Rightarrow,to path=(U) -- (D)]  \arrow[bend left=50]{r}[name=V,below]{}
\arrow[bend right=50]{r}[name=W]{} \& \bullet \arrow[Rightarrow,to path=(V) -- (W)]{}
\end{tikzcd} \]

It corresponds to a generic, non-regular morphism $\lambda:X \rightarrow Y$ where $Y$ is the regular polygraph which has:

\begin{itemize}

\item two $0$-cell $x$ and $t$.

\item two $1$-cell $g:x \rightarrow  t$ and $h:t \rightarrow x$.

\item two $2$-cell $\beta:g \rightarrow g$ and $\gamma:t \rightarrow t$.

\end{itemize}

with $x \in X$ being send to $x \in Y$, $f \in X$ to $g \#_0 h \in Y^*$ and $\alpha \in X$ to $\beta \#_0 \gamma$.

One has:

\[ \lambda \sigma_X = \delta(\beta \#_0 \gamma) - \delta(g \#_0 h) + \delta(x)= (\beta + \gamma - t) - ( g + h -t) +x = \sigma_Y - t \]

Hence it is an example of generic (non-regular) morphism between regular polygraphs which does not preserve $\sigma$.

For a counter example to the converse implication, one can define $Y'$ to be exactly as $Y$ but with $x$ and $t$ identified, and then the composite $\lambda': X \rightarrow Y \rightarrow Y'$ is no longer generic but does satisfy $\lambda' \sigma_{X} = \sigma_{Y'}$.
}

\block{\label{Prop_CaracGenericRegular}Finally, one might wonder if one can obtain a characterization of generic and regular morphisms between regular polygraphs which will be an actual converse to proposition \ref{Prop_GenRegularPresSigmaAndPos}. We have not been able to do so (nor to find a counterexample) but one has the following proposition:

\Prop{ Let $X$ be a regular polygraph such that for any cell $a$ of $X$ the corresponding plex $\underline{a} \rightarrow X$ is a monomorphism.

Let $f :X^* \rightarrow Y^*$ be a morphism of $\infty$-category between regular polygraph such that $f :\Z X \rightarrow \Z Y$ preserve alternate-positivity and satisfies $f \sigma_X = \sigma_Y$. 

Then $f$ is generic and regular.

}

This already cover some case of interest: when $X$ is $\partial D_n$ (this case will be used in \ref{Prop_tensor_product_of_regular}) and when $X$ is a polyplex.

Because of proposition \ref{Prop_GenRegularPresSigmaAndPos}, we only need to show that the map is regular. The proof will be by induction on the number of cells of $X$, but one need first the following lemma:
}

\block{\label{Lem_imageSubPol}\Lem{Let $f:X^* \rightarrow Y^*$ be a morphism preserving $\sigma$ and alternate positivity. Let $V$ be a subpolygraph of $X$, and $V_f$ be the set of cells of $Y$ which appears non trivially in (the $\delta$-basis decomposition of) $f(\delta_v)$ for $v \in V$. Then $V_f$ is a sub-polygraph of $Y$ and $f$ restrict to a map from $V^*$ to $V_f^*$.}

\Dem{We will use the criterion given in proposition \ref{Prop_deltaSeeSubPolygraph}: we need to show that for all generator $v \in V_f$, $\delta(sv)$ and $\delta(tv)$ are in $\Z V_f$ where $sv$ and $tv$ denote the source and target of $v$.

If $v \in V_f$ then $\delta_v$ appear non-trivially in a $f(\delta_x)$ for $x \in V$. Because $f$ preserves alternate-positivity and $\sigma$, one can write:

\[ (-1)^x f(\delta_x) = \sum_i (-1)^{a_i} \delta_{a_i} \]

with the $a_i$ pairwise distinct and $a_0=v$. We call $k+1$ the dimension of $v$, and one applies $\pi_{k}^{\epsilon}$ on both side to obtain that:

\[ (-1)^x f(\pi_k^{\epsilon}(\delta_x)) = \sum_i (-1)^{a_i} \delta(\pi_k^{\epsilon}(a_i)) \]

Let $w$ be a coefficient appearing non-trivially in $(-1)^{a_0} \delta(\pi_k^{\epsilon}(a_0))=(-1)^v \delta(\pi_k^{\epsilon}(v))$, because $Y$ is regular, the second point of proposition \ref{Prop_pushsigma} shows that $w$ appear there with coefficient $(-1)^{w}$. One distinguishes two possibilities: 

\begin{itemize}

\item Either $w$ appears non-trivially in the expression of $f(\delta(\pi_k^{\epsilon}(x)))$, hence it appears non-trivially in a $f(\delta_y)$ for $y$ appearing in $\pi_k^{\epsilon}(\delta_x)$, but such a $y$ is in $V$, and hence this proves that that $w \in V_f$.

\item Either $w$ is canceled out in the expression of $f(\delta \pi_k^{\epsilon}(x))$, so there is (at least) one $i$ such that $w$ appear with sign $-(-1)^w$ in $(-1)^{a_i} \pi_k^{\epsilon}(\delta_{a_i})$. But if $a_i$ dimension greater $>k$ then, again because of the second point of proposition \ref{Prop_pushsigma},  $(-1)^{a_i}\pi_k^{\epsilon}(\delta_{a_i})$ only contains term of the form $(-1)^x \delta_x$, so $w$ must appear in a $\pi_k^{\epsilon}( \delta_{a_i})$ for $i \leqslant k$, but this means that $w$ is itself one of the $a_i$ so in this case it already appears non-trivially in the expression of $f(\delta_x)$, hence $w \in V_f$.

\end{itemize}

In both case, we have shown that for $v \in V_f$, every generator appearing in $\delta(\pi_k^{\epsilon}(v))$ is in $V_f$ which proves that $\pi_k^{\epsilon}(v)$ (i.e. $sv$ or $tv$) is in $\Z V_f$ and hence that $V_f$ is indeed a polygraph by proposition \ref{Prop_deltaSeeSubPolygraph}. The fact that $f$ send $V^*$ to $V_f^*$ follow immediately from the first point of proposition \ref{Prop_deltaSeeSubPolygraph}.

}
}

\block{We can now prove proposition \ref{Prop_CaracGenericRegular}:

\Dem{

Because of the last claim in proposition \ref{Prop_GenRegularPresSigmaAndPos}, we only need to prove that $f$ send cell to regular arrows.
Let $a$ be any cell of $X$, and let $\underline{a} \subset X$ be the corresponding plex, which by assumption is a sub-polygraph of $X$. One can form its ``image'' in $Y$ in the sense of lemma \ref{Lem_imageSubPol}, and one obtains a subpolygraph $Y_a \subset Y$, such that the map $X^* \rightarrow Y^*$ restrict to a map $Y_a \rightarrow \overline{a}$, which satisfies the assumption of proposition \ref{Prop_CaracGenericRegular} (indeed, $Y_a$ only contains the cells of $Y$ which will appears in the image of $\sigma_{\underline{a}}$ and it preserves alternate-positivity as a restriction of a map that does.

In particular, composing with the map $D_n \overset{a}{\rightarrow} \underline{a} \rightarrow Y_a$ gives a map $D_n \rightarrow Y_a$ which preserves $\sigma$ and alternate positivity, hence, because of proposition \ref{Prop_delta_sigma_Charac_polyplex} and \ref{Prop_Sigmacriterion_boundary} this shows that $Y_a$ is a polyplex with spherical boundary, hence it is the polyplex representing the image of $a$ in $Y$ and it is indeed regular.
}

}

\subsection{On the Gray tensor product}
\label{subsec_Gray_tensor_product}

\blockn{We start with the following ``black box theorem'' containing everything we need about the Gray tensor product of $\infty$-categories (which we do not claim is new in any way). In fact, one could prove that these properties completely characterize the Gray tensor product.}

\block{\Th{There is a monoidal structure on the category of $\infty$-categories given by the the so-called Gray tensor product which satisfies the following properties:

\begin{itemize}

\item The category of strict $\infty$-categories is monoidal closed.

\item The linearization functor $C \mapsto \Z C$ from $\infty$-category to strict $\infty$-groups is monoidal.

\item There is an $(n+m)$-arrow $t_{n,m}$ of $D_n \otimes D_m$ such that under the isomorphism $\Z (D_n \otimes D_m) \simeq \Z(D_n) \otimes \Z(D_m)$,  $\delta(t_{n,m})$  is identified with $\delta(t_n) \otimes \delta(t_m)$.

Moreover the map:

\[  \partial D_n \otimes D_m \coprod_{\partial D_n \otimes \partial D_m} D_n \otimes \partial D_m \rightarrow D_n \otimes D_m\]

is obtained by freely adding the arrow $t_{n,m}$ with a specified boundary.

\end{itemize}

 }

The first point has been proved several times. The most interesting proof for our purpose is the one suggested by R.Steiner in \cite{steiner2004omega} and made complete and precise by D.Ara and G.Maltsiniotis in appendix $A$ of \cite{ara2016joint}. The proof relies on what is now often called ``Steiner's theory'' (introduced in \cite{steiner2004omega}) which show that there is a small class of strict of $\infty$-categories (containing for examples Orientals, Cubes, Globes etc.) which is equivalent to a certain category of chain complexes (a subcategory of directed augmented chain complexes to be precise), one then observe that this nice class of strict of $\infty$-categories is endowed with a tensor product coming from the tensor product of chain complexes and one can prove that this tensor product can be extended using an analogue of Day convolution to all strict $\infty$-categories. The general idea and the first part of the proof is entirely due to R.Steiner in \cite{steiner2004omega}, but his proof that the monoidal product indeed extend to all strict $\infty$-categories, contains some gap, or at least is imprecise. A complete and correct proof following this idea has been given by D.Ara and G.Maltsiniotis in \cite{ara2016joint}.

A completely different approach to define the Gray tensor product (and older) is to use instead a cubical description of strict $\infty$-category. The Gray tensor products is then the one induces by the tensor product of cubes, see \cite{crans1995pasting}.

It is important to note that this makes the Gray tensor product quite difficult to compute. The only access point we have to it is the tensor product of ``Steiner categories'' (i.e. those which can be described using Steiner theory) so only some basic diagrams.

The second point has been proved in terms of the tensor product of chain complexes in \cite[Prop. A.18]{ara2016joint}.

The last point is an easy exercise for anyone familiar with Steiner's theory: as all the $D_n$ and $\partial D_n$ are strong Steiner $\infty$-categories one can compute explicitly their tensor product as the image under Steiner's equivalences of the tensor products of corresponding chain complexes. It is shown in \cite{steiner2004omega} (see theorem 6.1) that Steiner's categories are always given by polygraphs whose generating cells are given by the basis of the corresponding complexes, so that one can easily check by hand that the map:

\[  \partial D_n \otimes D_m \coprod_{\partial D_n \otimes \partial D_m} D_n \otimes \partial D_m \rightarrow D_n \otimes D_m\]

is an inclusion of polygraphs, missing only one cell $t_{n,m}$ of dimension $n+m$, and that $\delta(t_{n,m})=\delta(t_n) \otimes \delta(t_m)$.
}

\block{\Def{Let $X$ and $Y$ be two strict $\infty$-categories, let $x\in X_n$ and $y \in Y_m$ be respectively a $n$-arrow of $X$ and a $m$-arrow of $Y$.

Let $x \otimes y$ be the $(n+m)$-arrow of $X \otimes Y$ defined as the image of $t_{n,m} \in D_n \otimes D_m$ under the tensor product
 \[ D_n \otimes D_m \overset{x \otimes y}{\rightarrow} X \otimes Y \]

}}

\block{\Prop{Under the isomorphism $\Z(X \otimes Y) \simeq \Z X \otimes \Z Y$ one has, for any cells $x \in X$ and $y \in Y$:

\[ \delta(x \otimes y) =  \delta(x) \otimes \delta(y) \]

}

\Dem{Let $x: D_n \rightarrow X$ and $y:D_m \rightarrow Y$ be the morphisms corresponding to $x$ and $y$. Then because linearization functor is monoidal one has that the square :

\[
\begin{tikzcd}[ampersand replacement=\&]
\Z (D_n \otimes D_m) \arrow{d}{(x \otimes y)} \arrow{r}{\sim}
 \& \Z D_n \otimes \Z D_m \arrow{d}{x \otimes y}
\\
\Z (X \otimes Y) \arrow{r}{\sim} \& \Z X \otimes \Z Y 
\end{tikzcd}
\]

But under the upper isomorphism $\delta(t_{n,m})$ corresponds to $\delta(t_n) \otimes \delta(t_m)$ by our black box theorem, and on the left side: 

\[ (x \otimes y)(\delta(t_{n,m})) =  \delta((x\otimes y)(t_{n,m})) = \delta(x \otimes y),\]

while on the right side:

 \[ (x \otimes y)(\delta( t_n) \otimes \delta(t_m)) = x(\delta(t_n)) \otimes y(\delta(t_m)) = \delta(x) \otimes \delta(y) \]
which gives us the result.

}
}

\blockn{The following proposition is by no way due to us: it is a ``folk result'' whose proof is well known by all experts\footnote{We can distinctly remember M.Lucas presenting a form of the following proof in the cubical setting during a talk in septembre 2017, as well as G.Maltsiniotis and D.Ara mentioning around the same time having written this proof either in unpublished notes or in a paper in preparation.} of the field; it appears as theorem 1.35 of \cite{hadzihasanovic2017algebra}. }

\block{\label{prop:TensProdOfPG}\Prop{\begin{itemize}

\item If $X$ and $Y$ are polygraphs, then $X \otimes Y$ is a polygraph whose generating cells are the $x \otimes y$ for $x$ and $y$ the generating cells of $X$ and $Y$.

\item If $f$ and $g$ are polygraphic morphisms $f :X \rightarrow X'$ and $g: Y \rightarrow Y'$ then $f \otimes g : X \otimes Y \rightarrow X' \otimes Y'$ is a polygraphic morphism too.

\end{itemize}

 }

\Dem{We first notice that the second point follows immediately from the first: $(f\otimes g)(x \otimes y) = f(x) \otimes g(y)$ hence $f \otimes g$ send generators to generators and hence is polygraphic.

We will now prove the first point for finite polygraphs, by induction on $n \times m$ for $n$ and $m$ the number of cells of $X$ and $Y$ respectively. For $n \times m = 0$, it means that one of the polygraph is the initial object and because the Gray tensor product is monoidal closed it commutes to the initial object, which proves the results. 
Assume that $X$ and $Y$ are polygraph with respectively $n$ and $m$ cells and that the result is known for all tensor product with $n' \times m'< n \times m$. Take $x$ and $y$ be some top dimensional cells of $X$ and $Y$ respectively, so that one has sub-polygraphs $X' = X-\{x\} \subset X$ and $Y' = Y-\{y\} \subset Y$. Our induction hypothesis applies to all the product of polygraphs:

\[ X' \otimes Y \quad X \otimes Y' \quad  X' \otimes Y' \]

And the argument for showing the second point apply to both the maps $X' \otimes Y' \rightarrow X \otimes Y'$ and $X' \otimes Y' \rightarrow X' \otimes Y$ making them both polygraphic inclusion. From this one deduce easily that:

\[ (X' \otimes Y) \coprod_{X' \otimes Y'} (X \otimes Y') \]

is a polygraph whose generator are given by pair $(a,b)$ with $a$ a cell of $X$, $b$ a cell of $Y$ and either $a$ different from $X$ or $b$ different from $y$, with the canonical map to $X \otimes Y$ sending the pair $(a,b)$ to $a \otimes b$.

One then recalls that $X$ is a pushout of $X'$ along $\partial D_i \hookrightarrow D_i$ and $Y$ is a pushout of $Y'$ along $\partial D_i \hookrightarrow D_i$ hence it follow formally (from the usual property of the ``pushout-product'' of maps) that the diagram:

\[
\begin{tikzcd}[ampersand replacement=\&]
 \displaystyle \partial D_n \otimes D_m \coprod_{\partial D_n \otimes \partial D_m} D_n \otimes \partial D_m \arrow{r} \arrow{d} \&  D_n \otimes D_m \arrow{d} \\
 \displaystyle (X' \otimes Y) \coprod_{X' \otimes Y'} (X \otimes Y') \arrow{r} \& X \otimes Y
\end{tikzcd}
\]

is a pushout, with the vertical arrow being given by the map ``$x \otimes y$'' and its restriction.

But our black box theorem says that the upper map is obtained by freely adding the cell $t_{n,m}$ so that $X \otimes Y$ is freely obtained from the polygraph $\displaystyle (X' \otimes Y) \coprod_{X' \otimes Y'} (X \otimes Y')$ by freely adding a single cell, whose image in $X$ is $x \otimes y$. And this complete our proof in the case of finite polygraphs. The case of general polygraph follows by just writing them as directed colimit of finite polygraphs.
}
}

\blockn{At this point we have obtained a new interesting combinatorial property of our $\delta$-counting function: if $X$ and $Y$ are polygraphs, then using the identification of the cells of $X \otimes Y$ with pair of cell $x \otimes y$ for $x$ and $y$ cells of $X$ and $Y$ to identity $\Z X \otimes \Z Y$ with $\Z (X\otimes Y)$ one has that $\delta(f \otimes g) = \delta(g) \otimes \delta(g)$. When looking at $\delta$ from the ``naive'' perspective of a counting function defined by induction on expressions, this is a very non-trivial property, which will allow us to obtain some results about the Gray tensor product of polygraphs. The next proposition does not really exploit the function $\delta$, but only the ``top dimensional part'' which lives more naturally in the chain complex and was well understood long before this work, but proposition \ref{Prop_tensor_product_of_regular} below will makes a real use of the $\delta$ function.}

\block{\Prop{The tensor product of two positive polygraphs is a positive polygraphs.}

\Dem{Let $P$ and $Q$ be two positive polygraphs, any cell of the tensor product $P \otimes Q$ is of the form $p \otimes q$ for $p$ and $q$ cells of $P$ and $Q$. Let $i$ and $j$ be the dimension of $p$ and $q$ and $n=i+j$ the dimension of $p \otimes q$ using the second special case of proposition \ref{Lem_pi_tens_degminus1} one has that:

\[ \delta( \pi_{n-1}^{\epsilon}(p \otimes q )) = \delta(\pi^{\epsilon}_{i-1} p) \otimes \delta(q) +  \delta(p) \otimes \delta( \pi^{|i|\epsilon}_{j-1}(q)) - \delta(\pi^{\epsilon}_{i-1} p) \otimes \delta( \pi^{|i|\epsilon}_{j-1}(q))  \]

The last terms contains no cell of dimension $(n-1)$ and both the first and second contains at least one generator of dimension $(n-1)$. This proves that the source and target of $p\otimes q$ is a non-identity arrow of dimension $n-1$. As this is true for all the cells of $P \otimes Q$, it proves that it is a positive polygraph.
}
}

\block{\label{Lem_tensorStabilityofPosandSigm}\Lem{\begin{itemize}

\item Let $X$ and $Y$ be two polygraphs, with $\sigma_X \in \Z X$, $\sigma_Y \in \Z Y$ and $\sigma_{X \otimes Y} \in \Z (X \otimes Y)$, then under the identification of $\Z(X \otimes Y)$ with $\Z X \otimes \Z Y$ one has:

 \[\sigma_X \otimes \sigma_Y = \sigma_{X \otimes Y} \]

\item Let $f:X \rightarrow Y$ and $g:W \rightarrow Z$ be two morphisms of polygraphs such that both the linearization of $f$ and $g$ preserve $\sigma_{\bullet}$ and/or alternate-positivity, then their tensor product $(f \otimes g) : X \otimes W \rightarrow Y \otimes Z$ also does.

\item The map $D_{n+m} \rightarrow D_n \otimes D_m$ preserves $\sigma$ and alternate positivity.

\end{itemize}}

\Dem{
\begin{itemize}

\item As the cells of $X \otimes Y$ are exactly the $x \otimes y$ for $x$ and $y$ cells of $X$ and $Y$ respectively and the dimension of $x \otimes y$ is the sum of the dimension of $x$ and $y$ one has:

 \begin{multline*}\sigma_P \otimes \sigma_Q = \left( \sum_{p \in P} (-1)^p \delta(p) \right) \otimes \left( \sum_{q \in Q} (-1)^q \delta(q) \right)\\ = \sum_{\mytop{p \in P}{q \in Q}} (-1)^{p+q} \delta(p \otimes q) = \sigma_{P\otimes Q} \end{multline*}

\item The first part of the lemma shows that:

\[ (f \otimes g) (\sigma_{X \otimes W}) = (f \otimes g)(\sigma_X \otimes \sigma_W) \]

But then:

\[(f \otimes g)(\sigma_X \otimes \sigma_W) = f \sigma_X \otimes g \sigma_W = \sigma_Y \otimes \sigma_Z= \sigma_{Y \otimes Z} \]

This proves that $(f \otimes g)$ preserves $\sigma$. We now need to show that it preserves alternate positivity. A general cell of $X \otimes W$ can be written as $p=x \otimes w$ for $x$ and $w$ cells of $X$ and $W$, and $dim(p)=dim(x)+dim(w)$. So, one has that:

\[ (-1)^p \delta_p = ((-1)^x \delta_x) \otimes ((-1)^w \delta_w)  \]

Hence if $f$ and $g$ preserves alternate-positivity, the image by $(f \otimes g)$ of $(-1)^{p} \delta_p$ is the tensor product of two alternate-positive elements. But if $y$ and $z$ are cells of $Y$ and $Z$, then $(-1)^y \delta_y \otimes (-1)^z \delta_z$ is alternate positive in $\Z(Y \otimes Z)$ for the exact same reason and hence any tensor product of alternate-positive elements is alternate positive.

\item Let $v : D_{k} \rightarrow D_n \otimes D_m$ (for $k=n+m$) the map corresponding to the top cell of $D_n \otimes D_m$. The generator of $D_k$ are denoted $[i^+]$ and $[i^-]$ for $i<k$ and $[k]=[k^+]=[k^-]$ in dimension $k$. One has that $[i^{\epsilon}] = \pi_i^{\epsilon}([k])$ hence as $v([k])=[n] \otimes [m]$ one has $v([i^{\epsilon}]) = \pi_i^{\epsilon}([n] \otimes [m])$. This last terms is in general difficult to computes, but its image in the linearization is well understood by the fact that linearization is monoidal and the description of tensor product of $\infty$-groups given in subsection \ref{subsec:TensorProdOfInfinityGroup}. In particular, using the formula of proposition \ref{Lem_pi_tens_degminus1}, in the linearization, and keeping the notation $[i^{\epsilon}]$ for the image of the generating cells by $\delta$:

\begin{align*}
v([i^{\epsilon}]) & =\pi_i^{\epsilon}([n] \otimes [m]) \\  & = \sum_{\mytop{a+b=i}{\mytop{a \leqslant n}{b \leqslant m} }}\left( \pi_a^{\epsilon} [n] \otimes \pi_b^{|a|\epsilon}[m] \right) -  \sum_{\mytop{a+b=i-1}{\mytop{a < n}{b< m }}} \left( \pi_{a}^{\epsilon} [n] \otimes \pi_b^{-|a| \epsilon} [m] \right) \\ 
& = \sum_{\mytop{a+b=i}{\mytop{a \leqslant n}{b \leqslant m} }}\left( [a^{\epsilon}] \otimes [b^{|a|\epsilon}] \right) -  \sum_{\mytop{a+b=i-1}{\mytop{a < n}{b< m }}} \left( [a^{\epsilon}] \otimes [b^{-|a| \epsilon}] \right)
\end{align*}

This already proves that $v((-1)^i [i^{\epsilon}])$ is alternate-positive: once we multiply both side by $(-1)^i$ the first sum only contains terms of dimension $i$ with sign $(-1)^i$ and the second only terms of dimension $(i-1)$ with sign $-(-1)^i$. Hence $v$ preserves alternate positivity. It remains to prove that $v$ preserves $\sigma$. One has that $v(\sigma_{D_k}))$ is the sum of $(-1)^k v([k]) =(-1)^{n+m} ([n]\otimes[m])$ and of all the $(-1)^i v([i^{\epsilon}])$ for $i<k$ and $\epsilon=+$ or $-$. Note that each of this terms is a sum of terms of the form $(-1)^x x$ for $x$ a cell of $D_n \otimes D_m$, hence what we need to do is to prove that each cell of $D_n \otimes D_m$ appears once and exactly once in this sum, we have different type of cell to treat separately:

\begin{itemize}

\item All the terms of the form $(-1)^{a+b}[a^{\epsilon}] \otimes [b^{\delta}]$ for $a<n$ and $b<m$ appears exactly once: if $\delta=|a| \epsilon$ they can only appears in the first sum of the expression of $(-1)^i v([i]^{\epsilon})$ for $i=a+b$ while if $\delta=-|a| \epsilon$ they can only appears in the second sum of the expression of $(-1)^i v([i]^{\epsilon})$ for $i=a+b+1$, and as $a <n$ and $b<m$ one indeed have $a+b+1 <m+n$.

\item The terms of the form $(-1)^{n+b}[n] \otimes [b^{\delta}]$ for $b <j$ can only appears in the first the sum of the expression of $(-1)^i v([i^{\epsilon}])$ for $i = n+b <n+m$ and $\epsilon = |n| \delta$, and so it appears exactly once.

\item The terms of the form $(-1)^{a+m}[a^{\delta}] \otimes [m]$ for $a <n$ can only appears in the first sum of the expression of $(-1)^i v([i^{\epsilon}])$ for $i = a+m <n+m$ and $\epsilon = \delta$.

\item The term $(-1)^{i+j} [i] \otimes [j]$ only appears as the single terms of $(-1)^k v([k])$.

\end{itemize}

\end{itemize}

}

}

\block{\label{Prop_tensor_product_of_regular}\Prop{\begin{itemize}

\item[(A)] The tensor product of two regular plexes is a regular plex.

\item[(B)] The tensor product of two regular polygraphs is a regular polygraph.

\item[(C)] The tensor product of two regular polyplexes is a regular polyplex, with the universal arrow being the tensor product of the universal arrow.

\item[(D)] The tensor product of two regular polyplexes with spherical boundary is a regular polyplex with spherical boundary.

\end{itemize}

}

I conjecture that more generally, in any good class of polygraphs stable under tensor product, the tensor product of two plexes is a plex and that the tensor product of two polyplexes is a polyplex as above. But I haven't been able to prove it for positive polygraph.

\Dem{We denote $A_n$ the proposition ``The tensor product of two regular plexes of dimension $i$ and $j$, with $i+j \leqslant n$ is a regular plex''. And similarly for $B_n$ , $C_n$ and $D_n$.

$A_0$ is trivial. We assume $A_n$ and we will gradually proves $B_n,C_n,D_n$ and finally $A_{n+1}$ which will concludes the proof by induction.

\bigskip

\textbf{Proof of $B_n$:} if $P \otimes Q$ is a product of two polygraphs of dimension $i$ and $j$ with $i+j \leqslant n$, then any cell of $P \otimes Q$ is of the form $p \otimes q$ for $p$ and $q$ cells of $P$ and $Q$. The cells $p$ and $q$ are represented by regular plexes $\underline{p} \overset{\chi_p}{\rightarrow} P$ and $\underline{q} \overset{\chi_q}{\rightarrow} Q$ of dimension smaller than $i$ and $j$. The cell $p \otimes q$ is the image of the top cell of the tensor product $\underline{p} \otimes \underline{q}$, but because we assumed $A_n$, $\underline{p} \otimes \underline{q}$ is a regular plex, hence $p\otimes q$ is indeed the image of a regular plex.

\bigskip

\textbf{Proof of $C_n$ and $D_n$:} Let $\underline{p}$ and $\underline{q}$ be two regular polyplexes of dimension $i$ and $j$. One can assume that $i+j=n$. We just proved that $\underline{p} \otimes \underline{q}$ is a regular polygraph, and it is endowed with an arrow $p \otimes q \in (\underline{p} \otimes \underline{q})^*$. Each of the maps $D_i \rightarrow \underline{p}$ and $D_j \rightarrow \underline{q}$ preserves $\sigma$, hence by lemma \ref{Lem_tensorStabilityofPosandSigm} the maps $D_n \rightarrow D_i \otimes D_j \rightarrow \underline{p} \otimes \underline{q}$ also does. This proves that the pair $(\underline{p} \otimes \underline{q}, p \otimes q)$ is a regular polyplex because of proposition \ref{Prop_delta_sigma_Charac_polyplex}. If one further assume that $\underline{p}$ and $\underline{q}$ have spherical boundary, it means that $D_i \rightarrow \underline{p}$ and $D_j \rightarrow \underline{q}$ preserve alternate-positivity, hence by the same argument as above, so does the map $D_n \rightarrow \underline{p} \otimes \underline{q}$ which proves that $\underline{p} \otimes \underline{q}$ has spherical boundary by proposition \ref{Prop_Sigmacriterion_boundary}.

\textbf{Proof of $A_{n+1}$:}

Let now $\underline{p}$ and $\underline{q}$ be plexes of dimension $i$ and $j$ such that $i+j=n+1$ (the case $i+j \leqslant n$ being already known). Let $p:D_i \rightarrow \underline{p}$ and $q: D_j \rightarrow \underline{q}$ be the two natural morphisms. First any cell of $\underline{p} \otimes \underline{q}$ which is not the top dimensional cell appears as a cell of a tensor product of regular polygraphs of dimension $i'$ and $j'$ with $i' +j' \leqslant n$, hence appears as a cell of a regular polygraph (i.e. is represented by a regular plex). This show that the polygraph $Y$ defined as $\underline{p} \otimes \underline{q}$ minus its top dimensional cell is already known to be a regular polygraph.

One more time, let:

\[ r:D_{n+1} \overset{v}{\rightarrow} D_i \otimes D_j \overset{p \otimes q}{\rightarrow} \underline{p} \otimes \underline{q} \]

Exactly as in the proof of $C_n$ and $D_n$, one has that $D_{n+1} \rightarrow \underline{p} \otimes \underline{q}$ preserves $\sigma_{\bullet}$ and alternate-positivity. One cannot applies proposition \ref{Prop_delta_sigma_Charac_polyplex} and \ref{Prop_Sigmacriterion_boundary} as above directly as we do not yet know that $\underline{p} \otimes \underline{q}$ is regular. But it is immediate to deduce from this that $\partial D_n \rightarrow Y$ preserve $\sigma$ and alternate positivity as the top cell of $D_n$ was sent to the only cell of $\underline{p} \otimes \underline{q}$ which is not in $Y$. This proves (for example by \ref{Prop_CaracGenericRegular}) that $\partial D_n \rightarrow Y$ is a generic and regular morphism. In particular the top cell of $\underline{p} \otimes \underline{q}$ has source and target that are regular, hence $\underline{p} \otimes \underline{q}$ is regular, and this concludes the proof.
}
}

\subsection{Join of $\infty$-categories}
\label{subsec:joinsAndCones}
\blockn{There is one last things we will need in order to prove the Simpson conjecture: some knowledge about the ``Join'' construction of $\infty$-category introduced by D.Ara and G.Maltsiniotis in \cite{ara2016joint} and its behavior with respect to our theory of plex, polyplex and regular polygraphs. This theory is extremely parallel to the theory of the Gray tensor product, and we will only make a very marginal use of it: we will only use joins of the form $* \star X$, and they will only appears in the proof of theorem \ref{Th_Spaces_and_regular_equiv} to prove that $*$ is terminal in the homotopy category of regular polygraphs. This section will also provide a nice proof using the join that Street's Orientals are plexes (proposition \ref{Prop_street_orientals}), but this can probably be achieved more directly without talking about the join. And it will also provide a very simple proof that the category of regular plexes is a weak test category, but unfortunately we have not been able to make any use of this proof to simplify the rest of the paper. For these reasons we will not give as much details about the construction of this join operation as the proof are very similar to those in section \ref{subsec_Gray_tensor_product} and that we will not use these results much.}

\block{\Def{An augmented strict $\infty$-group $G$ is a strict $\infty$-group endowed with a morphisms to the strict $\infty$-group $\Z = \Z *$.}

i.e. it is a linear map $e: G \rightarrow \Z$ such that for each $g$, $e(g) = e(\pi_0^+ g)=e(\pi_0^- g)$, which implies more generally $e(g) = e(\pi_k^{\epsilon} g)$ for all $k,\epsilon$.

As $*$ is the terminal $\infty$-category, for every $\infty$-category $C$, $\Z C$ has a canonical augmentation, which is defined by $e(\delta c)= 1$ for all $c \in C$, and any functor between $\infty$-category induce a morphism of augmented $\infty$-group.
}

\block{\Def{If $e: G \rightarrow \Z$ is an augmented $\infty$-group, one defines a new $\infty$-group $\Sigma G$, as an abelian group it is $G \oplus \Z$ with the $\pi$ operation defined as follow:

\[ \pi_k^{\epsilon}(g \oplus z ) = (\pi_{k-1}^{\epsilon}) g \oplus z \]
\[ \pi_0^{-}(g \oplus z) = 0 \oplus z \]
\[ \pi_0^{+}(g \oplus z) = 0 \oplus \left( z+e(g) \right)\]

}

One easily check that $\Sigma G$ is indeed a globular group (hence a strict $\infty$-group). Conversely, if $H$ is a strict $\infty$-group and $H_0 \simeq \mathbb{Z}$, then $\{v \in H | \pi_0^-(v) = 0 \}$ forms an augmented strict $\infty$-group with $ \pi^{\epsilon}_k (v) $ being defined as the former $\pi^{\epsilon}_{k+1}$, and with augmentation $\pi^+_0$.

These two construction induces an equivalence between the category of augmented strict $\infty$-groups and the (non full) subcategory of strict-$\infty$ groups endowed with an identification $G_0 \simeq \mathbb{Z}$ (and morphisms preserving the identification).}

\block{ As $(G \otimes H)_0 \simeq G_0 \otimes H_0$, the tensor product of strict $\infty$-groups restrict to a tensor product on this category of strict $\infty$-groups with an identification $G_0 \simeq \mathbb{Z}$. In particular, this induces a new monoidal structure denoted $\star$ on the category of augmented strict $\infty$-groups by transport along the equivalence described above. The unit of this monoidal structure is the $\infty$-group $0$.

If $G$ and $H$ are two strict $\infty$-groups, the underlying group of $G \star H$ is given by $G \oplus (G \otimes H) \oplus H$, and is functorial in $G$ and $H$ in the natural way. We will not really take time to explicitly describe globular structure on $G \star H$ as we will not need it in the present paper, we just want to says that because of the functoriality and the fact that $0$ is the unit of this tensor product, the two components $G$ and $H$ in  $G \oplus (G \otimes H) \oplus H$ are actually isomorphic copies of $G$ and $H$. An element of the terms $G \otimes H$ of the form $g \otimes h$ will be denoted $g \star h$. If $g$ and $h$ are respectively element of $G_n$ and $H_m$ then $g \star h$ is in $(G \star H)_{n+m+1}$.
}

\blockn{We will also introduce the joint operation on $\infty$-category by a black-box theorem similar to the one for the Gray tensor product:}

\block{\Th{(Ara-Maltsiniotis, \cite{ara2016joint})There is a monoidal structure $\star$ on the category of $\infty$-categories such that:

\begin{itemize}
\item The empty category is the unit object.

\item The linearization functor from $\infty$-categories to augmented $\infty$-groups is monoidal for the monoidal structure $\star$ on both categories.

\item For each object $X$, the functor:

\[ X \star \_  : \infty - Cat \rightarrow (\infty - Cat) _{X \backslash}\]

has a right adjoint.

\item There is a cell $t_{n,m}$ in $D_n \star D_m$, and $D_n \star D_m$ is obtained from $ D_n \star \partial D_m \coprod_{\partial D_n \star \partial D_m} \partial D_n \star D_m$ by freely adding the cell $t_{n,m}$ with specified boundary. Under the identification of $\Z(D_n \star D_m)$ with $\Z D_n \star \Z D_m$ the element $\delta(t_{n,m}) = \delta t_n \star \delta t_m$.

\end{itemize}

}

Note that as $ X \star \_ $ send the empty $\infty$-category to $X$ it cannot have a right adjoint itself, and the functor mentioned in the theorem does make sense.
\Dem{The first three points are proved in \cite{ara2016joint} section 7. The last property follows very directly from their construction of the Join in terms of Steiner's theory as $D_n \star D_m$ and all its sub-objects are described explicitly by Steiner complexes.}
}

\block{Following the exact same procedure as in section \ref{subsec_Gray_tensor_product}, one can show that given two polygraphs $X$ and $Y$, $X \star Y$ is a polygraph with generating the cells of $X$ and the cells of $Y$ (which gives sub-polygraphs isomorphic to $X$ and $Y$) and for each cell $x \in X$ of dimension $n$, and $y \in Y$ of dimension $m$, $x \star y$ is a cell of $X \star Y$ of dimension $n+m+1$. Moreover the linearization maps send the cells of $X \star Y$ to the generators of $\Z X \star \Z Y$ with the same name which allows to get some information on the source and targets of these cell $x \star y \in X \star Y$.  }

\block{\label{discus:descriptionOfCQ}In the present paper we will only be interested in the functor $X \mapsto * \star X$, where $*$ the unit $\infty$-category. We will denote it by $X \mapsto C X$, as homotopically speaking, it corresponds to a cone construction, and we will give a more explicit description only in this case. On an augmented strict $\infty$-group $e:G \rightarrow \Z$, $C G = (\Z *) \star G $ can be described as follow (we leave out the details of this computation):

\begin{itemize}

\item As an abelian group $CG$ identifies with $G \oplus G \oplus (\mathbb{Z} *)$. The first component corresponds to the natural morphisms $G \rightarrow CG$, the last to the morphisms $(\mathbb{Z} *) \rightarrow CG$. 

\item The function $\pi_i^{\epsilon}$ are given on the first and last component by how they act in $G$ and $\Z *$. We identify the corresponding element of $CG$ with the elements of $G$ and $\Z *$.

\item The elements of the middle component are denoted $T g$, the function $\pi_n^{\epsilon}$ are defined as:

\[ \begin{array}{r c l} \pi_n^{-}(T g) & =&  \left\lbrace \begin{array}{l l} T(\pi_{i-1}^+ g) & \text{ If $i>0$}\\ e(g) * & \text{ If $i=0$.}  \end{array} \right. \\
 \pi_n^{+}(T g) & = & \left\lbrace \begin{array}{l l} T(\pi_{i-1}^- g) + \pi_i^+ g - \pi_{i-1}^- g  & \text{ If $i>0$}\\ \pi_0^+ g & \text{ If $i=0$.}  \end{array} \right. \end{array} \]

\end{itemize}

}

\block{If $X$ is a polygraph then $CX$ is a polygraph which contains $X \subset CX$, as an additional $0$-cell $* \in CX$, and for each $n$-cell $x \in X$ it has one additional $(n+1)$-cell $T(x)$. We do not have a nice description of what are the source and target of $T(x)$ in general, but we know that the linearization of $CX$ is $C \Z X$ so we can compute a few things about these $T(x)$ using their value in the linearization. One can also describe them in low dimension:

\begin{itemize}

\item if $x \in X$ is a $0$-cell, then $T(x): * \rightarrow x$.

\item If $f:x \rightarrow y$ is a $1$-cell then $T(f)$ is a $2$-cell:

\[\begin{tikzcd}[ampersand replacement=\&]
* \arrow{d}[swap]{Tx} \arrow{dr}{Ty}[swap]{}[below,name=start]{}  \&  \\
x \arrow{r}[swap]{f} \arrow[Rightarrow,from=start] \& y
\end{tikzcd}\]

indeed one has in the linearization $\pi^-_1 \delta_{T f} = \delta_{Ty}$ and $\pi_1^+ \delta_{T f} = \delta_{Tx}+ \delta_f - \delta_x$, by the formula given for $C \Z X$.

\item And more generally, $T x$ is a cell completing the interior of a ``cone'' formed by $T$ applied to the boundary of $x$, and whose base is the cell $x$. But the precise expression highly depends on the shape of $x$ and we do not know any nice way to describe it in general.

\end{itemize}
}

\block{\label{Prop_Cplex=plex}\Prop{If $\underline{p}$ is a regular plex with universal arrow $p$, then $C \underline{p}$ is a regular plex with universal cell $Tp$.}

\Dem{Let $n$ be the dimension of $\underline{p}$, and consider $t:D_{n+1} \rightarrow C \underline{p} $ corresponding to the arrow $T(p)$. The first step of the proof is to show that:

\[ t (\sigma_{D_{n+1}}) = \sigma_{C \underline{p}} \]

And that $t$ preserves alternate positivity. 
Indeed (not writing all the $\delta$) :

\[ t(\sigma_{D_{n+1}})  = (-1)^{n+1} (T p) + \left( \sum_{i=0}^n (-1)^i \pi_i^{+} Tp \right) + \left( \sum_{i=0}^n (-1)^i \pi_i^{-} Tp \right)  \]

As $\underline{p}$ is a plex, one can write that:

\[ \sigma_{\underline{p}} =  (-1)^np + \sum_{i=0}^{n-1} (-1)^i (\pi_i^- p + \pi_i^+ p) \]

Hence, with some reordering, the second term of the sum becomes:

\[\begin{array}{r c l} \sum_{i=0}^n (-1)^i \pi_i^{+} Tp 
&=&
\pi_0^+ p +  \sum_{i=1}^n (-1)^i \left( T(\pi^-_{i-1} p) + \pi_i^+ p - \pi_{i-1}^- g \right)
\\ &=&
\sigma_{\underline{p}} + \sum_{i=0}^{n-1} (-1)^{i+1} T(\pi^-_i p) 
\end{array}  \]

Similarly, the third term becomes: 

\[\begin{array}{r c l} \sum_{i=0}^n (-1)^i \pi_i^{-} Tp & = & * + \sum_{i=1}^n (-1)^i T (\pi_{i-1}^+ p)\\
&=& * + \sum_{i=0}^{n-1} (-1)^{i+1} T(\pi_{i}^+ p)
\end{array}\]

Regrouping everythings, one gets:

\[ \begin{array}{r c l} \chi (\sigma_{D_{n+1}}) &=& (-1)^{n+1} T p + \sigma_{\underline{p}} + * - \sum_{i=0}^{n-1} (-1)^{i+1} (T\pi_{i}^+ p + T \pi_{i}^- p)\\
&=& * + \sigma_{\underline{p}} - T \sigma_{\underline{p}} \end{array}
 \]
which is indeed $\sigma_{C \underline{p}}$ as claimed: the cells of of $C \underline{p}$ are exactly the cells of $\underline{p}$, the cell $*$ and the cells $Tx$ for $x$ cells of $\underline{p}$, but with a shift of dimension, hence the minus signs. And all the terms appeared directly with the correct sign, hence showing that $\chi$ preserve alternate-positivity

One cannot directly conclude that $C \underline{p}$ is a plex using proposition \ref{Prop_delta_sigma_Charac_polyplex} and \ref{Prop_Sigmacriterion_boundary} as we do not now yet whether $C \underline{p}$ is a regular polygraphs. So instead we proceed by induction following the same scheme as we did for proposition \ref{Prop_tensor_product_of_regular}:
 
If we now that for $\underline{p}$ an $n-1$-dimensional regular plex, $C \underline{p}$ is a $n$-dimensional regular plex, then this implies that for each regular polygraph $P$ in dimension $n-1$, $C P$ is regular as well: indeed any plex appearing in $CP$ is either a cell of $P$ or the cell $*$ (in which case it is regular) or a cell $T x$ for $x \in P$, whose corresponding plex if $C \underline{x}$ , which is known to be regular.

Now if $\underline{p}$ is an $n$-dimensional regular plex, then this arguments also shows that $(C \underline{p}) - \{ T p \}$ is a regular polygraphs, and the computation above show that the restriction of $\chi$ as a map $\partial D_{n+1} \rightarrow (C \underline{p}) - \{ T p \}$ preserve the $\sigma$ and alternate-positivity. Hence proposition \ref{Prop_CaracGenericRegular} show that it is a generic and regular map, and hence when we glue the final cell $T p$ at its top one get a regular plex. 

}

}

\blockn{This has a very interesting corollary, which sadly we will not be able to use, but we thought it was interesting pointing it out:}

\block{\label{Cor_RegPg_is_WeakTest}\Cor{The category of regular plexes is a weak test category}

\Dem{ The construction $C$ above is a ``separating Décalage'' on the category of regular plexes in the sense of \cite{cisinski2011categorie}, and the existence of such a décalage proves that this category is a weak test category.}

Note that it is a very reasonable conjecture that $C$ generally send plex to plex without the regularity assumption. If this in indeed the case then the category of positive plex would as well be a weak test category contrary to what one could have guessed. But the point is that the model structure that we will construct below on the category of positive polygraphs does not have the same weak equivalence as those coming from the fact that this category is (conjecturally) a weak test category, while for the category of regular polygraph this will be the case. We believe this observation explain why the fact that regular plexes form a weak test category is not more useful here.

}

\blockn{Another corollary that we will use a lot more is:}

\block{\label{Prop_street_orientals}\Prop{There is downward closed\footnote{we mean that if $v$ is a plex and $v \rightarrow \Ocal(n)$ is any polygraphic morphisms, then $v$ is also an orientals.} full subcategory of the category of regular plex isomorphic to the semi-simplicial category, whose objects are given by Street's Orientals $\Ocal(n)$. }

By the semi-simplicial category, we mean the category of finite non-empty totally ordered set and order preserving injection between them. 

Street's Orientals have been first defined by Street in \cite{street1987algebra}. An other (equivalent) definition has been given by R.Steiner in \cite{steiner2004omega} using his theory of chain complexes. Moreover, in section 8 of \cite{ara2016joint}, D.Ara and G.Maltsiniotis have shown that the orientals can be alternatively defined in terms of the join as:

\[\begin{array}{r c l} 
\Ocal(0) & = &  *\\
\Ocal(n) & = & C \Ocal(n-1)\end{array} \]

\Dem{It is immediate from the inductive definition of Ara and Maltsiniotis and our proposition \ref{Prop_Cplex=plex} that Orientals are regular plexes. It is also rather clear, both from Street's initial definition and Steiner's definition that any cell of dimension $n$ that appears in $\Ocal(m)$ is itself of shape $\Ocal(n)$, which shows that Orientals are downward closed in the category of plexes. Finally, Steiner description show that the cells of dimension $k$ of the Orientals $\Ocal(n)$ are given by subsets with $k+1$-elements of  $\{0,\dots,n\}$, and given a cell $x$ of dimension $k$ of $\Ocal(n)$, corresponding to a subset $X \subset \{0,\dots,n\}$ of cardinal $k+1$ , the map from $\Ocal(k)$ to $\Ocal(n)$ corresponding to $x$ send a cell corresponding to subset of $\{0,\dots,k\}$ to the corresponding subset of $X$ (identified with $\{0,\dots,k\}$ by the only order preserving bijection). Which is just a convoluted way of saying that polygraphic morphisms between the $\Ocal(n)$ corresponds (functorially) to order preserving injection  between the $\{0,\dots,n\}$. I.e. that they form a category equivalent to the semi-simplicial category.

}
}

\section{Homotopy theoretic aspect}

\blockn{In this section we take $\C$ to be a class of polygraphs which is either the class of regular polygraphs, or the class of positive polygraphs. In particular, it is a $\otimes$-stable algebraic class of polygraphs as in section \ref{section_algebraic_class}. I.e. one has a monoidal bi-closed category of $\C$-polygraphs, a monoidal bi-closed category of $\C$-categories and a ``free $\C$-category'' functor, which we will denote $X^*$, and which is a monoidal functor from $\C$-polygraphs to $\C$-categories. We also have a ``Nerve'' functor (or forgetful functor) from $\C$-categories to $\C$-polygraphs, which is right adjoint to $(\_)^*$.

Note that in particular the $X^*$ of this section, is different from the one we have used until now: Indeed before $X^*$ denoted the free $\infty$-category on a polygraphs (so the case $\C=$ all polygraphs), but now we use it not denotes the free $\C$-category. If $\C = \Pb^{reg}$ is the class of regular polygraph, this means than $X^*$ is now only the regular $\infty$-category of regular arrows of ``$X^*$'', instead of the ordinary $\infty$-category of all arrows of ``$X^*$''.}

\subsection{A groupoid folk weak model structure on $\C$-polygraphs and $\C$-categories}
\label{subsec:FolkMS}

\blockn{We will construct weak model structures (in the sense of \cite{henry2018weakmodel}) on the category of $\C$-polygraphs and of $\C$-categories using the Gray tensor product on both these categories and the theorem $3.2.1$ of \cite{henry2018weakmodel}. So essentially one only needs to give generating cofibrations, generating anodyne maps and check that the pushout-product condition with respect to the Gray tensor product are satisfied, and that one can define a weak cylinder object for the unit of the Gray tensor product.

The model structure that we will construct is an analogue of the folk model structure, but in a version intended to model $\infty$-groupoids and not $\infty$-categories (i.e. a left Bousefield localization of the folk model structure for which the local fibrant objects are those in which every arrow is weakly invertible).
}

\blockn{We start with the category of $\C$-polygraphs.}

\block{\begin{itemize}

\item The cofibrations will be the monomorphisms. A natural choice of generating cofibrations are all the map $\partial \underline{c} \hookrightarrow \underline{c}$ where $\underline{c}$ is a $n$-plex in $\C$ and $\partial \underline{c}$ is the sub-polygraph containing all cells except the unique $n$-dimensional cell $c$.

\item The generating anodyne maps will be all the maps:

 \[ \Lambda^a \underline{c} \hookrightarrow \underline{c} \]

such that $\underline{c}$ is a $n$-plex, $a$ is an $(n-1)$-cell of $\underline{c}$ and $ \Lambda^a \underline{c}$ is a sub-polygraph of $\underline{c}$ which contains all the cells except $a$ and the top dimensional cell of $\underline{c}$.

\end{itemize}

This of course should remind the reader of the generating cofibrations and generating trivial cofibrations of the simplicial model structure which have exactly this form. In fact simplicial horn inclusion are exactly our generating anodyne maps in the case where the target is an orientals.

}

\blockn{As for other examples treated in \cite{henry2018weakmodel}, the $\Lambda^a \underline{c} \hookrightarrow \underline{c}$ will not be ``Generating acyclic cofibrations'' in the sense that not every acyclic cofibration will be a retract of a transfinite composition of pushout of these generating map. The left class of the weak factorization system generated by  $J = \{ \Lambda^a \underline{c} \hookrightarrow \underline{c} \}$ will be called anodyne map. One will define fibrant object as the object having the lifting property against $J$ and fibration between fibrant objects will be exactly the map having the right lifting property against $J$. Anodyne maps will all be acyclic cofibrations, but the converse need not be true. In fact, our first attempt at producing a cylinder in \ref{discuss_contruc_cylinder} is very likely to be an example of an acyclic cofibration which is not anodyne. Such set $J$ are also sometime called ``Pseudo-generating set of acyclic cofibration''.
}

\blockn{We start by two easy lemmas to get a better understanding of the pushout of our generating anodyne maps.}

\block{\label{Lem_occurenceOfCellsInPlexAndPolyplex}\Lem{\begin{itemize}

\item Let $\underline{p}$ be an $n$-polyplex. Then each $n$-cell of $\underline{p}$ appears exactly once in the universal cell $p$.

\item Let $\underline{p}$ be an $n+1$-plex. Each $n$-cell of $\underline{p}$ appears exactly once in either $\pi_n^{+} \underline{p}$ or in $\pi_n^{-} \underline{p} $ and not in the other.

\end{itemize}

}

We remind the reader that the number of time an $n$-cell appears in an $n$-arrow of a polygraph is well defined, for example because of \ref{Prop_delta_count_in_top_Dim}.

\Dem{The first claim is shown by induction on polyplexes: it is clear if $\underline{p}$ is a plex, if $\underline{p}$ is a composite $\underline{p} = \underline{q} \#_k \underline{h}$, one can be assume that $k<n$ otherwise the composition is trivial. Hence the $n$-cells of $\underline{p}$ are exactly the disjoint union of the $n$-cell of $\underline{q}$ and of $\underline{h}$. Each $n$-cell of $\underline{q}$ appears exactly once in the universal arrow of $\underline {q}$ and each $n$-cell of $\underline{h}$ appears exactly once in the universal arrow of $\underline{h}$, hence they all appears exactly once in their composite, which is the universal arrow of $\underline{p}$.

For the second claim, any $n+1$-plex $\underline{p}$ is obtained by gluing together the polygraph $\pi_n^+ \underline{p}$ and $\pi_{n}^- \underline{p}$ along their common boundary and adding a single $(n+1)$-dimensional cell, as this common boundary is $(n-1)$-dimensional, this gluing does not affect the $n$-cells, hence the set of $n$-cells of $\underline{p}$ is the disjoint union of the set of $n$-cell of the polyplexes $\pi_n^+ \underline{p}$ and $\pi_n^- \underline{p}$. The first point of the lemma then concludes the proof.

}
}

\block{\label{Lem_Anodyne_maps_and_Nb_of_Cells}\Lem{Let $u \hookrightarrow v$ be a polygraphic extension in $\C$, such that $v$ has two cells $x$ and $\theta$ in addition of those of $u$, such that:

\begin{itemize}

\item $x$ is of dimension $n$ and $\theta$ is of dimension $n+1$

\item $x$ appears exactly once in either the source or the target of $\theta$, and not in the other.

\end{itemize}

Then $u \hookrightarrow v$ is a pushout of one of the generating anodyne maps $\Lambda^x \underline{\theta} \hookrightarrow \underline{\theta}$.

}

\Dem{Take $\chi_{\theta}:\underline{\theta} \rightarrow v$ be the plex representing $\theta$. We claim that there is a unique $n$-cell $x_0 \in \underline{\theta}$ (called $x$ in the statement of the lemma) whose image in $v$ is $x$: indeed, because of lemma \ref{Lem_occurenceOfCellsInPlexAndPolyplex}, each $n$-cell of $\underline{\theta}$ appears exactly once in either of $\pi_n^{\epsilon} \theta$. Hence if they were no cell of $\underline{\theta}$ sent to $x$, then the $\pi_n^{\epsilon} \theta$ would have expressions not involving $x$, and if they were two $n$-cells of $\underline{\theta}$ sent to $x$, then $x$ will have to appears (non trivially) more than once in an expressions of $\pi_n^{\epsilon}(\theta)$.

In particular the map $\Lambda^{x_0} \underline{\theta} \rightarrow \underline{\theta} \rightarrow v$ does not have $x$ or $\theta$ in its image, hence factors into $u$, and the pushout $u \coprod_{\Lambda^{x_0} \underline{\theta}} \underline{\theta}$ is exactly adding the cells $x$ and $\theta$ to $u$ so is isomorphic to $v$.
}
}

\block{\label{discuss_contruc_cylinder}We now construct the weak cylinder object for the unit of the tensor product. The unit for the Gray tensor product is the polygraph $*$ with only one cell. One needs to construct a weak cylinder object for $*$. Obviously, one take $D_1 := * \rightarrow *$ as $I *$, with the obvious map $* \coprod * \hookrightarrow I *$. The map $* \hookrightarrow I *$ is one of our generating anodyne map. There is no map $D_1 \rightarrow *$, so instead one needs to construct an object $D *$ detecting trivial path, i.e. such that one has a map:

 \[ (I *)\coprod_{* \coprod *} * \rightarrow D * \]

with the map $* \hookrightarrow D*$ being an acyclic cofibrations. The first natural candidate for $D * $ is the following polygraph:

\begin{itemize}
\item It has one $0$-cell $*$.
\item It has one $1$-cell $p: * \rightarrow *$.
\item It has one $2$-cell $\theta: p \#_0 p \rightarrow p$.
\end{itemize}

The idea being that as our model structure makes every cell invertible, or more precisely ``divisible'', one can assert that a given cell (here $p: * \rightarrow *$) is an identity by just saying that it is idempotent (an invertible idempotent always being an identity).  It is not totally trivial that $* \rightarrow D * $ is an acyclic cofibration, because it is not an anodyne map, so instead we will use the following polygraph $D'*$:

\begin{itemize}
\item It has two $0$-cell $*$ and $t$.
\item It has two $1$-cell $w:t \rightarrow *$ and $p:*\rightarrow *$.
\item It has one $2$-cell $\theta : t \#_0 p \rightarrow p$.
\end{itemize}

The cofibration  $(I *)\coprod_{* \coprod *} * \rightarrow D' *$ sending the unique $1$-cell to $p$. $* \rightarrow D' *$ is anodyne: it is obtain by two pushout of generating anodyne map, one adding $t$ and $w$, and one adding $p$ and $\theta$, hence it qualifies to complete our weak cylinder object.
}

\blockn{The pushout product condition for generating cofibrations between polygraphs have been checked during the proof of \ref{prop:TensProdOfPG}. We only need to check the pushout product condition for a generating anodyne map with a generating cofibration.}

\block{\label{Lem:pushoutCondForPG}\Lem{Let $i$ be a generating cofibration and $j$ be a generating anodyne maps. Then $i \otimes' j$ and $j \otimes' i$ are pushout of generating anodyne maps.}

\Dem{We will prove it for $i \otimes j'$, the other one being completely similar (even a little simpler).

Let $i :\partial \underline{c} \hookrightarrow \underline{c}$ and $j' : \Lambda^a \underline{\theta} \hookrightarrow \underline{\theta}$. As usual, we denote by $\theta$ and $c$ the maximal cell of $\underline{\theta}$ and $\underline{c}$.

$i \otimes j'$ is the polygraphic map $P \hookrightarrow \underline{c} \otimes \underline{v}$, where $P$ is the pushout:

\[ P := \partial \underline{c} \otimes \underline{\theta} \coprod_{\partial \underline{c} \otimes \Lambda^a \underline{\theta} } \underline{c} \otimes \Lambda^a \underline{\theta} \]

but we know very well what are the cells of all these polygraph, so $P$ is the subpolygraph of $\underline{c} \otimes \underline{v}$ of cells $x \otimes y$ such that either $x \in \partial \underline{c}$ or $y \in \Lambda^a \underline{\theta}$. In particular there is only two cells that are not in $P$: $c \otimes a$ and $c \otimes \theta$. If $i$ is of dimension $n$ and $a$ of dimension $m$, they are respectively of dimension $n+m$ and  $n+m+1$. In order to conclude, using lemma \ref{Lem_Anodyne_maps_and_Nb_of_Cells}, one just need to shows that $c \otimes a$ appears exactly once in either the source or the target of $c \otimes \theta$. But the number of occurrences of cells in tensor products is essentially the only things that one know how to compute with tensor product, and so this is easy to check (for example using \ref{Lem_pi_tens_degminus1}) :

\[ \pi^{\epsilon}_{n+m} (\delta_c \otimes \delta_{\theta} ) = \pi^{\epsilon}_{m-1} \delta_c \otimes \delta_{\theta} + c \otimes \pi^{|m|\epsilon}_n \delta_{\theta} + \dots \]

where the $+ \dots$ is a terms of dimension lower than that of $c \otimes a$. As $c \otimes a$ has dimension $n+m$ its number of occurrences in  $\pi^{\epsilon}_{n+m} (c \otimes \theta ) $ is the same as the coefficient of $\delta_c \otimes \delta_a$ in the expression above.  $\delta_c \otimes \delta_a$ cannot appear in $ \pi^{\epsilon}_{m-1} \delta_c \otimes \delta_{\theta}$, so it only appear in $c \otimes \pi^{|m|\epsilon}_n \delta_{\theta}$, and there it appears exactly if $a$ appears in $\pi_n^{|m|\epsilon} (\theta)$, i.e. it appears exactly once, for exactly one value of $\epsilon$, so this concludes the proof.
}
}

\blockn{At this point we have:}

\block{\label{Th:PolygraphModelStructure}\Th{There is a weak model structure on the category of $\C$-polygraph such that:

\begin{itemize}
\item Every object is cofibrant and cofibrations are the monomorphisms\footnote{Constructively, the cofibrations are the level wise complemented monomorphisms.}.

\item Fibrant object are those with the right lifting property against all the generating anodyne maps, and fibrations between fibrant objects are characterized by the lifting property against generating anodyne maps.

\end{itemize}

}

This weak model category can is in fact a right semi-model structure in the sense of \cite{barwick2010left} (we need to define more precisely the notion of fibration between non-fibrant objects in order to make this true). But this observation is not of any use for the present paper, and we will not prove it.

\Dem{This follows immediately from theorem $3.2.1$ of \cite{henry2018weakmodel}. Our generating cofibrations and anodyne maps are map between finitely presentable objects in a $\omega$-presentable category so one immediately have a weak factorization system by the small object argument. One knows that it is a monoidal categories because of \ref{prop:TensProdOfPG}, and closed because the tensor product commutes to colimits in each variables and the category is presentable. We have a weak cylinder objects for the unit (~\ref{discuss_contruc_cylinder}) and we checked the pushout-product conditions (~\ref{Lem:pushoutCondForPG}). This shows that there is a weak model structure.}
}

\blockn{And one also easily get:}

\block{\Th{There is a weak model structure on the category of $\C$-category such that:

\begin{itemize}

\item Cofibrant objects are the polygraphs, cofibrations between cofibrant objects are the polygraphic extension\footnote{and their retract if one wants cofibrations to be stable under retract. It has been shown by F.Metayer in \cite{metayer2008cofibrant} that retract of polygraphs are again polygraphs.}.

\item Fibrant objects are the $\C$-categories with the right lifting property against the polygraphic generating anodyne maps. Fibrations between fibrant objects are also characterized by this lifting property.

\end{itemize}

}

One can take the set of $\partial D_n \rightarrow D_n$, or alternatively, all the $\partial c \rightarrow c$ for $c$ a plex as generating cofibrations.

\Dem{One take as generating cofibrations and generating anodyne maps the image of these in the category of $\C$-polygraphs. As the functor from $\C$-polygraphs to $\C$-categories is monoidal and commutes to colimits the pushout-product condition immediately follows from those in the category of $\C$-polygraphs. Also the weak cylinder object for $*$ in the category of $\C$-polygraphs gives a weak cylinder object for $*$ in the category of $\C$-categories. So the theorem follows also from theorem $3.2.1$ of \cite{henry2018weakmodel}.
}
}

\block{\Def{A Simpson $\C$-groupoids is a fibrant object in this weak model category of $\C$-category.}

This does corresponds to the idea of category with strict composition and associative and weak units and inverse: the underlying algebraic structure has all the composition defined by $\C$-polyplexes and the are strictly compatible (a $\C$-composition of $\C$-composition is given by a single $\C$-composition whose diagram is just the composition of diagrams), and fibrant objects additionally have the property that every cell is ``divisible'' in a relatively strong sense (that is what the lifting property against all the generating anodyne map means) which implies both the existence of weak inverse and weak units (weak units are characterized as the idempotent).

The homotopy category of the category of $\C$-categories (in the sense of \cite{henry2018weakmodel}), is hence equivalent to the localization of the category of Simpson groupoid at all trivial fibrations, which are the map having the lifting property against all the $\partial D_n \rightarrow D_n$, so this is a very reasonable definition of the homotopy category of $\infty$-groupoid with strict composition and weak units and inverse. Depending on if $\C= \Pb^+$ or $\C= \Pb^{reg}$ one has more or less compositions operations that are made strict.

Using the $\pi$-sets characterization of equivalences in weak model categories presented in section $2.5$ of \cite{henry2018weakmodel} one can also gets a characterization of weak equivalences between Simpson $\C$-groupoids as the map inducing bijection on the $\pi_n$ for an appropriate definition of $\pi$-sets that is ``left as an exercise'', but taking the generating cofibration to be the $\partial D_n \hookrightarrow D_n$ one does get a single $\pi$-set for each $n$ that will be the exact analogue of the usual homotopy groups.

Our definition also require a rather high number of ``divisibility conditions'' on our Simpson groupoids. It would be interesting to try to find smaller set of conditions that are sufficient to implies that a given $\C$-category is a Simpson groupoid. For example something like ``having units and weak $n$-inverse of $n$-arrows''. We will not discuss this question in the present work.
}

\subsection{The Quillen equivalence between $\C$-polygraphs and $\C$-categories}
\label{subsec:QuillenEquivCatPolyg}

\blockn{In this subsection we show that the free $\C$-category functor and the Nerve functor:

\[ (\_)^* : \C - \text{Polygraphs} \leftrightarrows \C-\text{Categories} : N \]

induce a Quillen equivalence between the two weak model structure constructed in the previous subsection. This is actually a Quillen equivalence in a very strong sense: its unit and co-unit of adjunction are respectively acyclic cofibrations and acyclic fibrations (as soon as the domain and co-domain are correctly fibrant and cofibrant).
}

\block{\label{Lem_IsQuillenAdjunction}\Lem{The adjunction $(\_)^* \dashv N$ is a Quillen adjunction.}

\Dem{This is immediate: Both the generating cofibrations and generating anodyne map between $\C$-categories were taken to be the image under $(\_)^*$ of those between $\C$-polygraphs.}
}

\block{\label{Lem_counit_trivFib}\Lem{For any $\C$-category $X$, the counit map $N(X)^* \rightarrow X$ is an $I$-fibration.}

We mean that it has the right lifting property against all the cofibrations. I.e. if $X$ is assumed to be fibrant, then this is an acyclic fibration.

\Dem{Consider a lifting problem of the form:

\[\begin{tikzcd}[ampersand replacement=\&]
A^* \arrow[hook]{d}{i} \arrow{r}{u} \& N(X)^* \arrow{d}{\mu_X} \\
B^* \arrow[dotted]{ur}{\exists ?} \arrow{r}{d} \& X
\end{tikzcd}\]

With $i$ a polygraphic extension, one can assume $A$ and $B$ finite.

Our first observation, is that if $u$ is a polygraphic map then one can check that if $\tilde{d} : B \rightarrow N(X)$ is the map corresponding to $d$ under the adjunction  $(\_)^* \dashv N$, then $\tilde{d}^* : B^* \rightarrow N(X)^*$ is a diagonal filling. In fact it is the unique polygraphic diagonal filling.

If the map $u$ is not polygraphic, one can factor it into $A^* \rightarrow A'^* \rightarrow N(X)^*$ such that the first map is arbitrary and the second map is polygraphic (either using the generic factorization, or more generally the fact that $A^*$ is finite and $N(X)^*$ is a directed colimits of its finite sub-polygraphs), one can then form a diagram (with $B'$ a pushout):

\[\begin{tikzcd}[ampersand replacement=\&]
A^* \arrow{r} \arrow[hook]{d} \&  A'^* \arrow[hook]{d} \arrow{r}{u^*} \& N(X)^* \arrow{d}{\mu_X} \\
B^* \arrow{r} \& B'^* \arrow{r} \arrow[dotted]{ur} \& X
\end{tikzcd}\]

$A' \rightarrow B'$ is a polygraphic extension because it is a pushout of a polygraphic extension, hence the dotted diagonal filling exists due to the first half of the proof.
}
}

\block{\label{Lem_Xstar_fibrant}\Lem{If $X$ is a fibrant $\C$-polygraph, then $X^*$ is a fibrant $\C$-category.}

\Dem{Consider a lifting problem:

\[\begin{tikzcd}[ampersand replacement=\&]
A^* \arrow[hook]{d}{j^*} \arrow{r}{u} \& X^* \\
B^* \arrow[dotted]{ur}{\exists ?}\& \\ 
\end{tikzcd}\]

With $j$ a generating anodyne map. Then one can factor the map $u$ into a polygraphic map: $A^* \rightarrow A'^* \overset{v^*}{\rightarrow} X$. And then one can form the pushout:

\[\begin{tikzcd}[ampersand replacement=\&]
A^* \arrow[hook]{d}{j^*} \arrow{r} \& A'^* \arrow[hook]{d}{j'^*} \arrow{r}{v^*} \& X^* \\
B^*  \arrow{r} \& B'^* \arrow[dotted]{ur} \& \\ 
\end{tikzcd}\]

Lemma \ref{Lem_Anodyne_maps_and_Nb_of_Cells} shows that $j'$ is again a pushout of a generating anodyne map (its adding two cells with the appropriate number of appearances in one into the source and target of the other, because $j$ is) hence there exists a polygraphic dotted arrow because $X$ is fibrant as a polygraph, which produce in particular filling for the diagram of $\C$-categories.
}
}

\block{\label{Prop_unit_triv}\Prop{If $X$ is a fibrant polygraph, then the unit of adjunction $\eta_X : X\rightarrow N(X^*)$ is an acyclic cofibration.}

The hypothesis that $X$ is fibrant cannot be weakened: if $X = \bullet \rightarrow \bullet \rightarrow \bullet$, then $X^*$ has only one additional regular/positive arrow which is the composite, so that $N(X^*)$ is a triangle without $2$-cell and the map $X \rightarrow N(X^*)$ is not acyclic.

\Dem{The map $X \rightarrow N(X^*)$ is a monomorphisms (because the functor $(\_)^*$ is faithful) hence a cofibration\footnote{Constructively, it is also level wise complemented as one can always decide if an arrow of a polygraph is a generator or not.}. By lemma \ref{Lem_Xstar_fibrant}, $X^*$ is fibrant and hence as $N$ is a right Quillen functor (lemma \ref{Lem_IsQuillenAdjunction}), $N(X^*)$ is fibrant. Hence $\eta_X$ is a cofibration between fibrant objects. In order to show that it is an equivalence, we will check that it has the weak right lifting property against the generating cofibrations:

\[\begin{tikzcd}[ampersand replacement=\&]
\partial \underline{c} \arrow[hook]{d} \arrow{r}{u} \& X \arrow{d}{\eta_X} \\
\underline{c} \arrow{r}{d} \& N(X^*) \\
\end{tikzcd}\]

$u$ corresponds to a pair of parallel arrows ($s$,$t$) in $X^*$ of shape given by the source and target of $\underline{c}$. The map $\underline{c} \rightarrow N(X^*)$ is the same a $\C$-morphism $c \rightarrow X^*$ so corresponds to an arrow $e:s \rightarrow t$ in $X^*$

As $X$ is fibrant, one can find cells $e':s \rightarrow t$ and $\theta:e \rightarrow e'$: indeed this corresponds to the lifting property of $X$ against a generating anodyne cofibration which add two cells $\widetilde{e'}$ and $\widetilde{\theta}$ with $\widetilde{e'}$ being the target of $\theta$ and not appearing in its source.

This can be packaged into a diagram:

\[\begin{tikzcd}[ampersand replacement=\&]
\partial \underline{c}
\arrow[hook]{dr}
\arrow[hook]{ddd} \arrow{rr}{u} \&  \& X \arrow{ddd}{\eta_X} \\
\& \underline{c} \arrow[hook]{d} \arrow{ru}{u'} \& \\
\& \underline{w} \arrow{dr}{\theta} \& \\
\underline{c}
\arrow[hook]{ur}
\arrow{rr}{d} \& \& N(X^*) \\
\end{tikzcd}\]

Where $u'$ is the arrow that extend $u$ by sending the universal cell to $e'$. $w$ is the polygraphic extension of $\partial \underline{c}$ by two cell $c_1$ $c_2$ corresponding to the universal cell of $c$, and one additional $\widetilde{\theta} : c_1 \rightarrow c_2$. And the map $\underline{w} \rightarrow N(X^*)$ is the map $\underline{w}^* \rightarrow X^*$ which send $\widetilde{\theta}$ to $\theta$ (so in this new presentation we have forgotten that $\theta$ was actually a cell instead of an arrow).

If one can prove that $\underline{w}$ is a weak cylinder object this conclude the proof. One clearly have that the natural map $\underline{c} \coprod_{\partial \underline{c}} \underline{c} \rightarrow \underline{w}$ is a cofibration, $\underline{c} \hookrightarrow \underline{w}$ is one of the generating anodyne map (it is adding $c_2$ and $\theta : c_1 \rightarrow c_2$). To conclude we need to construct an object $D_{\partial \underline{c}} \underline{c}$ detecting trivial homotopies, but this can be done exactly as for the construction of cylinder object of $*$ in \ref{discuss_contruc_cylinder}.
}

}

\block{\Th{The adjunction $(\_)^* \dashv N$ is a Quillen equivalence between the category of $\C$-polygraphs and the category of $\C$-categories.}
\Dem{We have proved in lemma \ref{Lem_IsQuillenAdjunction} that it is a Quillen adjunction. In lemma \ref{Lem_counit_trivFib} that if $X$ is a fibrant $\C$-category then the co-unit $\mu_X$ is an acyclic fibration and finally proposition \ref{Prop_unit_triv} show that the unit on a bifibrant object is an acyclic cofibration. This is sufficient to show that one has a Quillen equivalence: If $X$ is a cofibrant polygraph and $X \overset{\sim}{\hookrightarrow} Y$ is a fibrant replacement, then $Y^*$ is a fibrant replacement of $X^*$ because of lemma \ref{Lem_Xstar_fibrant}, and $X \overset{\sim}{\hookrightarrow} Y \overset{\sim}{\hookrightarrow} N(Y^*)$ is a weak equivalence. If $Z$ is a fibrant $\C$-category, then $N(Z)$ is cofibrant and $N(Z)^* \rightarrow Z$ is an acyclic fibration.

}
}

\subsection{The Quillen equivalence between  regular polygraphs and spaces}
\label{subsec:QuillenequivSpace}

\blockn{In this final section, we will prove that the category $\Pb^{reg}$ of regular polygraphs, with the weak model structure of theorem \ref{Th:PolygraphModelStructure} is Quillen equivalent to the category of spaces. While we will not directly use the fact that the category of regular plexes is a weak test category (corollary \ref{Cor_RegPg_is_WeakTest}), the functor that will induce an equivalence between $\Pb^{reg}$ is the usual functor that send an object $X \in \Pb^{reg}$, seen as a presheaf on the category of regular plex and maps it to its category of elements: $Plex^{reg}/ X$. So we will also show that this model structure we have on this category and its equivalence to the category of spaces is compatible to the Grothendieck style homotopy theory coming from the fact that the category of regular plex is a weak test category.

\bigskip

In order to get a Quillen functor instead of a mere homotopical functor, we will not use the homotopy theory of category as a target, but instead one will take the simplicial Nerve of the category $N(Plex/X)$ and call it the geometric realization of $X$, denoted by $R(X)$.

\bigskip

It appears that due to the directed nature of our category of plex, the simplicial set $R(X)$ has the following property:

\begin{itemize}

\item A cell $x$ of $R(X)$ is non-degenerate if and only if the corresponding functor $x:[n] \rightarrow Plex/X$ is strictly increasing (i.e. $x(i+1)$ is of greater dimension than $x(i)$).

\item The faces of a non-degenerate cells are always non-degenerate.

\item if $f :X \rightarrow Y$ is a morphisms in $\Pb^{reg}$ then the induce morphisms $R(f)$ send non-degenerate cells to non-degenerate cells.

\end{itemize}

This implies that $R$ actually factor into the category of semi-simplicial sets and the functor which freely add units. One will keep denoting $R(X)$ the semi-simplicial version of the geometric realization whose $n$-cells are only the ``non-degenerate'' functor $[n] \rightarrow Plex/X$ as above.

It will be convenient to work with semi-simplicial sets instead of simplicial sets because one has a left adjoint functor $\Ocal$ from semi-simplicial sets to $\Pb^{reg}$ induced by Street's Orientals.

As ordinary simplicial sets will actually never appear in the present paper, one will denote by $\Delta$ the semi-simplicial category, and $\widehat{\Delta}$ the category of semi-simplicial sets.

}

\block{The functor $R:\Pb^{reg}\rightarrow \widehat{\Delta}$ can be written in level $n$ as:

\[ R(X)_n = \coprod_{e_0 \rightarrow \dots \rightarrow e_n} X(e_n) \]

Where the coproduct is over all the strictly increasing chains $e_0 \rightarrow \dots \rightarrow e_n$ of plexes. In particular it is clear on that expression that $R(X)$ commutes to all colimits, and hence is a left adjoint functor. 

In order to go any further one will need the following key lemma. The failure of this lemma for more general positive (non regular) polygraphs is what makes the general form of the Simpson conjecture still out of reach.

}

\block{\label{Lemma_Key_weakequivalence}\Lem{Let $W$ be a class of monomorphisms in $\Pb^{reg}_n$ which contains isomorphisms, satisfies $2$-out-of-$3$ among monomorphisms, and the ``cube lemma'' stated below. Then the following are equivalents:

\begin{itemize}

\item[$(A_n)$] The generating anodyne map in $\Pb^{req}_n$ are in $W$.

\item[$(B_n)$] Any $f:\underline{p}\rightarrow \underline{q}$ between plexes in $\Pb^{reg}_n$ is in $W$.

\end{itemize}

And moreover they further imply:

\begin{itemize}

\item[$(C_n)$] Any $f: \underline{p} \rightarrow \underline{q}$ from a plex to a polyplex in $\Pb^{reg}_n$ is in $W$.

\item[$(D_n)$] For any regular $n$-polyplex $\underline{p}$ with spherical boundary and $x$ a $n$-cell of $\underline{n}$ the map $\partial \underline{p} \rightarrow \underline{p}-\{x\}$ is in $W$.

\item[$(E_n)$] For any $f:\underline{p} \rightarrow \Lambda^x \underline{\theta}$ with $\theta$ a plex in $\Pb^{reg}_{n+1}$ and $\underline{p}$ a plex in $\Pb^{reg}_n$, $f \in W$.

\end{itemize}

}

In our case, the ``cube lemma'' will be the following property:

If one has a diagram:

\[\begin{tikzcd}[ampersand replacement=\&]
B\arrow[hook]{d}{\in W} \& \arrow[hook]{l} A \arrow[hook]{r} \arrow[hook]{d}{\in W} \& C \arrow[hook]{d}{\in W} \\
B'\& \arrow[hook]{l} A' \arrow[hook]{r} \& C'
\end{tikzcd}\]

Where the vertical map are monomorphisms in $W$ and the horizontal map are monomorphisms then the natural comparison map:

\[ B \coprod_A C \hookrightarrow B' \coprod_{A'} C' \]

is also in $W$. This conditions implies in particular that map in $W$ are stable under pushout along monomorphisms by considering:

\[\begin{tikzcd}[ampersand replacement=\&]
C \arrow{d}{Id} \& \arrow[hook]{l}A \arrow{d}{Id} \arrow[hook]{r} \& A \arrow{d}{\in W} \\
C \& \arrow[hook]{l}A \arrow[hook]{r} \& B
\end{tikzcd}\]

Note that it is also trivial that  $C_n$ also implies $A_n$ and one will see that $E_n$ also implies $A_n$. But those does not seems to have any interest.

\Dem{We will proceed by induction: one assumes that $A_{n-1},B_{n-1},C_{n-1}, D_{n-1}$ and $E_{n-1}$ are all satisfied, and we will prove that the propositions $A_n,B_n$ are equivalent and also implies $C_n, D_n$ and $E_n$. 

\bigskip

$A_n \Rightarrow B_n$:

Let $\underline{p} \rightarrow \underline{q}$ be a morphism between plexes in $\Pb^{reg}_n$. If $\underline{p}$ has the dimension of $\underline{q}$, then this an isomorphisms hence is in $W$. If $\underline{p}$ is of dimension strictly less than $\underline{q}$, then it factors into one of the generating anodyne maps $\underline{p} \rightarrow \Lambda^x \underline{q} \hookrightarrow \underline{q}$. The first map is in $W$ because of $E_{n-1}$ and the second because of $A_n$. This implies $B_n$.

\bigskip

$B_n \Rightarrow C_n$:

We prove by induction on the regular polyplex $\underline{q}$, that for any plex $\underline{p}$, any map $\underline{p} \rightarrow \underline{q}$ is in $W$. If $\underline{q}$ is of dimension $<n$, then this is already known (by$C_{n-1}$), if $\underline{q}$ is an $n$-plex, then this is $B_n$. If $\underline{q} = \underline{h} \#_i \underline{k}$, with $\underline{h}$ and $\underline{p}$ satisfying the assumption, then $\pi^+_i \underline{h} \subset \underline{h}$ and $\pi_i^- \underline{k} \subset \underline{k}$ are monomorphisms and are in $W$ because of the induction hypothesis applied to any plex of the subpolygraph (and $W$ satisfies $2$-ou-of-$3$). hence as $ \underline{q} = \underline{h} \coprod_{\pi_i^- \underline{k}} \underline{k}$, this implies that both map $\underline{h} \rightarrow \underline{q}$ and $\underline{k} \rightarrow \underline{q}$ are in $W$ because $W$ is stable under pushout along monomorphisms, as the map $\underline{p} \rightarrow \underline{q}$ factor into of one of the two component (by a map which is in $W$ by induction) this proves $C_n$.

\bigskip

$C_n \Rightarrow D_n$:

Let $\underline{p}$ be an $n$-plex and $x$ a cell of $\underline{p}$. One can always consider a ``vertical decomposition'' of $\underline{p}$ in the form $\underline{p} = \underline{s} \#_{n-1} \widetilde{\underline{x}} \#_{n-1} \underline{t}$ where $\widetilde{\underline{x}}$ is a polyplex with a single $n$-cell $x$ (and cells of lower dimension) and $\underline{s}$ and $\underline{t}$ are arbitrary polyplexes. One apply the construction of lemma \ref{Lemma_Stuff} to $\widetilde{\underline{x}}$, one get a $(n-1)$-polyplex $\underline{x'}$ with a $(n-1)$-cell $x'$ such that $\underline{x'} - \{x' \}$ identify exactly to $\widetilde{\underline{x}}$ minus its unique $n$-cell $x$ as well as all the inner cell of $\pi^-_{n-1}x$ and $\pi^+_{n-1} x$, and moreover $\pi_{n-2}^{\epsilon} \widetilde{\underline{x}} = \pi_{n-2}^{\epsilon} \underline{x'}$.

One can then see that:

\[ \underline{p} - \{x \} = \underline{s} \coprod_{\underline{x'}-\{x'\}} \underline{t} \]

Indeed, $\underline{s}$ and $\underline{t}$ are sub-objects of $\underline{p}$, any cell of $\underline{p}$ other than $x$, is either in $\underline{s}$ or in $\underline{t}$ (because any cell of $\widetilde{\underline{x}}$ other than $x$ is in the source or target of $\widetilde{\underline{x}}$, and hence either in $\underline{s}$ or in $\underline{t}$) and finally the intersection of $\underline{s}$ and $\underline{t}$ in $\underline{p}$ is exactly this object $\underline{x'}-\{x'\} \subset X$. Indeed $\underline{s} \wedge \underline{t}$ are exactly the cells that are both in the source and target of $\widetilde{\underline{x}}$, and they corresponds exactly to cells of $\widetilde{\underline{x}}$ that are not $x$ or inner cells of the source and target of $x$ because of proposition \ref{Prop_unique_inner_cell}. One can then form the following ``cube diagram'':

\[\begin{tikzcd}[ampersand replacement=\&]
\pi_{n-1}^- \underline{p} \arrow[hook]{d}{\sim} \& \arrow[hook]{l} \partial \underline{x'} = \partial_{n-2} \underline{p} \arrow[hook]{d}{\sim} \arrow[hook]{r}  \& \pi_{n-1}^+ \underline{p} \arrow[hook]{d}{\sim} \\
\underline{s} \& \arrow[hook]{l} \underline{x'}-\{x'\}  \arrow[hook]{r}\& \underline{t}
\end{tikzcd}\]

Where $\partial_{n-2} \underline{p}$ denotes the shared $n-2$ dimensional boundary of $\pi_{n-1}^- \underline{p}$ and $\pi_{n-1}^- \underline{p}$. The leftmost and rightmost vertical map are equivalences because of $C_n$ (by taking any plex in the domain and applying $2$-out-of-$3$). The middle vertical map is an equivalence because of $D_{n-1}$. The pushout of the first line is $\partial \underline{p}$ and the pushout of the lower line is $\underline{p}-\{x \}$ as observed above. So the cube axioms shows that $\partial \underline{p} \rightarrow \underline{p}-\{x\}$ is an equivalence.

\bigskip

$D_n \text{ and } C_n \Rightarrow E_n :$

The domain of a generating anodyne map in $\Pb^{reg}_{n+1}$ is given by $X=\underline{p'} \coprod_{\partial \underline{p}} (\underline{p}- \{x \})$ where $\underline{p}$ and $\underline{p'}$ are $n$-polyplexes with same spherical boundary $\partial \underline{p} = \partial \underline{p'}$ and $x$ is any $n$-cell of $\underline{p}$. In particular $D_n$ and the stability of $W$ under pushout shows that $\underline{p'} \hookrightarrow X \in W$.

This is essentially enough to conclude: Any plex of $X$ that factor into $\underline{p'}$ induces an equivalence with $X$ because of $C_n$ and the stability under composition of $W$. This is sufficient to show the result for any zero dimensional plex as any such plex is connected by a series of $1$-polyplex to the $0$-source (or eventually $0$-target if we are in dimension $1$) of $X$ which is in $\underline{p'}$ and any other plex contains one of the zero dimensional plex so the $2$-out-of-$3$ property for $W$ allows to conclude.

\bigskip

$E_n \Rightarrow A_n$:

Note that $E_n \Rightarrow B_n$ as any plex $\underline{p}$ in $\Pb^{reg}_n$ is the domain of a generating anodyne map $\underline{p} \hookrightarrow \underline{p} \cup \{x,\theta\}$ in $\Pb^{reg}_{n+1}$ where $x$ is parallel to the top dimensional cell $p$ of $\underline{p}$ and $\theta$ is a $k+1$-cell between $x$ and $p$.

If $f :X \rightarrow \underline{p}$ is a generating anodyne map in $\Pb^{reg}_n$, then any plex $\underline{x} \rightarrow X$ is in $W$ because of $E_n$ and $\underline{x} \rightarrow \underline{p}$ is in $W$ because of $B_n$. So $2$-out-of-$3$ for $W$ implies that $f:X \rightarrow \underline{p} \in W$.

}
}

\block{\Prop{$R: \Pb^{reg} \rightarrow \widehat{\Delta}$ is a left Quillen functor.}

\Dem{It clearly send monomorphisms to monomorphisms. Now the set of monomorphisms in $\Pb^{reg}$ that are sent to equivalence in $\widehat{\Delta}$ contains iso, is stable under pushout, composition, $2$-out-of-$3$ and the cube lemma, and it contains all the map between plexes because all plexes are sent to contractible objects. Hence it satisfies all the condition of lemma \ref{Lemma_Key_weakequivalence} and so also contains the generating anodyne cofibrations.

}
}

\block{In the other direction, proposition \ref{Prop_street_orientals} shows that Street's Orientals form a downward closed, full category of the category of regular plex isomorphic to the semi-simplicial category. This identifies the category of semi-simplicial sets with the category of regular polygraphs whose cells are all orientals. One denotes by $\Ocal : \widehat{\Delta} \rightarrow \Pb^{reg}$  this inclusion of semi-simplicial sets into $\Pb^{reg}$. It commutes to all limits and all colimits, in particular it is a left adjoint functor.}

\block{\Lem{$\Ocal : \widehat{\Delta} \rightarrow \Pb^{reg}$ is a left Quillen functor.}

\Dem{It is a left adjoint functor, it clearly send monomorphisms to monomorphisms and the Horn inclusion are a special case of our generating anodyne maps.}
}

\block{We will conclude the proof by showing, using the following theorem that these two left adjoint functors are homotopically inverse from each other:}

\block{
\label{Th_self_equiv_of_directed_cat}
\Th{ Assume that:

\begin{itemize}

\item $\Dcal$ is a directed category. I.e. there is a height function from the set of objects of $\Dcal$ to the integer such that every non-identity arrow has a target of higher height than its source.

\item There is a weak model structure on the presheaf category $\widehat{\Dcal}$ where every object is cofibrant and the cofibrations are the monomorphisms\footnote{constructively we want the cofibration to be the levelwise complemented monomorphisms.}.

\item For every object $x \in \Dcal$, the map $x \rightarrow 1$ to the terminal object in $\widehat{\Dcal}$ is a weak equivalence for the model structure of $\widehat{\Dcal}$.

\item $L : \widehat{\Dcal} \rightarrow \widehat{\Dcal}$ is a left Quillen functor.

\item There is a specific object $* \in \Dcal$ such that every object of $\Dcal$ has a map from $*$.

\item $L(*)= *$

\end{itemize}

Then $Ho(L) : Ho(\widehat{\Dcal}^{cof}) \rightarrow Ho(\widehat{\Dcal}^{cof})$ is isomorphic to the identity. In particular, $L$ is a left Quillen equivalence.
}

The last two assumptions can be replaced by anything that will ensure that $L(d) \rightarrow 1$ is an equivalence for all $d \in \Dcal$.

\Dem{

Let $L : \widehat{\Dcal} \rightarrow \widehat{\Dcal}$ be such a functor. One denotes by $L^+$ the functor $L \coprod Id$. We will construct another left adjoint functor $W: \widehat{\Dcal} \rightarrow \widehat{\Dcal}$ endowed with a natural transformation: $\tau : L^+ \rightarrow W$ with the following property:

\begin{itemize}

\item[(1)] For every cofibration $x \hookrightarrow y$ in $\widehat{\Dcal}$ the map induced by the naturality square of $\tau$:

\[ L^+y \coprod_{L^+ x } W(x) \hookrightarrow W(y)\]

is a cofibration.

\item[(2)] The restriction of $\tau$ : $L \rightarrow W$ and $Id \rightarrow W$ are objectwise weak equivalences in $\widehat{\Dcal}$.

\end{itemize}

One easily see that the first condition implies that $W$ send cofibrations to cofibrations, and the second show that it send trivial cofibrations to trivial cofibrations. So in particular $W$ will also be a left Quillen functor. 

Finally the two restrictions $\tau_1 : L \rightarrow W$ and $\tau_2 : Id \rightarrow W$ being weak equivalences implies that all these functor are isomorphic when acting on the homotopy category of cofibrant objects, hence this proves the proposition.

So it remains to construct $W$. This will be done by induction on the height: we chose some height function $h:\Dcal \rightarrow \mathbb{N}$ (on the two cases of interest to us, the dimension is a good choice, the number of cells would also work), we let $\Dcal_n$ be the full subcategory of $\Dcal$ of object of height at most $n$. Note that $\widehat{\Dcal_n}$ is the full subcategory of object of $\widehat{\Dcal}$ that have no cells of height higher than $n$. By cell we mean element of $X(d)$ for some $d$, and the height of a cell is the height of $d$. 

As usual, if $d \in \Dcal$ one denotes by $\partial d$ the sub object of $d$ in $\widehat{\Dcal}$ of all non-identity morphisms to $d$. If $d$ is of height $n$, then $\partial d$ is an object of $\widehat{\Dcal_{n-1}}$.

The category $\Dcal_0$ (or more precisely, the first non-empty category $\Dcal_i$) is reduced to our object $*$. Indeed as $*$ admits maps to any other object his dimension is strictly smaller than the dimension of any other object, so we can decide to call it $0$. On $*$ one defines $W$ by sending the object $*$ to an interval object $I *$ for $*$ and it clearly satisfies the conditions.

We now assume that $W$ has been constructed on $\widetilde{\Dcal_{n-1}}$. One first defines the values of $W$ on objects of $\Dcal$ of height $n$. Let $d$ be such an object.

We define $W(d)$ as follow one takes a factorization as a cofibration followed by a weak equivalence\footnote{One could says trivial fibrations. But only asking weak equivalences will often allow to have a construction that only involve finite objects at each stage. Although this plays no role in the present paper.}:

\[ W(\partial d) \coprod_{L^+(\partial d)} L^+(d) \hookrightarrow W(d) \overset{\sim}{\rightarrow} 1 \]

One observes that $\tau_1: L(d) \rightarrow W(d)$ and $\tau_2: d \rightarrow W(d)$ are both weak equivalences because all these objects are equivalent to the terminal object. Indeed, $L(d) \rightarrow 1$ is a weak equivalence because $*=L(*) \hookrightarrow L(d)$ is a weak equivalence as the image of a trivial cofibration between cofibrant object by a left Quillen functor.

The functoriality $W(d') \rightarrow W(d)$ for any other object $d'$ is given by $W(d')= \rightarrow W(\partial d) \rightarrow W(d)$, and there is other arrow on which to define the functoriality. And the natural transformation $L^+(d) \rightarrow W(d)$ is the obvious coproduct inclusion.

This gives us a functor $W: \Dcal_n \rightarrow \widehat{\Dcal}$ with a natural transformation $L^+_{| \Dcal_n} \rightarrow W$ which admits a left Kan extension $\widehat{\Dcal_n} \rightarrow \Dcal$ also denoted $W$ with a natural transformation also denoted $\tau$ from $L^+$, which extend the previously constructed $W$. This new $W$ can be described as follow: in order to compute $W(X)$ one first compute $W(X^{n-1})$ where $X^{n-1}$ is the restriction of $X$ to $\Dcal_{n-1}$ and then for each cell of $X$ of height $n$, one take a pushout: 

\[ W(x) = W(X-\{x\}) \coprod_{W(\partial d)} W(d) \]

One will prove that this extension of $W$ has the desired property when applied to objects with a finite number of $n$-cells, by induction on the number of their $n$-cells, it then easily passes to directed limits to get the result for object with an infinite number of cells.

\begin{itemize}

\item If a natural transformation satisfies the property of condition $(1)$ for morphisms $f$ and $g$ in $\widehat{\Dcal}$ then one can check it satisfies it for their composite as well. So it is enough to check it for inclusion $x \hookrightarrow y$ that only add one cell, i.e. inclusion $Z=X-\{x\} \hookrightarrow X$, i.e. that the morphism:

\[ L^+(X) \coprod_{L^+(Z)} W(Z) \rightarrow W(X) \]

As $X= Z \coprod_{\partial d} d$, one has $L^+(X) = L^+(Z) \coprod_{L^+ \partial d} L^+ d$, hence the map above can be rewritten as:

\[ L^+(d) \coprod_{L^+(\partial d )} W(Z) \rightarrow W(X) \]

But one also has $W(X) = W(d) \coprod_{W(\partial d)} W(Z)$

So in the end, the map is the pushout along $W(\partial d) \rightarrow W(Z)$ of the map:

\[ W(\partial d) \coprod_{L^+(\partial d)} L^+(d) \hookrightarrow W(d)\]
which by definition of $W(d)$ is a cofibration.

\item We now prove that $\tau_1 : L(X) \rightarrow W(X)$ is an acyclic cofibration. One assume it is one for $Z=X-\{x\}$ and we prove it is one for $X$.

One has a cube diagram:

\[\begin{tikzcd}[ampersand replacement=\&]
L(Z) \arrow{d}{\sim} \& \arrow[hook]{l} L(\partial d)\arrow{d}{\sim} \arrow[hook]{r} \& L(d) \arrow{d}{\sim} \\
W(Z) \& \arrow[hook]{l} W(\partial d) \arrow[hook]{r} \& W(d)
\end{tikzcd}\]

The first vertical map is an equivalence because of our induction hypothesis. The second because it is in dimension $<n$ and the last one has been observed to be an equivalence before.

The exact same proof applies to $\tau_2$ and this conclude the proof.

\end{itemize}

}}

\block{\label{Th_Spaces_and_regular_equiv}\Th{Both functors $\Ocal \circ R$ and $R \circ \Ocal$ are homotopically equivalent to the identity functor respectively on $\Pb^{reg}$ and $\widehat{\Delta}$. In particular they are Quillen equivalences inverse of each other.}

\Dem{$R \circ \Ocal$ and the category $\widehat{D}$ immediately satisfies the assumption of our theorem \ref{Th_self_equiv_of_directed_cat}. For the case of $\Ocal \circ R$ the last remaining non-trivial condition is the fact that every map from $d \rightarrow 1$ is a weak equivalence. Our lemma \ref{Lemma_Key_weakequivalence} show that in $\Pb^{reg}$ any morphism between plexes is a weak equivalence, so it is enough to show that the map $* \rightarrow 1$ is a weak equivalence.

This will be proved using the cone construction of subsection \ref{subsec:joinsAndCones}.

Indeed, one easily see by induction on cells of $X$ that for any object $X$, the natural ``initial point'' map $* \rightarrow C X$ is a weak equivalence. Indeed it is an isomorphisms when $X$ is empty, and every time one add one $a$ cell to $X$, $CX$ gets two new cells, $a$ and $Ta$, and $a$ appears exactly once in the target of $Ta$ and not in its source (see the formula in \ref{discus:descriptionOfCQ}), so that the map $C(X-\{x\}) \rightarrow CX$ is a trivial cofibration.

Applying this to the terminal object itself gives us a map:   $\epsilon: 1 \rightarrow C1$ which in the homotopy category produce a morphisms $1 \rightarrow *$. The object $1$ is still terminal in the homotopy category, so the composite $1 \rightarrow * \rightarrow 1$ is always the identity in the homotopy category. The composite $* \rightarrow 1 \rightarrow *$ on the other hand corresponds to the morphisms $* \rightarrow C 1$ which send $*$ to the unique $0$-cell of $1$, and $C1$ being identified to $*$ along its other $0$-cell corresponding to its initial point. But these two points are canonically homotopic, by the $1$-cell $T *$, hence one indeed have that $*$ is isomorphic to $1$ in the homotopy category.

}
}

\subsection{Further conjectures}
\label{subsec:Conjecutres}

\block{\label{Conj:SimpsonGeneral}We can also state a more precise form of the general version of C.Simpson's conjecture:

\Conj{(The General Simpson conjecture for $\infty$-groupoid) The inclusion of $\widehat{\Delta}$ into $\Pb^+$ induced by Street's orientals is a left Quillen equivalence, or equivalently, the inclusion of $\Pb^{reg} \rightarrow \Pb^+$ is a left Quillen equivalence.}

As we already know that $\Pb^+$ is Quillen equivalent to the model category of non-unital $\infty$-category, this probably constitute the strongest version of C.Simpson conjecture's that we can hope for $\infty$-groupoids. It corresponds to the strictification of absolutely all composition operations that are not directly identities. Our result of \cite{henry2017positive} shows that there is no ``Eckmann-Hilton'' type obstruction for this result to hold, but a proof of this conjecture would require new ideas and typically a more clever ``geometric realization functors'' for positive polygraphs as we already suggested at the end of \cite{henry2017positive}.

The statement above is essentially equivalent to the conjecture we made at the end of \cite{henry2017positive}, but is now considerably more precise.

Note that non-unital $\infty$-category are algebra for a Batanin operad (in the sense of \cite{batanin1998monoidal} or \cite{leinster2004higher}), in fact even a sub-terminal Batanin operad. It seems that one can also encode a notion of weak unit in such an operads and we believe that a proof of the conjecture above would also provide a first example of a Batanin operad for which one can prove the homotopy hypothesis. A non-cellular operad though. We have no opinion on whether this is a sign that the conjecture above is very hard or a real hope to make some progress on the homotopy hypothesis.

}

\block{Finally there is the question of extending this to the $\infty$-categorical version of the Simpson conjecture, i.e. that (informally) any weak $\infty$-category\footnote{Here $\infty$ means $(\infty,n)$ for $n$ arbitrary large, or even ``$(\infty,\infty)$'' if we are optimistic.} is equivalent to one where all composition operations are strict, except units. Here again, there will be (at least) two version of the conjecture to consider: a ``general'' and a ``regular'' one, that strictify different set of operations.

 There is one additional difficulty of our approach when trying to move to this $\infty$-categorical context: for $\infty$-groupoids, unit could be characterized as idempotent. But for general $\infty$-categories unit should be characterized as ``invertible idempotent'', but one cannot defines invertible without having a notion of unit, so it is not going to be possible to have similarly a model structure on regular or non-unital $\infty$-categories such that the Fibrant object are the ones with weak units. We see two possible approach for this problem:

\begin{itemize}

\item One could add some operations attaching to each $n$-arrow $f$ a $n+1$-arrow $I_f:f \rightarrow f$ not satisfying any equation, and ask that in a fibrant objects these $I_f$ behave as weak units.

\item One could work with ``stratified'' or ``marked'' objects, similarly to J.Lurie model structure on Marked simplicial set from \cite{lurie2009higher}, or D.Verity model structure on stratified simplicial sets whose fibrant objects are the weak complicial sets (see \cite{verity2008weak}). I.e. one would work with marked/stratified version of polygraphs and $\infty$-categories. Weak units could then be characterized as marked idempotent and fibrant objects should have weak units and have all their marked arrow satisfying appropriate divisibility conditions.

\end{itemize}

We have a net preference for the second solution which some more flexible and simpler to handle, but both are probably equally good.

}

\block{We believe that the extension to $\infty$-categorical context and the passage from the regular to the general Simpson conjecture are completely orthogonal difficulties. We do not have many ideas serious idea at the moment on how to attack the general Simpson conjecture, but we believe the following approach could lead to a proof of the ``regular $\infty$-categorical'' Simpson conjecture (although one could start with the maybe easier $(\infty,1)$-categorical version:

\begin{enumerate}

\item One defines ``marked'' or ``stratified'' regular polygraphs and regular $\infty$-categories, as regular polygraphs and regular $\infty$-categories endowed with a set of cell/arrow that are marked (and morphisms should send marked cell/arrow to marked ones). Marked cells are thought of as isomorphisms or equivalence.

\item One extend the Gray tensor product to these marked objects. The goal is two have a tensor product that will represent the appropriate $\infty$-categorical product, so it should be the ``pseudo-Gray'' tensor product and not the lax one. i.e. all cell $f \otimes g$ for $f$ and $g$ of dimension $\geqslant 1$ should be marked.

\item One needs to come up with a good notion of ``generating anodyne maps'' for marked regular polygraphs that would generalize both the one used by D.Verity in \cite{verity2008weak} for stratified simplicial sets and the one of the present paper. One of the criterion is that they should satisfy the pushout-product axiom against the cofibrations (which are still the monomorphisms) for this extended Gray tensor product.

\item We have already shown in \cite{henry2018weakmodel} that there is a weak model structure on stratified semi-simplicial sets which is equivalent to Verity model structure on stratified simplicial sets.

\item One shows that the equivalences we constructed in the present paper can be lifted to equivalences between the marked version of the three model structures involved. 

\end{enumerate}

Clearly points $3.$ and $5.$ are the most vague and difficult ones.

}

\end{document}